# SMARANDACHE NON-ASSOCIATIVE RINGS

**W. B. Vasantha Kandasamy**

Department of Mathematics
Indian Institute of Technology, Madras
Chennai – 600036, India
e-mail: *vasantha@iitm.ac.in*
web: *http://mat.iitm.ac.in/~wbv*

**Definition:**

**Generally, in any human field, a *Smarandache Structure* on a set A means a weak structure W on A such that there exists a proper subset B ⊂ A which is embedded with a stronger structure S.**

**These types of structures occur in our everyday's life, that's why we study them in this book.**

**Thus, as a particular case:**

**A *Non-associative ring* is a non-empty set R together with two binary operations '+' and '.' such that (R, +) is an additive abelian group and (R, .) is a groupoid. For all a, b, c ∈ R we have (a + b) . c = a . c + b . c and c . (a + b) = c . a + c . b.**

**A *Smarandache non-associative ring* is a non-associative ring (R, +, .) which has a proper subset P ⊂ R, that is an associative ring (with respect to the same binary operations on R).**

# CONTENTS

# PREFACE

An associative ring is just realized or built using reals or complex; finite or infinite by defining two binary operations on it. But on the contrary when we want to define or study or even introduce a non-associative ring we need two separate algebraic structures say a commutative ring with 1 (or a field) together with a loop or a groupoid or a vector space or a linear algebra. The two non-associative well-known algebras viz. Lie algebras and Jordan algebras are mainly built using a vector space over a field satisfying special identities called the Jacobi identity and Jordan identity respectively. Study of these algebras started as early as 1940s. Hence the study of non-associative algebras or even non-associative rings boils down to the study of properties of vector spaces or linear algebras over fields.

But study of non-associative algebras using loops, that is loop algebras (or loop rings) was at its peak only in the 1980s. But till date there is no separate book on loop rings. So in this book we have given more importance to the study of loop rings and above all Smarandache loop rings. Further in my opinion the deeper study of loop rings will certainly branch off into at least five special type of algebras like Moufang algebras, Bol algebras, Bruck algebras, Alternative algebras and P-algebras. The author does not deny there are several research papers on these five special algebras but each of these algebras can be developed like Lie algebras and Jordan algebras.

Another generalized class of non-associative algebras are from groupoid rings. The class of loop rings are strictly contained in the class of groupoid rings. Several important study in the direction of Smarandache notions has been done in this book. Groupoid rings are a very new concept, till date only very few research papers exist in the field of Groupoid rings. The notion of Smarandache Lie algebras and Smarandache Jordan algebras are introduced. There are several informative and exhaustive books on Lie algebras and Jordan algebras. Further as the study of these two algebras heavily exploits the concepts of vector spaces and matrix theory we do not venture to approach Smarandache Lie algebras or Smarandache Jordan algebras in this direction. But in this book we approach Smarandache Lie algebras and Smarandache Jordan algebras as a study in par with associative ringsfor such a study has not yet been approached by any researcher in this way. This book is distinct and different from other non-associative ring theory books as the properties of associative rings are incorporated and studied even in Lie rings and Jordan rings, an approach which is not traditional.

This book has seven chapters. First chapter is not only introductory but also devotes an entire section to introduce the basic concepts in Smarandache vector spaces. Chapter two studies about Smarandache loop rings and in chapter three, the Smarandache groupoid rings is introduced and analysed. Chapter four is devoted to the introduction of several new concepts in Lie algebras and Smarandache Lie algebras. Chapter five defines Smarandache Jordan algebras and gives several new notions so far not studied about in Jordan algebras and Smarandache Jordan algebras. Overall, we have introduced nearly 160 Smarandache notions relating to non-associative rings. The sixth chapter gives suggestions for future study. The last chapter gives about 150 problems which will help any researcher.

I first deeply acknowledge my indebtness to Dr. Minh Perez of American Research Press, whose constant mental support and encouragement made me write this book.

I thank my husband Dr. Kandasamy and my daughters Meena and Kama, who tirelessly worked for very long hours, typing and proof-reading this book.

I dedicate this book to Mahatma Jotirao Phule (1827-1890), a veteran social revolutionary whose progressive ideas and activism – encompassing the triple concerns of caste, class and gender – heralded a definitely different future for India's oppressed masses. At the age of twenty, Phule established the first and only school for the girls of untouchable and oppressed castes, at a period of time in history where women and the oppressed castes were denied all education. His greatness remains in his challenging the ruthless gender-based and caste-based domination over education, a domination which still continues and requires our intensive struggle.

**Chapter 1**

# BASIC CONCEPTS

This chapter on basic concepts has three sections. First section mainly describes the specific notions of vector space and bilinear forms. In section two we develop the definition of Smarandache K-vectorial space given in [43]. In fact we develop several important properties in this direction, which has so far remained as an undeveloped concept. In the final section we just recall the definition of semigroups, groupoids, S-semigroups, S-groupoids and loops.

## 1.1 Basics of Vector space and bilinear forms

In this section we just recall the concept of vector space, bilinear form, symmetric and skew symmetric forms. As we expect the reader to be well versed with concepts in algebra we do not recall the definition of fields or rings.

**DEFINITION 1.1.1:** *A vector space or a linear space consists of the following:*

1. *a field K of scalars.*
2. *a set V of elements called vectors.*
3. *a rule or operation called vector addition which associates with each pair of vectors $\alpha, \beta$ in V a vector $\alpha + \beta$ in V, called the sum of $\alpha$ and $\beta$ in such a way that*

    a. *addition is commutative $\alpha + \beta = \beta + \alpha$.*
    b. *addition is associative that is $\alpha + (\beta + \gamma) = (\alpha + \beta) + \gamma$.*
    c. *There is a unique vector 0 in V, called the zero vector such that $\alpha + 0 = \alpha$ for all $\alpha \in V$.*
    d. *For each vector $\alpha$ in V there is a unique vector $-\alpha$ in V such that $\alpha + (-\alpha) = 0$.*

4. *An operation called scalar multiplication, which associates with each scalar c in K and a vector $\alpha$ in V a vector $c\alpha$ in V called the product of c and $\alpha$ in such a way That*

    a. *$1.\alpha = \alpha$ for all $\alpha \in V$.*
    b. *$(a_1 a_2)\alpha = a_1(a_2 \alpha)$.*
    c. *$a(\alpha + \beta) = a\alpha + a\beta$.*
    d. *$(a_1 + a_2)\alpha = a_1\alpha + a_2\alpha$, $a_1, a_2, a \in K$ and $\alpha, \beta \in V$.*

Thus we see a vector space unlike the abstract structures rings or fields is a composite structure consisting of two algebraic structures called fields and vectors, which forms an additive abelian group. Here it is pertinent to mention that 'vector' does not imply it is a vector quantity, only to mention that it is different from the elements of the field which are termed as scalars.

**DEFINITION 1.1.2:** *A subspace of the vector space V is a subset W of V which is itself a vector space over F with the operations of vector addition and scalar multiplication on V.*

The reader is advised to learn the concept of basis, linear transformation, linear operator etc. We say a vector space is finite dimensional if there exists a finite set of linearly independent elements $(v_1, v_2, \ldots, v_n)$ in V which generate V, otherwise it is infinite dimensional that is there does not exist a finite set of linearly independent elements in V which can span V.

The following definitions will be used in the later chapters, hence we recall them.

**DEFINITION 1.1.3:** *A bilinear form on a vector space V (real or complex) is a mapping $[x, y] \to (x, y)$ of $V \times V$ into the set of (real or complex) numbers such that*

$$(ax + bx', y) = a(x, y) + b(x', y) \text{ and}$$
$$(x, ay + by') = a(x, y) + b(x, y')$$

*for all vectors $x, x', y, y'$ in V and scalars a, b in the field.*

**DEFINITION 1.1.4:** *A bilinear form is skew symmetric if $(x, y) = -(y, x)$ for all vectors x, y. A Hermitian form on a complex vector space V is a mapping $[x, y] \to (x, y)$ of $V \times V$ into the set of complex numbers such that $(x, y) = \overline{(y, x)}$ and $(ax + bx', y) = a(x, y) + b(x', y)$ so that $(x, ay + by') = \bar{a}(x, y) + \bar{b}(x, y)$ for all vectors $x, x', y, y'$ and scalars a, b. A bilinear form $(x, y)$ is called non-singular if for every $x_0 \in V$ $(x_0 \neq 0)$ the linear form $(x, x_0)$ is not identically 0. Otherwise the bilinear form is said to be singular. This form is said to be positive definite if $(x, x) > 0$ for all $x \neq 0$; $x \in V$.*

**DEFINITION 1.1.5:** *A subset B of an associative algebra A over a field F is called weakly closed if for every ordered pair (a, b), $a, b \in B$, there is defined an element $\tau(a, b) \in F$ such that $ab + \tau(a, b) ba \in B$. We assume the mapping $(a, b) \to \tau(a,b)$ is fixed and write $a \times b = ab + \tau(a, b) ba$.*

*A subset $\Pi$ of B is called a subsystem if $c \times d \in \Pi$ for every $c, d \in \Pi$, and $\Pi$ is a left ideal (ideal) if $a \times c \in \Pi$ ($a \times c$ and $c \times a \in \Pi$) for every $a \in B$ and $c \in \Pi$. GL (V) or gl (v) denotes the group of all non-singular linear transformations of V, this is a linear algebraic group.*

## 1.2 Smarandache vector spaces

In this section we introduce the concept of Smarandache K-vectoroial space as given in [43]. We further illustrate with examples. As no one has worked on the Smarandache properties of Smarandache K-vectorial spaces here we have developed some of the concepts, which will be very useful in the study of linear structure and Lie algebra. The notions of Smarandache basis and Smarandache transformations are defined.

**DEFINITION [43]:** *The Smarandache K-vectorial space (S-K-vectorial space) is defined to be K-vectorial space (A, +, .) such that a proper subset of A is a K-algebra (with respect to the same induced operation and another '×' operation internal on A where K is a commutative field).*

By a proper subset we understand a set included in A, different from the empty set, from the unit element if any, and from A.

*Example 1.2.1:* Let $V_1 = F_{m \times n}$ be a vector space over Q and $V_2 = Q[x]$ be a vector space over Q. Then $V = V_1 \times V_2$, the direct product of the vector spaces over Q. We see V is a S-Q-vectorial space. For $\{0\} \times V_2 \subset V$ is a proper subset of V is a Q-algebra with respect to the same induced operation and the '×' in $V_2$ being the polynomial multiplication. Here K = Q.

*Example 1.2.2:* $V = V_1 \times V_2 \times V_3$ where $V_1 = Q \times Q$ is a vector space over Q, $V_2 = Q[x]$ is a vector space over Q and $V_3$ = {set of all odd degree polynomials with coefficients from Q}. V is a S-Q-vectorial space and $\overline{V_1} = V_1 \times \{0\} \times \{0\}$ is a proper subset of V which is a K-algebra.

Similarly $\overline{V_2} = \{0\} \times V_2 \times \{0\}$ is a proper subset which is a Q-algebra over Q. Thus by this example we illustrate that there can be more than one proper subset which is a K-algebra, for K = Q in this example. Now we proceed on to define a Smarandache K-vectorial subspace.

**DEFINITION 1.2.1:** *Let A be a K-vectorial space. A proper subset X of A is said to be a Smarandache K-vectorial subspace (S-K-vectorial subspace) of A if X itself is a S-K-vectorial space.*

**THEOREM 1.2.1:** *Let A be a K-vectorial space. If A has a S-K-vectorial subspace then A is S-K-vectorial space.*

*Proof:* Follows from the very definition. Thus we need not say that A is a S-K-vectorial space for if A has a S-K-vectorial subspace then A itself is a S-K-vectorial space.

The study of Smarandache basis is really a new concept.

**DEFINITION 1.2.2:** *Let V be a finite dimensional vector space over a field K. Let B = {$v_1, v_2, …, v_n$} be a basis of V. We say B is a Smarandache basis (S-basis) of V if B has a proper subset say A, $A \subset B$ and $A \neq \phi$; $A \neq B$ such that A generates a subspace which is linear algebra over K that is if W is the subspace generated by A then W must be a K-algebra with the same operation of V.*

Thus it is interesting to note that all basis of a vector space need not in general be a S-basis.

**THEOREM 1.2.2:** *Let V be a vector space over the field K. If B is a S-basis then B is a basis of V.*

*Proof*: Obvious by the very definition of these concepts.

**DEFINITION 1.2.3:** *Let V be a finite dimensional vector space over a field K. Let B = $\{v_1, v_2, \ldots, v_n\}$ be a basis of V. If every proper subset of B generates a linear algebra over K then we call B a Smarandache strong basis (S-strong basis) for V.*

**DEFINITION 1.2.4:** *Let L be any vector space over the field K. We say L is a Smarandache finite dimensional vector space (S-finite dimensional vector space) of K if every S-basis has only finite number of elements in it. It is interesting to note that if L is a finite dimensional vector space then L is a S-finite dimensional space provided L has a S-basis.*

It can also happen that L need not be a finite dimensional space still L can be a S-finite dimensional space.

**THEOREM 1.2.3:** *Let V be a vector space over the field K. If A = $\{v_1, v_2, \ldots, v_n\}$ is a S-strong basis of V then A is a S-basis of V.*

*Proof*: Easily follows from the very definition of these concepts.

**THEOREM 1.2.4:** *Let V be a vector space over the field K. If A = $\{v_1, v_2, \ldots, v_n\}$ is a S-basis of V. A need not in general be a S-strong basis of V.*

*Proof*: By an example. Let V = Q[x] be the set of all polynomials of degree less than or equal to 10. V is a vector space over Q. Clearly A = $\{1, x, x^2, \ldots, x^{10}\}$ is a basis of V. In fact A is a S-basis of V for take B = $\{1, x^2, x^4, x^6, x^8, x^{10}\}$. Clearly B generates a linear algebra. But all subsets of A do not form a S-basis of V, so A is not a S-strong basis of V but only a S-basis of V. Hence the claim.

Now we proceed on to define Smarandache eigen values and Smarandache eigen vectors of a vector space.

**DEFINITION 1.2.5:** *Let V be a vector space over the field F and let T be a linear operator from V to V. T is said to be a Smarandache linear operator (S-linear operator) on V if V has a S-basis, which is mapped by T onto another S-basis of V.*

**DEFINITION 1.2.6:** *Let T be a S-linear operator defined on the space V. A characteristic value c in F of T is said to be a Smarandache characteristic value (S-characteristic value) of T if the characteristic vector of T associated with c generate a subspace which is a linear algebra, that is the characteristic space associated with c is a linear algebra. So the eigen vector associated with the S-characteristic values will be called as Smarandache eigen vectors (S-eigen vectors) or Smarandache characteristic vectors (S-characteristic vectors).*

All other related concepts can be developed. But as this book is on non-associative ring we may not require all these properties we have just freshly defined the essential ones for making use of them in later chapters.

## 1.3 Basic definitions of other algebraic structures

In this section we just define the concept of semigroups, groupoids, S-semigroups, S-groupoid, S-loops and loops. For in chapters 2 and 3 the non-associative rings are introduced mainly using these algebraic structures loops and groupoids viz loop rings and groupoid rings. These rings are non-associative having several properties other than the properties enjoyed by Lie algebras and Jordan algebras.

**DEFINITION 1.3.1:** *Let S be a non-empty set on which is defined a binary operation '.' such that*

1. *$a.b \in S$ for all $a, b \in S$. (closure axiom)*
2. *$a. (b.c) = (a.b).c$ for all $a, b, c \in L$. (associative law)*

*(S, .) is called the semigroup. The semigroup in general may not have identity. If the semigroup S has an element $e \in S$, such that $s.e = e.s = s$ for all $s \in S$, then we call (S, .) a semigroup with unit or a monoid.*

**DEFINITION 1.3.2:** *Let G be a non-empty set on which is defined a binary operation '.' such that S is closed under '.' But '.' on S is not in general associative, then we say (S, .) is a groupoid. Thus a groupoid is a non-associative semigroup.*

A groupoid may or may not have the identity. If it has identity then we call the groupoid to be a groupoid with identity. Thus trivially all semigroups are groupoids and all groupoids in general need not be semigroups.

*Example 1.3.1:* Let $Z^+$ be the set of integers. Define an operation '.' on $Z^+$ by $a.b = na + mb$ where m and n are any pair of chosen numbers such that $(m, n) = 1$ ($m \neq 1$, $n \neq 1$). Clearly $(Z^+,.)$ is a groupoid which is not a semigroup but $Z^+$ under usual multiplication and addition are semigroups.

**DEFINITION 1.3.3:** *Let L be a non-empty set, we say L is a group if on L is defined a binary operation '$*$' such that*

1. *L is closed under '$*$' i.e. $a * b \in L$ for all $a, b \in L$.*
2. *$a * (b * c) = (a * b) * c$ for all $a, b, c \in L$.*
3. *There exist an element $e \in L$ such that $a * e = e * a = a$ for all $a \in L$.*
4. *For every $a \in L$ there exist a unique element $a^{-1} \in L$ such that $a * a^{-1} = a^{-1} * a = e$.*

A non-associative group is called a loop.

**DEFINITION 1.3.4:** *The Smarandache semigroup (S-semigroup) is defined to be a semigroup A such that a proper subset of A is a group.*

**DEFINITION 1.3.5:** *A Smarandache groupoid G (S-groupoid) is a groupoid, which has a proper subset $S \subset G$ which is a semigroup under the operations of G.*

**DEFINITION 1.3.6:** *A loop L is said to be a Smarandache loop (S-loop) if L has a proper subset, which is a group under the operations of L.*

**DEFINTION 1.3.7:** *Let R be a non empty set on which is defined two binary operations '+' and '.' ; (R, +, .) is a non-associative ring, if*

1. *(R, +) is an additive abelian group with 0*
2. *(R, .) is closed with respect to '.' and a . (b.c) ≠ (a.b) . c for some a, b, c ∈ R.*
3. *a. (b + c) =  a.b + a.c (The distributive law are true) for all a, b, c, ∈ R.*

*Let R be a non-assoiciative ring. R is said to be a Smarandache non-associative ring (SNA-ring) if S contains a proper subset P such that P is an associative ring under the operations of R.*

(We as a matter of convention put ab for a.b).

Loop rings are constructed analogous to group rings. Group rings are associative where as loop rings are non-associative rings. Similarly analogous to semigroup rings, groupoid rings are constructed replacing a semigroup by a groupoid. Groupoid rings are also non-associative rings. Thus the study of there two classes of non-associative rings is carried out in chapters 2 and 3 respectively.

**Chapter 2**

# LOOP RINGS AND SMARANDACHE LOOP RINGS

This chapter introduces the notion of Smarandache loop rings. It has five sections. In section one we just recall the definition of loops and Smarandache loops and give some of its basic properties which are very essential for the study of later sections in this chapter.

In section two we recall the properties of loop rings and introduce the basic definitions and notions of Smarandache loop rings. In fact several properties are defined and illustrated with examples. Section three defines several types of Smarandache elements. Around twenty theorems are proved about these special Smarandache elements. Smarandache substructures are analysed and defined in section four. Over twenty definitions and thirty-five theorems about these Smarandache substructures are given. In the final section we give a few general properties of loop rings and their Smarandache analogues.

## 2.1 Introduction of Loops and Smarandache loops

In this section we just recall the definition of loops and Smarandache loops, illustrate them with examples and enumerate its basic properties.

**DEFINITION 2.1.2**: *A non-empty set L is a said to form a loop, if in L is defined a binary operation called product and denoted by '.' such that for all a, b $\in$ L we have a.b $\in$ L, L contains an element e such that a.e = e.a = a for all a $\in$ L; e is called the identity element of L. For every ordered pair (a, b) $\in$ L $\times$ L there exists a unique pair (x, y) $\in$ L $\times$ L such that ax = b and ya = b. The binary operation '.' need not in general be associative for a loop.*

Thus all groups are loops and loops in general are not groups.

*Example 2.1.1*: Let (L, .) be a loop given by the following table:

| . | e | $a_1$ | $a_2$ | $a_3$ | $a_4$ | $a_5$ |
|---|---|---|---|---|---|---|
| e | e | $a_1$ | $a_2$ | $a_3$ | $a_4$ | $a_5$ |
| $a_1$ | $a_1$ | e | $a_4$ | $a_2$ | $a_5$ | $a_3$ |
| $a_2$ | $a_2$ | $a_4$ | e | $a_5$ | $a_3$ | $a_1$ |
| $a_3$ | $a_3$ | $a_2$ | $a_5$ | e | $a_1$ | $a_4$ |
| $a_4$ | $a_4$ | $a_5$ | $a_3$ | $a_1$ | e | $a_2$ |
| $a_5$ | $a_5$ | $a_3$ | $a_1$ | $a_4$ | $a_2$ | e |

**DEFINITION 2.1.2**: *Let (L, .) be a loop, we say L is a commutative loop if a.b = b.a for all a, b $\in$ L. Loop L is said to be of finite order if L has finite number of elements;*

*and we denote it by o(L) or |L|. If L has infinite number of elements then we say L is of infinite order.*

**DEFINITION 2.1.3**: *Let (L, .) be a loop and S a proper subset of L. We say (S, .) is a subloop if S itself is a loop under '.' A subloop S of L is said to be a normal subloop of L if*

1. $xS = Sx$.
2. $(Hx) y = H (xy)$.
3. $y (xH) = (yx) H$ for all $x, y \in L$.

*A loop L is said to be a simple loop if it does not contain any nontrivial normal subloops.*

**DEFINITION 2.1.4**: *If x and y are elements of a loop L, the commutator (x, y) is defined by $xy = (yx) (x, y)$.*

*The commutator subloop of a loop L (denoted by L') is the subloop generated by all of its commutators, that is $L' = \langle \{x \in L / x = (y, z) \text{ for some } y, z \in L\} \rangle$, where $L' \subseteq L$, $\langle L' \rangle$ denotes the subloop generated by L'.*

**DEFINITION 2.1.5**: *Let L be a loop, if x, y, z are elements of a loop L an associator(x, y, z) is defined by $(xy)z = (x(yz))$, (x, y, z) the associator subloop of a loop L (denoted by A (L)) is the subloop generated by all of its associators, that is, $\langle \{x \in L / x = (a, b, c) \text{ for some } a, b, c \in L\} \rangle$.*

**DEFINITION 2.1.6**: *A loop L is said to be a Moufang loop if it satisfies any one of the following identities*

$$(xy)(zx) = (x (yz)) x.$$
$$((xy) z) = x (y (xy)).$$
$$x (y (xz)) = ((xy) x) z$$

*for all $x, y, z \in L$.*

*These identities will be termed as Moufang identities.*

**DEFINITION 2.1.7**: *A loop L is called a Bruck loop if $(x (yx)) z = x (y (xz))$ and $(xy)^{-1} = x^{-1} y^{-1}$ for all $x, y, z \in L$.*

**DEFINITION 2.1.8**: *A loop (L, .) is called a Bol loop if $((xy) z) y = x ((yz) y)$ for all $x, y, z \in L$.*

**DEFINITION 2.1.9:** *A loop L is said to be right alternative if $(xy) y = x (yy)$ for all $x, y \in L$ and left alternative if $(xx)y = x(xy)$ for all $x, y \in L$. L is said to be an alternative loop, if it is both a right and a left alternative loop.*

**DEFINITION 2.1.10:** *A loop (L, .) is called a weak inverse property loop (WIP-loop) if $(xy)z = e$ imply $x(yz) = e$ for all $x, y, z \in L$.*

**DEFINITION 2.1.11:** *A loop L is said to be semi alternative if (x, y, z) = (y, z, x) for all x, y, z ∈ L, where (x, y, z) denotes the associator of elements x, y, z ∈ L.*

**DEFINITION 2.1.12**: *Let L be a loop. The left nucleus $N_\lambda$ = {a ∈ L / (a, x, y) = e for all x, y ∈ L} is a subloop of L. The middle nucleus $N_\mu$ = {a ∈ L / (x, a, y) = e for all x, y ∈ L} is a subloop of L. The right nucleus $N_p$ = {a ∈ L / (x, y, a) = e for all x, y ∈ L} is a subloop of L where e is the identity elemtent of L.*

*The nucleus N(L) of a loop L is the subloop given by $N(L) = N_\mu \cap N_\lambda \cap N_p$.*

**DEFINITION 2.1.13**: *The Moufang center C(L) is the set of elements of the loop L which commute with every element of L, that is C (L) = {x ∈ L / xy = yx for all y ∈ L}.*

*The center Z(L) of a loop L is the intersection of the nucleus and the Moufang center that is $Z(L) = C(L) \cap N(L)$.*

**DEFINITION 2.1.14**: *A map θ from a loop L to a loop L' is called a homomorphism if θ(ab) = θ(a) θ(b) for all a, b ∈ L.*

**DEFINITION 2.1.15:** *Let (L, .) be a loop. For α ∈ L define a right multiplication $R_a$ as the permutation of the loop (L, .) as follows. $R_a$: x → x.a we will call the set {$R_a$ / a ∈ L} the right regular representation of (L, .) or briefly, the representation of L.*

**DEFINITION 2.1.16**: *For any predetermined a, b in the loop L, a principal isotope (L, ∗) of the loop (L, .) is defined by x ∗ y = X.Y where X.a = x and b.Y = y. Let L be a loop. L is said to be G-loop if it is isomorphic to all of its principal isotopes.*

Now we define a new class of loops of even order built using $Z_n$.

**DEFINITION 2.1.17**: *Let $L_n(m)$ = {e, 1, 2, 3, …, n} be the set where n > 3, n is odd and m is a positive integer such that (m, n) = 1 and (m – 1, n) = 1 with m < n. Define on $L_n(m)$ a binary operation '.' as follows.*

$$e.i = i.e = i \text{ for all } i \in L_n(m)$$
$$i.i = i^2 = e \text{ for all } i \in L_n(m)$$
$$i.j = t \text{ where } t = (mj - (m-1) i) \pmod{n}$$

*for all i, j ∈ $L_n(m)$; i ≠ j, i ≠ e and j ≠ e.*

*Then ($L_n(m)$, .) is a loop.*

*Let $L_n$ denote the class of all loops $L_n(m)$ for a fixed n and various m's satisfying the conditions; m < n; (m, n) = 1 and (m – 1, n) = 1; that is $L_n$ = {$L_n(m)$ / n > 3, n odd, m < n, (m, n) = 1 and (m – 1, n) = 1}.*

***Example 2.1.2***: Let $L_7$ (3) = {e, 1, 2, 3, 4, 5, 6, 7} be a loop in $L_7$ given by the following table:

| . | e | 1 | 2 | 3 | 4 | 5 | 6 | 7 |
|---|---|---|---|---|---|---|---|---|
| e | e | 1 | 2 | 3 | 4 | 5 | 6 | 7 |
| 1 | 1 | e | 4 | 7 | 3 | 6 | 2 | 5 |
| 2 | 2 | 6 | e | 5 | 1 | 4 | 7 | 3 |
| 3 | 3 | 4 | 7 | e | 6 | 2 | 5 | 1 |
| 4 | 4 | 2 | 5 | 1 | e | 7 | 3 | 6 |
| 5 | 5 | 7 | 3 | 6 | 2 | e | 1 | 4 |
| 6 | 6 | 5 | 1 | 4 | 7 | 3 | e | 2 |
| 7 | 7 | 3 | 6 | 2 | 5 | 1 | 4 | e |

Clearly $L_7(3)$ is a WIP-loop and it is non-commutative and of order 8.

Now we proceed on to define the concept of Smarandache loops and give some of the Smarandache properties enjoyed by them as that alone will make this chapter self-contained.

**DEFINITION 2.1.18**: *A Smarandache loop (S-loop) is defined to be a loop such that a proper subset A of L is a subgroup (with respect to the same induced operations) that is $\phi \neq A \subset L$.*

*Example 2.1.3*: Let L be a loop given by the following table:

| . | e | 1 | 2 | 3 | 4 | 5 |
|---|---|---|---|---|---|---|
| e | e | 1 | 2 | 3 | 4 | 5 |
| 1 | 1 | e | 4 | 2 | 5 | 3 |
| 2 | 2 | 4 | e | 5 | 3 | 1 |
| 3 | 3 | 2 | 5 | e | 1 | 4 |
| 4 | 4 | 5 | 3 | 1 | e | 2 |
| 5 | 5 | 3 | 1 | 4 | 2 | e |

Clearly {e, 1}, {e, 2}, {e, 3}, {e, 4} and {e, 5} are subgroups of L. Hence $L_5(3)$ is a S-loop.

**THEOREM [72]**: *Every power associative loop is a S-loop.*

*Proof*: Obvious by the very definition of power associative loop as every element generates a group. Hence every power associative loop is a S-loop.

**THEOREM [72]**: *Every diassociative loop is a S-loop.*

*Proof*: Follows from the fact that every pair of elements in a loop generates a subgroup. Hence the claim.

**THEOREM [72]**: *Every loop $L_n(m)$ in $L_n$, n > 3, n an odd integer (m, n) = 1 and (n, m – 1) = 1 with m < n is a S-loop.*

*Proof*: Since $L_n(m) \in L_n$ are power associative loops as every element generates a group, we see each $L_n(m) \in L_n$ is a S-loop.

**DEFINITION [72]**: *Smarandache Bol loop (S-Bol loop) is defined to be a loop L such that a proper subloop (which is not a subgroup) A of L is a Bol loop (with respect to the operations of L). That is $0 \neq A \subset S$. Similarly we define Smarandache Bruck loop (S-Bruck loop), Smarandache Moufang loop (S-Moufang loop) and Smarandache right (left) alternative loop (S-right (left) alternative loop).*

In the definition we insists that A should be a subloop of L and not a subgroup of L. For every subgroup is a subloop but a subloop in general is not a subgroup. Further every subgroup will certainly be a Moufang loop, Bol loop, or a Bruck loop.

**THEOREM [72]**: *Every Bol loop is a S-Bol loop but every S-Bol loop is not a Bol loop.*

*Proof*: Left for the reader to prove.

**DEFINITION 2.1.19**: *Let L and L' be two S-loops, with A and A' its subgroups respectively. A map $\phi$ from L to L' is called a Smarandache loop homomorphism (S-loop homomorphism) if $\phi$ restricted to A is mapped onto the subgroup A' in L', that is $\phi: A \to A'$ is a group homomorphism.*

The concept of S-loop isomorphism and automorphism are defined in a similar way. For results about S-loops refer [72]. We conclude this section with the introduction of Smarandache subloop, Smarandache normal subloop and related concepts.

**DEFINITION 2.1.20**: *Let L be a loop a proper subset A of L is said to be a Smarandache subloop (S-subloop) of L if A is a subloop of L and A is itself a S-loop; that is A contains a proper subset B which is a group under the operations of L, i.e., we demand a S-subloop which is not a subgroup.*

**THEOREM 2.1.1**: *Let L be a loop, if L has a S-subloop then L is a S-loop.*

*Proof*: Straightforward by the definition.

**THEOREM 2.1.2**: *Let L be a S-loop every subloop of L need not be a S-subloop of L.*

*Proof*: By an example. Let $L_5(2) \in L_5$ be a loop. It is clearly a S-loop but $L_5(2)$ has no S-subloops.

**DEFINITIONS 2.1.21**: *Let L be a S-loop. If L has no subloops but only subgroups we call L a Smarandache subgroup loop (S-subgroup loop).*

**THEOREM 2.1.3**: *Let $L_n(m) \in L_n$ where n is a prime. Then the class of loops in $L_n$ are S-subgroup loop.*

*Proof*: Use simple number theoretic argument and use the fact n is a prime we see $L_n(m)$ has only subgroups and no subloops.

**DEFINITION 2.1.22**: *Let L be a loop. We say a nonempty subset A of L is a Smarandache normal subloop (S-normal subloop) if*

a. *A is itself a normal subloop of L.*

b. *A contains a proper subset B where B is a subgroup under the operations of L. If L has no S-normal subloop we say the loop L is Smarandache simple (S-simple)*

**THEOREM 2.1.4**: *Let L be a loop. If L has a S-normal subloop then L is a S-loop.*

*Proof*: Left for the reader as an exercise.

**DEFINITION 2.1.23**: *Let L be a finite S-loop. An element $a \in A$, $A \subset L$, A the subgroup of L is said to be Smarandache Cauchy element (S-Cauchy element) of L if $a^r = 1$ (r > 1), 1 is the identity element of L and r divides the order of L; otherwise a is not a S-Cauchy element of L.*

**DEFINITION 2.1.24**: *Let L be a finite S-loop. If every element in every subgroup is a S-Cauchy element then we say S is a Smarandache Cauchy loop (S-Cauchy loop).*

**THEOREM [72]**: *Every loop in the class of loops $L_n$ are S-Cauchy loops.*

*Proof*: Obvious by the very definitions.

**DEFINITION 2.1.25**: *Let L be a finite loop. If the order of every subgroup in L divides the order of the loop L then L is a Smarandache Lagranges loop (S-Lagranges loop).*

**DEFINITION 2.1.26:** *Let L be a finite loop. If there exist at least one subgroup in L whose order divides the order of L then L is Smarandache weakly Lagranges loop (S- weakly Lagranges loop).*

**DEFINITION 2.1.27**: *Let L be a loop. We say L is a Smarandache commutative loop (S-commutative loop) if L has a proper subset A such that A is a commutative group.*

*Example 2.1.4*: The loop $L_7(3) \in L_7$ is a S-commutative loop but clearly $L_7(3)$ is not a commutative loop.

**DEFINITION 2.1.28**: *If every subgroup of a loop L is commutative we call L a Smarandache strongly commutative loop (S-strongly commutative loop).*

The author is requested to obtain the relations between these concepts.

**THEOREM [72]**: *Let L be a power associative loop, L is then a S-commutative loop.*

*Proof*: Direct by the very definitions.

Thus we have a class of loops which are S-commutative.

**DEFINITION 2.1.29**: *Let L be a loop. L is a Smarandache cyclic loop (S-cyclic loop) if L contains atleast a proper subset A which is a cyclic group. We say L is*

*Smarandache strongly cyclic (S-strongly cyclic) if every subset which is a subgroup is a cyclic group.*

**THEOREM [72]**: *Let n be a prime, every loop in the class of loops $L_n$ are S-strongly cyclic loops.*

*Proof*: Follows from the fact when n is a prime, $L_n(m)$ for any m have only n subgroups each cyclic of order 2.

**THEOREM [72]**: *Let $L_n$ be a class of loops, n > 3, if $n = p_1^{\alpha_1} p_2^{\alpha_2} ... p_k^{\alpha_k}$, ($\alpha_i \geq 1$ for i = 1, 2, 3, …, k) then it contains exactly $F_n$ loops which are strictly non-commutative and they are*

> *S-strongly commutative loops*
> *S-strongly cyclic loops; where $F_n = \prod_{i=1}^{k}( p_i - 3 )p_i^{\alpha_i - 1}$*

*Proof*: Kindly refer [72].

**DEFINITION 2.1.30**: *Let L be a loop. If A (A a proper subset of L) is a S-subloop of L is such that A is a pseudo commutative loop then we say L is a Smarandache pseudo commutative loop (S-pseudo commutative loop) i.e., for every a, b ∈ A, we have x ∈ B such that a(xb) = b (xa) (or (bx) a) B a subgroup in A.*

*If a.x.b = b.x.a for all x ∈ A and a, b ∈ A we say L is a Smarandache strongly pseudo commutative loop (S-strongly pseudo commutative loop). Now we proceed on to define the concept of Smarandache commutator subloop (S-commutator subloop) of a loop L denoted by $L^s$. Let L be a loop; A ⊂ L be a S-subloop of L, $L^s = \langle \{x \in A / x = (y, x)$ for some $y, z \in A\}\rangle$ with the usual notation and if L has no S-subloops but L is a S-loop we replace A by L.*

**THEOREM [72]**: *Let L be a S-loop, which has no S-subloops then $L^1 = L^s$.*

*Proof*: Left for the reader to prove.

Similarly we have to define the concept of Smarandache associative subloop of a loop L.

**DEFINITION 2.1.31**: *Let L be a loop. We say L is a Smarandache associative loop (S-associative loop) if L has a S-subloop A such that, A contains a triple x, y, z (all three elements distinct and different from e, the identity element of A) such that x.(y.z) = (x.y).z. This triple is called the Smarandache associative triple (S-associative triple).*

*We say L is a Smarandache strongly associative loop (S-strongly associative loop) if in L every S-subloop has an S-associative triple.*

*Suppose if x, y ∈ A (A ⊂ L, L a loop, A a S-subgroup of L) we have (xy) x = x (yx) then we say the loop L is a Smarandache pair wise associative loop (S-pairwise associative loop).*

*The Smarandache associator subloop (S-associator subloop) of L is denoted by $L^A$, is the subloop generated by all the associators in A, where A is a S-subloop of L ($A \subset L$) i.e., $L^A = \langle \{x \in A \,/\, x = (a, b, c) \text{ for some } a, b, c \in L\} \rangle$. If L has no S-subgroup but L is a S-loop we replace A by L itself. Thus we have in that case the S-associator coincides with the associator of the loop L.*

**DEFINITION 2.1.32**: *Let L be a loop, A is a S-subloop in L, if we have an associative triple a, b, c in A such that (ax) (bc) = (ab) (xc) for some $x \in A$, we say L is a Smarandache pseudo associative loop (S-pseudo associative loop) if every associative triple in A is a pseudo associative triple of A for some x in A. If L is a S-loop having no S-subloops then we replace A by L. Suppose L is a S-pseudo associative loop, then PA $(L^s) = \langle \{t \in A \,/\, (ab)(tc) = (at)(bc) \text{ where } (ab)c = a(bc)$ for $a, b, c \in A\} \rangle$ denotes the Smarandache pseudo associator (S-pseudo associator).*

We define or recall the definition of these concepts as they will play a major role in the definition and study of properties of the Smarandache non-associative rings, which is very recent [73].

**PROBLEMS:**

1. Give an example of a S-loop of order 15.
2. Can a S-loop of order 5 be S-cyclic?
3. Given an example of a S-loop of order 21, which is S-commutative.
4. Find $L^s$ for the loop $L = L_9(5)$.
5. Does $L_{13}(5)$ have nontrivial $L^A$?

## 2.2 Properties of Loop Rings and Introduction to Smarandache loop rings

In this section we just recall the definition of loop rings and enumerate several of the properties enjoyed by loop rings and here give some of its basic Smarandache properties mainly using the well known class of loops $L_n$. Further we state that throughout this paper by a ring R we mean only a commutative ring with 1 or a field; it can be finite or infinite. $Z_n$ denotes the ring of integers modulo n, n a composite number, $Z_p$, the prime field of characteristic p, Z the set of integers which is a ring, in fact an integral domain, Q the field of rationals and R the field of reals. L will denote a loop non-necessarily abelian but never a group so by a loop we mean a non-associative group.

**DEFINITION 2.2.1**: *Let R be a commutative ring with identity 1 and let L be a loop. The loop ring of the loop L over the ring R denoted by RL consists of all finite formal sums of the form $\sum_i \alpha_i m_i$ (i – runs over a finite number) where $\alpha_i \in R$ and $m_i \in L$ satisfying the following conditions.*

a. $\sum_{i=1}^{n} \alpha_i m_i = \sum_{i=1}^{n} \beta_i m_i \Leftrightarrow \alpha_i = \beta_i$ for $i = 1, 2, 3, \ldots, n$

b. $\left[\sum_{i=1}^{n} \alpha_i m_i\right] + \left[\sum_{i=1}^{n} \beta_i m_i\right] = \sum_{i=1}^{n} (\alpha_i + \beta_i) m_i$

c. $\left(\sum_{i=1}^{n} \alpha_i m_i\right)\left(\sum_{j=i}^{n} \beta_i n_j\right) = \underset{m_k = m_i m_j}{\Sigma \gamma_k m_k}$ where $\lambda_k = \Sigma \alpha_i \beta_j$

d. $r_i m_i = m_i r_i$ for all $r_i \in R$ and $m_i \in L$.

e. $r \sum_{i=1}^{m} r_i m_i = \sum_{i=1}^{m} (rr_i) m_i$ for $r \in R$ and $\Sigma r_i m_i \in RL$.

RL is a non-associative ring with $0 \in R$ as its additive identity. Since $1 \in R$ we have $L = 1.L \subseteq RL$ and $Re = R \subseteq RL$, where $e$ is the identity of $L$.

We give one example of a loop ring.

***Example 2.2.1***: Let $Z_2 = \{0, 1\}$ be the prime field of characteristic two and L be the loop given by the following table:

| . | e | a | b | c | d |
|---|---|---|---|---|---|
| e | e | a | b | c | d |
| a | a | e | c | d | b |
| b | b | d | a | e | c |
| c | c | b | d | a | e |
| d | d | c | e | b | a |

$Z_2L$ is the loop ring, which is obviously a non-associative ring.

We give some of the basic definitions about loop rings.

**DEFINITION 2.2.2**: *A loop ring RL is said to be an alternative loop ring if $(xx)y = x(xy)$ and $x(yy) = (xy)y$ for all $x, y \in RL$. In a ring R let $x, y \in R$. A binary operation known as the circle operation is denoted by 'o' and is defined by $x \circ y = x + y - xy$.*

**DEFINITION 2.2.3**: *Let $x \in R$, R a ring. x is said to be right quasi regular (r.q.r) if there exist a $y \in R$ such that $x \circ y = 0$ and x is said to be left quasi regular (l.q.r) if there exists a $y' \in R$ such that $y' \circ x = 0$. An element is quasi regular (q.r) if it is both right and left quasi regular. y is known as the right quasi inverse (r.q.i) of x and y' is the left quasi inverse (l.q.r) of x. A right ideal or left ideal in R is said to be right quasi regular (l.q.r or q.r respectively) if each of its element is right quasi regular (l.q.r or q.r respectively).*

**DEFINITION 2.2.4**: *Let R be a ring. An element $x \in R$ is said to be a regular element if there exists a $y \in R$ such that $xyx = x$.*

Since loop rings are non-associative we need to define left regular element and right regular elements for loop rings.

**DEFINITION 2.2.5**: *Let RL be a loop ring of a finite loop L over the ring R. An element $x \in RL$ is said to be right regular (left regular) if there exists a $y \in RL$ ($y' \in RL$) such that $x(yx) = x$ or $((xy')x = x)$.*

*The Jacobson radical J(R) of a ring R is defined as follows:*

*$J(R) = \{a \in R \,/\, aR$ is right quasi regular ideal of $R\}$.*

*A ring R is said to be semi-simple if $J(R) = \{0\}$, where $J(R)$ is the Jacobson radical of R.*

**DEFINITION 2.2.6**: *Let RL be a loop ring of a loop L over the ring R. Let*

$$\alpha = \sum_{i=1}^{m} \alpha_i m \in RL.$$

*Support of $\alpha$ denoted by $\operatorname{supp} \alpha = \{m_i \,/\, \alpha_I \neq 0\}$. $|\operatorname{Supp} \alpha|$ denotes the number of elements in support of $\alpha$. The augmentation ideal of RL denoted by $W(RL) = \{\alpha = \Sigma\alpha_i m_i \in RL \,/\, \Sigma\alpha_i = 0\}$.*

Now we recall some of the basic properties starting from some of the properties enjoyed by the new class of loops over rings $Z_2$ or $Z_p$ or $Z$ or $Q$.

**THEOREM 2.2.1**: *Let L be a loop of order n and R be any commutative ring with 1 and of characteristic 0. RL be the loop ring of the loop L over R. If $\alpha = \Sigma\alpha_i m_i$ is right quasi regular, then $\Sigma \alpha_i \neq 1$ ($\alpha_i \in R$ and $m_i \in L$ for all $i = 1, 2, 3, \ldots, n$).*

*Proof*: It is a matter of routine once we take $L = \{m_1 = 1, m_2, \ldots, m_n\}$ to be the n elements of L and make use of the fact, for

$$\alpha = \sum_{i=1}^{n} \alpha_i m_i, \quad \beta = \sum_{j=1}^{n} \beta_j m_j$$

in RL; $\alpha + \beta - \alpha\beta = 0$. Equating the coefficients of the like terms and adding these coefficients we get $\Sigma\alpha_i + \Sigma\beta_j - (\Sigma\alpha_i)(\Sigma \beta_j) = 0$ or $\Sigma \alpha_i = \Sigma\beta_j (\Sigma\alpha_i - 1)$. If $\Sigma\alpha_i = 1$ then we get $\Sigma \alpha_i = 0$ a contradiction hence $\Sigma \alpha_i \neq 1$. When we replace the ring of integers R by Z we get if $\alpha = \Sigma\alpha_i m_i \in ZL$ is right quasi regular only in two cases when

$\Sigma \alpha_i = 0$ or $\Sigma \alpha_i = 2$. A similar argument will yield the result directly.

Finally in case of finite characteristic we arrive at a result only in case of the prime field of characteristic p; i.e., $Z_p$. Let $Z_pL$ be the loop ring, $\Sigma\alpha_i m_i$ is right quasi regular only when $\Sigma \alpha_i \not\equiv 1 \pmod{p}$.

**THEOREM 2.2.2**: *Let L be a loop and $Z_2L$ be the loop ring. ($Z_2 = (0, 1)$) the prime field of characteristic two). If $\alpha \in Z_2L$ is right quasi regular, then $|supp\ \alpha|$ is an even number.*

*Proof*: Left for the reader as an exercise.

**THEOREM 2.2.3**: *Let $Z_2$ be the prime field of characteristic two. L a commutative loop in which every element $m \in L$ is such that $m^2 = 1$, then $\alpha = \Sigma\ \alpha_i\ m_i \in Z_2\ L$ is right quasi regular if and only if $|supp\alpha|$ is an even number.*

*Proof*: Direct by using the fact $m^2 = 1$ for every $m \in L$ and L is a commutative loop.

Thus we have still modified results in case of field of prime characteristic p, p an odd prime.

*Example 2.2.2*: Let L be a loop given by the following table:

| . | e | a | b | c | d |
|---|---|---|---|---|---|
| e | e | a | b | c | d |
| a | a | e | c | d | b |
| b | b | d | a | e | c |
| c | c | b | d | a | e |
| d | d | c | e | b | a |

Z the ring of integers. ZL be the loop ring of a loop L over Z. $e + a \in ZL$ is quasi regular for $(e + a)\ o\ (e + a) = (e + a) + (e + a) - (e + a)(e + a) = 0$. Consider $x = a - d \in ZL$ with $y = e + a - c + d$ we get $(a - d)\ o\ (e + a - c + d) = a - d + (e + a - c + d) - (a - d)(e + a - c + d) = 0$ and

$(e + a - c + d)\ o\ (a - d) = e + a - c + d + a - d - (e + a - c + d)(a - d) = 0$. Thus $a - d$ is also quasi regular, now consider

$(e + a)\ 0\ (a - d) = e + a + a - d - (e + a)(a - d) = a + b$.

It is easily verified $a + b$ has no right quasi inverse and left quasi inverse in ZL. Thus we make the following important observation that in general the set of quasi regular elements of a ring need not form a closed set under the circle operation.

One of the nice results about Jacobson radical for loops rings is

**Result**: Let L be a finite loop. $Z_2$ the prime field of characteristic 2 and $Z_2L$ be the loop ring of the loop L over $Z_2$. Then $m_i + m_j \notin J(Z_2L)$ (where $|L| = n$, n is odd, $n \geq 7$; $Z_2 = [0, 1]$) for any $m_i, m_j \in L$.

**THEOREM 2.2.4**: *Let F be a field of characteristic zero and L a finite loop; FL the corresponding loop ring. If $\alpha = \Sigma\ \alpha_i\ m_i \in FL\ (\Sigma\alpha_i \neq 0)$ is a right (left) regular element with y as its right (left) relative inverse, then $\Sigma\ \beta_j\ = 1/\Sigma\ \alpha_i$.*

*Proof*: Let $\alpha$ be the right regular with y as its right relative inverse. Then x (yx) = x that is $\Sigma \alpha_i m_i ((\Sigma \beta_j m_j) (\Sigma \alpha_i m_i)) = \Sigma \alpha_i m_i$. Equating the coefficients of the like terms and adding all the equations we get $(\Sigma \alpha_i)((\Sigma \beta_j)(\Sigma \alpha_i)) = \Sigma \alpha_i$ or $\Sigma \beta_j = 1 / \Sigma \alpha_i$ (as $\Sigma \alpha_i \neq 0$).

**THEOREM 2.2.5**: *Let $Z_p L$ be the loop ring of a finite loop $L = \{m_1 = e, m_2, ..., m_n\}$ over the prime field $Z_p$. Then $x = k(e + m_2 + ... + m_n)$ ($0 \neq k \in Z_p$) is regular if and only if $p \nmid n$ where n is the order of the loop.*

*Proof*: Left for the reader as an exercise.

In view of this it can be easily proved that if p/n then $Z_p L$ is not regular.

**THEOREM 2.2.6**: *Let R be the field of reals and L a finite loop in which $m_i^2 = 1$ for every $m_i \in L$. Then the loop ring RL has a nontrivial idempotent $\alpha = \Sigma \alpha_i m_i$ only if $\alpha_1 < 1$.*

*Proof*: By simple number theoretic arguments we get the result as $\alpha_1^2 + \alpha_2^2 + ... + \alpha_n^2 = \alpha_1$ so $\alpha_1 < 1$.

**THEOREM 2.2.7**: *Let $Z_p L$ be the loop ring of a loop L over $Z_p$. Then $x = am_i + bm_j \in Z_p L$ ($a \neq 0$, $a \neq b$, $b \neq 0$, $m_i \neq m_j$ with p an odd prime) is an idempotent if and only if $m_i^2 = e$ and $m_j^2 = e$, $a = (p+1)/2$ and $b = (p-1)/2$.*

*Proof*: Simple calculations using number theoretical concepts, the proof is easily obtained. Hence left for the reader to prove.

We now give some relation between the quasi regular elements and idempotents, how the existence of one, results in the non-existence of the other which is given by the following theorem. The reader is expected to prove it.

**THEOREM 2.2.8**: *Let F be a field of characteristic two and let $L = \{m_1 = e, m_2, ..., m_n\}$ be the finite loop of order n. FL the loop ring of L over F. Then $\alpha = t(e + m_2 + ... + m_n)$ ($t \neq 0$) is a quasi regular element in FL if and only if $\alpha$ is not an idempotent element in L.*

Yet another result of some significance is.

**THEOREM 2.2.9**: *FL be the loop ring of a finite loop L in which there exists at least one $m_i$ in L such that $m_i^2 = e$ ($m_i \neq e$) and F be a field of characteristic zero. Then $\alpha = m(e + m_i) \in FL$ ($m \neq 0$) is quasi regular if and only if $\alpha$ is not an idempotent.*

*Proof*: Straightforward, hence left for the reader to prove.

Thus we can conclude if L is a loop such that for each $m_i \in L$ we have $m_i^2 = 1$ and Z be the ring of integers then in the loop ring ZL no non-zero idempotent can be right quasi regular. Further even in $Z_2 L$, L satisfying the condition $m_i^2 = 1$ for all $m_i$ in L and $|L| < \infty$ a non-zero idempotent element cannot be right quasi regular.

Now we give the relation between the J(FL) and W(FL) in a loop ring.

**THEOREM 2.2.10**: *Let L be a finite loop and let F be any field of characteristic zero or ring of integers. Then $J(FL) \subset W(FL)$.*

*Proof*: The reader is expected to prove.

Now a nontrivial question would be; can we ever have J(FL) = W(FL)? The answer is yes; given by the following theorem, which the reader is expected to prove.

**THEOREM 2.2.11**: *Let F be any field or ring of integers and L be a finite loop. Then W(FL) = J(FL) if and only if every element of W(FL) is right quasi regular.*

This is however illustrated by the following example.

*Example 2.2.3*: Let $Z_2 = \{0, 1\}$ be the prime field of characteristic two and L be a loop given by the following table:

| .     | e     | $x_2$ | $x_3$ | $x_4$ | $x_5$ | $x_6$ |
|-------|-------|-------|-------|-------|-------|-------|
| e     | e     | $x_2$ | $x_3$ | $x_4$ | $x_5$ | $x_6$ |
| $x_2$ | $x_2$ | e     | $x_5$ | $x_3$ | $x_6$ | $x_4$ |
| $x_3$ | $x_3$ | $x_5$ | e     | $x_6$ | $x_4$ | $x_2$ |
| $x_4$ | $x_4$ | $x_3$ | $x_6$ | e     | $x_2$ | $x_5$ |
| $x_5$ | $x_5$ | $x_6$ | $x_4$ | $x_2$ | e     | $x_3$ |
| $x_6$ | $x_6$ | $x_4$ | $x_2$ | $x_5$ | $x_3$ | e     |

It is easily verified $J(Z_2L) = W(Z_2L)$. If we replace $Z_2$ by $Z_3$ we see that $J(Z_3L) \subset W(Z_3L)$. Z be the ring of integers and let L be a finite loop, then $J(ZL) \subset W(ZL)$. Several such results can be obtained. Now we proceed on to define Smarandache loop ring we in fact will define two levels of Smarandache loop rings.

**DEFINITION 2.2.7**: *Let L be a loop and R any ring. The loop ring RL is a Smarandache loop ring (S-loop ring) if and only if L is a S-loop.*

Thus in the very definition itself we demand the loop to be a S-loop. Even if the loop ring RL is a SNA-ring still we will not accept the ring to be a S-loop ring.

A natural question would be what is a SNA-ring. We first define it in the following. We first define level I Smarandache NA rings.

**DEFINITION 2.2.8**: *Let (R, +, .) be a non-associative (NA) ring. If a proper subset P of R is an associative ring under the operation of '+' and '.' then we say R is a Smarandache NA ring (SNA-ring) of level I. That is if a ring is properly contained in a non-associative ring then we call the NA ring to be a Smarandache NA-ring (SNA-ring).*

*Example 2.2.4*: Let L be a loop given by the following table:

| . | e | $a_1$ | $a_2$ | $a_3$ | $a_4$ | $a_5$ | $a_6$ | $a_7$ |
|---|---|---|---|---|---|---|---|---|
| e | e | $a_1$ | $a_2$ | $a_3$ | $a_4$ | $a_5$ | $a_6$ | $a_7$ |
| $a_1$ | $a_1$ | e | $a_5$ | $a_2$ | $a_6$ | $a_3$ | $a_7$ | $a_4$ |
| $a_2$ | $a_2$ | $a_5$ | e | $a_6$ | $a_3$ | $a_7$ | $a_4$ | $a_1$ |
| $a_3$ | $a_3$ | $a_2$ | $a_6$ | e | $a_7$ | $a_4$ | $a_1$ | $a_5$ |
| $a_4$ | $a_4$ | $a_6$ | $a_3$ | $a_7$ | e | $a_1$ | $a_5$ | $a_2$ |
| $a_5$ | $a_5$ | $a_3$ | $a_7$ | $a_4$ | $a_1$ | e | $a_2$ | $a_6$ |
| $a_6$ | $a_6$ | $a_7$ | $a_4$ | $a_1$ | $a_5$ | $a_2$ | e | $a_3$ |
| $a_7$ | $a_7$ | $a_4$ | $a_1$ | $a_5$ | $a_2$ | $a_6$ | $a_3$ | e |

Z be the ring of integers, ZL be the loop ring, clearly ZL is a non-associative ring. Now 1.L = L ⊂ ZL; so ZL is a SNA-ring. A natural question is "Will every loop ring be a SNA-ring?" The answer is yes as every loop contains 1 and the ring R over which the loop ring is defined is an associative ring. Hence every loop ring is a SNA-ring of level II.

But is the loop ring a S-loop ring the answer is need not be. For if the loop L fails to be a S-loop, the loop ring will not be a S-loop ring even if it happens to be a SNA-ring.

**DEFINITION 2.2.9**: *Let R be a NA-ring. A non-empty subset S of R is said to be a SNA-subring of R if S contains a proper subset P such that P is an associative ring under the operations of R.*

**THEOREM [73]**: *Let R be a NA-ring. If R has a SNA-subring then R is a SNA-subring.*

*Proof*: Direct by the very definitions.

**THEOREM 2.2.12**: *Let (R, +, .) be a SNA-ring. Every subring of R need not in general be a SNA-subring of R.*

*Proof*: Let ZL be the loop ring given in example 2.2.4. Clearly if we take H = {e, $a_1$} then ZH is a subring of ZL which is not a SNA-subring of ZL. Hence the claim.

Now we proceed on to define the concept of Smarandache ideal in a SNA-ring.

**DEFINITION 2.2.10**: *Let R be a non-associative ring. A proper subset I of R is said to be a Smarandache NA right ideal (SNA-right ideal) of R if*

*I is a SNA-subring of R; i.e., J ⊂ I, J is a proper subset of I which is an associative subring under the operations of R.*

*For all i ∈ I and j ∈ J we have either ij or ji is in J. If I is simultaneously a SNA-right ideal and SNA-left ideal then we say I is a SNA-ideal of R.*

**THEOREM 2.2.13**: *Let R be a non-associative ring. If R has a SNA-ideal then R is a SNA-ring.*

*Proof*: Straightforward by the very definition of SNA-ideal. Another interesting result is that in general a SNA-ideal need not be an ideal.

In view of this we have the following.

**THEOREM 2.2.14**: *Let R be any non-associative ring. I be a SNA-ideal of R then I in general need not be an ideal of R.*

*Proof*: By an example, we can prove; the task of which is left for the reader.

**DEFINITION 2.2.11**: *Let R be a NA ring we say R is a Smarandache NA Moufang ring (SNA-Moufang ring) if R contains a subring S where S is a SNA-subring and for all x, y, z in S we have (x.y).(z.x) = (x.(y.z)).x that is the Moufang identity is true in S.*

It is to be noted that we do not demand the whole of R to satisfy Moufang identity. It is sufficient if a S-subring of R satisfies the identity, that is enough for the ring to be a SNA-Moufang ring.

**DEFINITION 2.2.12**: *Let R be a non-associative ring, R is said to be a Bol ring if R satisfies the Bol identity ((x∗y) ∗ z) ∗ y = x ∗ ((y ∗ z) ∗ y) for all x, y, z ∈ R.*

We now define the condition for it to be a SNA-Bol ring.

**DEFINITION 2.2.13**: *Let R be a non-associative ring, R is said to be SNA-Bol ring if R contains a subring S such that S is a SNA-subring of R and we have the Bol identity to be true for all x, y, z ∈ S ⊂ R.*

In view of this we have the following theorem:

**THEOREM 2.2.15**: *Let R be a non-associative ring which is a Bol ring. Then R is a SNA-Bol ring.*

*Proof*: Direct by the very definitions.

**THEOREM 2.2.16**: *Let R be a non-associative ring which is a SNA-Bol ring then R need not be a Bol ring.*

*Proof*: The reader is requested to construct an example to prove this theorem.

**DEFINITION 2.2.14**: *Let R be any non-associative ring. R is said to be a right alternative ring if (xy) y = x (yy) for all x, y ∈ R similarly left alternative if (xx) y = x (xy) for all x, y ∈ R. Finally R is an alternative ring if it is simultaneously both right alternative and left alternative.*

*Example 2.2.5*: Let Z be the ring of integers. L be a loop given by the following table:

| . | e | $g_1$ | $g_2$ | $g_3$ | $g_4$ | $g_5$ |
|---|---|---|---|---|---|---|
| e | e | $g_1$ | $g_2$ | $g_3$ | $g_4$ | $g_5$ |
| $g_1$ | $g_1$ | e | $g_3$ | $g_5$ | $g_2$ | $g_4$ |
| $g_2$ | $g_2$ | $g_5$ | e | $g_4$ | $g_1$ | $g_3$ |
| $g_3$ | $g_3$ | $g_4$ | $g_1$ | e | $g_5$ | $g_2$ |
| $g_4$ | $g_4$ | $g_3$ | $g_5$ | $g_2$ | e | $g_1$ |
| $g_5$ | $g_5$ | $g_2$ | $g_4$ | $g_1$ | $g_3$ | e |

Clearly it can be easily verified that the loop ring ZL is a right alternative ring.

**DEFINITION 2.2.15**: *Let R be a ring; R is said to be a SNA-right alternative ring if R has a subring S such that S is a SNA-subring of R and S is a right alternative ring that is (xy) y = x (yy) is true for all x, y $\in$ S. Similarly we can define SNA-left alternative ring. If R is simultaneously both SNA-right alternative and SNA-left alternative ring then we say R is a SNA-alternative ring.*

**DEFINITION 2.2.16**: *Let R be a non-associative ring. R is said to be a SNA-commutative ring if R has a proper subset S (S $\subset$ R) such that S is a commutative associative ring with respect to the operations of R.*

**DEFINITION 2.2.17**: *Let (R, +, .) and ($R_1$, +, .) be any two SNA-rings. A map $\phi$ : R to $R_1$ is a Smarandache NA-ring homomorphism (SNA-ring homomorphism) if $\phi$ : A $\rightarrow$ $A_1$ is a ring homomorphism where A $\subset$ R and $A_1$ $\subset$$R_1$ are just associative subrings of R and $R_1$, respectively. Clearly Ker $\phi$ is an ideal of A and never an S-ideal of R. In fact Ker $\phi$ may not be even an ideal of R.*

**DEFINITION 2.2.18**: *Let R be a SNA-ring if R has no SNA-ideals then we say R is SNA-simple ring.*

We define SNA-quotient ring for a given ring R and an SNA-ideal I; R / I is a Quotient ring as in the case of rings. One is not in a position to say whether the SNA-quotient rings of a SNA-maximal ideal is simple? The answer to this question is little complicated as we have the very definition of SNA-ideals only relative to SNA-subrings and not relative to the whole ring R. So one may presume that the quotient ring may have SNA-ideals relative to some other SNA-subrings.

Thus when we say the ring R is SNA-simple it is assumed the ring has no SNA-ideals related to any of the SNA-subrings. One may be wanting to know the concept of Smarandache maximal and SNA-minimal ideals. For the sake of helping the reader we define explicitly these types of ideals.

**DEFINITION 2.2.19**: *Let R be a ring. S a SNA-subring of R. We say an SNA-ideal I $\subset$ S is maximal if we have an SNA-ideal, K relative to S such that I $\subset$ K $\subset$ S then either I = K of K = S similarly an ideal J $\subset$ S is minimal if we have another ideal L $\subset$ S such that L $\subset$ J then L = (0) or L = J.*



Unlike in the case of minimal and maximal ideals we see in case of SNA-ideals the maximality or the minimality is only a relative concept for it follows by the very definition of them. Similarly we say a SNA-ideal is principal if it is generated by a single element. We say a SNA-ideal J of S, S a SNA-subring of R is prime if $xy \in J$ then x or y is in J.

**PROBLEMS:**

1. Let $L = \{m_1 = e, m_2, \ldots, m_n\}$ be a finite loop with $m_i^2 = e$ for all $m_i \in L$. $Z_2$ be the field of characteristic 2 and $Z_2L$ be the loop ring. If $\alpha \in Z_2 L$ is an idempotent of $Z_2L$ prove

    a. |supp $\alpha$| is an even number and $e \notin$ supp $\alpha$ (e the identity element of L).
    b. |supp $\alpha$| is an odd number and $e \in$ supp $\alpha$.
    c. If L is commutative then $Z_2 L$ has only trivial idempotents.

2. Let L be a loop given by the following table $L_7 (4)$.

    | . | e | 1 | 2 | 3 | 4 | 5 | 6 | 7 |
    |---|---|---|---|---|---|---|---|---|
    | e | e | 1 | 2 | 3 | 4 | 5 | 6 | 7 |
    | 1 | 1 | e | 5 | 2 | 6 | 3 | 7 | 4 |
    | 2 | 2 | 5 | e | 6 | 3 | 7 | 4 | 1 |
    | 3 | 3 | 2 | 6 | e | 7 | 4 | 1 | 5 |
    | 4 | 4 | 6 | 3 | 7 | e | 1 | 5 | 2 |
    | 5 | 5 | 3 | 7 | 4 | 1 | e | 2 | 6 |
    | 6 | 6 | 7 | 4 | 1 | 5 | 2 | e | 3 |
    | 7 | 7 | 4 | 1 | 5 | 2 | 6 | 3 | e |

    Does $Z_2L_7(4)$ have idempotents other than 1 and 0? Justify your answer.

3. Can $Z_2L_7(4)$ given in problem 2 have nontrivial

    a. quasi regular elements?
    b. regular elements?
    c. right quasi regular elements which are not left quasi regular?

4. Let $L_{15}(8)$ be a loop and Z be the ring of integers. $ZL_{15}(8)$ be the loop ring of the loop $L_{15}(8)$ over Z.

    a. Can $ZL_{15} (8)$ have idempotents? If so find them.
    b. Does $ZL_{15}(8)$ have right quasi regular elements?
    c. Prove no non-zero idempotent element of $ZL_{15}(8)$ can be right quasi regular.

5. Let $Z_3 = \{0, 1, 2\}$ be the prime field of characteristic 3 and $L_9(8)$ be the loop. $Z_3L_9(8)$ be the loop ring of the loop $L_9 (8)$ over $Z_3$, find



a. W ($Z_3L_9(8)$).
   b. J ($Z_3L_9(8)$).

6. For the loop ring $Z_3$ L where L is given by the following table:

    | . | e | 1 | 2 | 3 | 4 | 5 |
    |---|---|---|---|---|---|---|
    | e | e | 1 | 2 | 3 | 4 | 5 |
    | 1 | 1 | e | 3 | 5 | 2 | 4 |
    | 2 | 2 | 5 | e | 4 | 1 | 3 |
    | 3 | 3 | 4 | 1 | e | 5 | 2 |
    | 4 | 4 | 3 | 5 | 2 | e | 1 |
    | 5 | 5 | 2 | 4 | 1 | 3 | e |

   Find

      a. SNA-subrings of $Z_3L$.
      b. SNA-ideals if any in $Z_3L$.
      c. SNA-principal ideal if any.
      d. SNA-maximal / minimal ideal if any.

7. Find for the loop ring ZL where the loop L is given in problem (6) and Z the ring of integers SNA-subrings which has no SNA-ideals related to them.

8. Does there exists for loop ring $Z_7L_{13}$ (5) and $Z_9L_{15}$ (8) maps $\phi$ such that $\phi$ is a SNA-homomorphism? Find Ker $\phi_1$ if $\phi_1 : Z_7 \ A_1 \rightarrow Z_9 \ A_2$ where $A_1 = \{e, 7\}$ and $A_2 = \{e, 2\}$ are subgroups of the loops $L_{13}(5)$ and $L_{15}(8)$ respectively. Can $\phi_1$ be an isomorphism? Justify your claim.

9. Find the ideals of the loop ring $Z_2L_5$ (3).

      a. Do the two sided ideals of $Z_2L_5$ (3) form a modular lattice?
      b. Find only the right ideals of $Z_2L_5$ (3). What is its lattice structure?

10. For the loops rings $Z_4L_9$ (2) and $Z_4L_9$ (8) find

      a. The set of right ideals and its lattice structure.
      b. The set of left ideals and their lattice structure.
      c. Two sided ideals and their lattice structure.

11. Find for the loop rings $Z_4L_9$ (2) and $Z_4L_9(8)$, J ($Z_4L_9(2)$) and J ($Z_4L_9$ (8)).

12. Find for the loop rings given in problem 11, W ($Z_4L_9(2)$) and W ($Z_4L_9$ (8)).

## 2.3 Smarandache Elements in Loop Rings

In this section we define the concept of Smarandache zero divisors, Smarandache nilpotents, Smarandache idempotents, Smarandache units, Smarandache normal



elements, Smarandache semi-idempotents and Smarandache quasi-regular elements in NA-rings. Since the concept of these Smarandache elements do not in any way demand the NA-ring to be SNA-ring we proceed on to introduce these notions in any NA-rings.

**DEFINITION 2.3.1**: *Let R be a NA-ring with unit. We say $x \in R \setminus \{1\}$ is a Smarandache unit (S-unit) in R if there exist a $y \in R$ with $xy = 1$.*

*There exist $a, b \in R \setminus \{x, y, 1\}$ such that*

  a.  $xa = y$ or $ax = y$.
  b.  $yb = x$ or $by = x$ and
  c.  $ab = 1$.

**THEOREM 2.3.1**: *Let R be a NA ring. If x is a SNA-unit then x is a unit.*

*Proof*: Direct from the very definition.

It is left for the reader to construct an example of a unit which is not a SNA-unit.

**THEOREM 2.3.2**: *Let R be a NA-ring. All units in R in general need not be SNA-units of R.*

*Proof*: To prove this theorem we construct the following example. Let $Z_4 = \{0, 1, 2, 3\}$ be the ring of integers modulo 4. Let L be a loop given by the following table:

| . | e | $g_1$ | $g_2$ | $g_3$ | $g_4$ | $g_5$ |
|---|---|---|---|---|---|---|
| e | e | $g_1$ | $g_2$ | $g_3$ | $g_4$ | $g_5$ |
| $g_1$ | $g_1$ | e | $g_3$ | $g_5$ | $g_2$ | $g_4$ |
| $g_2$ | $g_2$ | $g_5$ | e | $g_4$ | $g_1$ | $g_3$ |
| $g_3$ | $g_3$ | $g_4$ | $g_1$ | e | $g_5$ | $g_2$ |
| $g_4$ | $g_4$ | $g_3$ | $g_5$ | $g_2$ | e | $g_1$ |
| $g_5$ | $g_5$ | $g_2$ | $g_4$ | $g_1$ | $g_3$ | e |

Let us take the loop ring $Z_4L$ identifying 1 with e we assume $1 \in Z_4L$. $\alpha = 1 + 2g_1 \in Z_4L$ $(1 + 2g_1)^2 = 1$ is a unit in $Z_4L$ but it is not a S-unit of $Z_4L$ Thus every unit in general is not a S-unit.

**DEFINITION 2.3.2**: *Let R be a non-associative ring. An element $a \in R \setminus \{0\}$ is called a Smarandache idempotent (S-idempotent) of R if $a^2 = a$. There exists $b \in R \setminus \{a\}$ such that*

  a.  $b^2 = a$.
  b.  $ab = b$ or $ba = (ba = a)$ $(ab = a)$.

*'or' in the mutually exclusive sense.*



**THEOREM 2.3.3**: *Let R be a non-associative ring. If $x \in R$ is a S-idempotent then it is an idempotent of R.*

*Proof*: Straightforward by the very definition.

To find the converse we study and analyze the following example:

*Example 2.3.1*: Let $Z_5 = \{0, 1, 2, 3, 4\}$ be the prime field of characteristic 5. Let $L_5(3)$ be the loop given by the following table:

| . | e | $g_1$ | $g_2$ | $g_3$ | $g_4$ | $g_5$ |
|---|---|---|---|---|---|---|
| e | e | $g_1$ | $g_2$ | $g_3$ | $g_4$ | $g_5$ |
| $g_1$ | $g_1$ | e | $g_4$ | $g_2$ | $g_5$ | $g_3$ |
| $g_2$ | $g_2$ | $g_4$ | e | $g_5$ | $g_3$ | $g_1$ |
| $g_3$ | $g_3$ | $g_2$ | $g_5$ | e | $g_1$ | $g_4$ |
| $g_4$ | $g_4$ | $g_5$ | $g_3$ | $g_1$ | e | $g_2$ |
| $g_5$ | $g_5$ | $g_3$ | $g_1$ | $g_4$ | $g_2$ | e |

Consider $\alpha = 1 + g_1 + g_2 + g_3 + g_4 + g_5$ in $Z_5L_5(3)$. Clearly $\alpha^2 = \alpha$. Thus $\alpha$ is an idempotent. Take $\beta = 3 + 3g_1$, $\beta^2 = (3 + 3g_1)^2 = 3 + 3g_1$. $(3 + 3g_1)(1 + g_1 + g_2 + g_3 + g_4 + g_5) = 1 + g_1 + g_2 + g_3 + g_4 + g_5$ Take $X = 4 + 4g_1 + 4g_2 + 4g_3 + 4g_4 + 4g_5$.

$X^2 = \alpha$, further $X\alpha = X$. Thus $\alpha$ is a S-idempotent but X is not an idempotent. $\beta = 3 + 3g_1$, take $Y = 2 + 2g_1$ we have $Y^2 = 3 + 3g_1$, $\beta Y = Y$.

Thus $\beta$ is a S-idempotent. Hence this example shows us the existence of idempotents, which are S-idempotent in a NA-ring. Now we give a nice theorem for the existence of S-idempotents in case of loop rings $Z_pL_n(m)$ for some suitable n.

**THEOREM 2.3.4**: *Let $Z_p = \{0, 1, 2, .., p-1\}$ be the prime field of characteristic p ($p > 2$). $L_n(m)$ be the loop where we choose $n = p$. Then $\alpha = (1 + g_1 + \ldots + g_p)$ and*

$$\beta = \left(\frac{p+1}{2} + \frac{p+1}{2}g_i\right);$$

$g_i \in L_p(m)$ are S-idempotents of the loop ring $Z_pL_p(m)$ for any loop $L_p(m) \in L_p$.

*Proof*: Let $\alpha = (1 + g_1 + \ldots + g_p)$ be in $Z_pL_p(m)$, clearly $\alpha^2 = \alpha$, so $\alpha$ is an idempotent in $Z_pL_p(m)$. Take $X = (p - 1) + (p - 1)g_1 + \ldots + (p - 1)g_p$. Clearly $X^2 = \alpha$ and $X\alpha = X$. Hence $\alpha$ is a S-idempotent take

$$\beta = \frac{p+1}{2} + \frac{p+1}{2}g_i;$$

$g_i \in L_p(m)$. Easily checked $\beta^2 = \beta$ so $\beta$ is an idempotent in $Z_pL_p(m)$. Consider



$$Y = \frac{p-1}{2} + \frac{p-1}{2} g_i;$$

$\beta Y = Y$. Thus elements of the form $\alpha$ and $\beta$ in $Z_p L_p(m)$ are S-idempotents.

**DEFINITION 2.3.3**: *Let R be a NA ring. A element $a \in R \setminus \{0\}$ is said to be a Smarandache zero divisor (S-zero divisor) if $a.b = 0$ for some $b \neq 0$ in R, and there exists $x, y \in R \setminus \{0, a, b\}$ $x \neq y$; such that*

  a. *$a.x = 0$ or $x.a = 0$.*
  b. *$b.y = 0$ or $y.b = 0$.*
  c. *$x.y \neq 0$ or $y.x \neq 0$.*

**THEOREM 2.3.5**: *Let R be a NA-ring. If $x \in R \setminus \{0\}$ is a S-zero divisor then x is a zero divisor in R.*

*Proof*: Direct by the very definition of S-zero divisors. Now it may happen that a S-zero divisor may satisfy all the three conditions a, b, c still it may behave in an entirely different way.

We define Smarandache pseudo zero divisors.

**DEFINITION 2.3.4**: *Let R be a NA-ring Let $x \in R \setminus \{0\}$ be a zero divisor in R. i.e., there exists $y \in R \setminus \{0\}$ with $x.y = 0$. We say x is a Smarandache pseudo zero divisor (S-pseudo zero divisor) if their exist $a \in R \setminus \{x, y, 0\}$ with*

$$a.y = 0 \text{ or } a.x = 0.$$
$$a^2 = 0.$$

**THEOREM 2.3.6**: *Let R be a NA-ring. If $x \in R$ is a S-pseudo zero divisor of R then x is a zero divisor of R.*

*Proof*: Follows from the fact that if x is a S-pseudo divisor then it is a zero divisor.

We illustrate them with examples.

*Example 2.3.2*: Let $Z_2 = \{0, 1\}$ be the prime field of characteristic two. $L_5(3)$ be the loop given by the following table:

| . | e | 1 | 2 | 3 | 4 | 5 |
|---|---|---|---|---|---|---|
| e | e | 1 | 2 | 3 | 4 | 5 |
| 1 | 1 | e | 4 | 2 | 5 | 3 |
| 2 | 2 | 4 | e | 5 | 3 | 1 |
| 3 | 3 | 2 | 5 | e | 1 | 4 |
| 4 | 4 | 5 | 3 | 1 | e | 2 |
| 5 | 5 | 3 | 1 | 4 | 2 | e |



We shall denote i by $g_i$ and e by 1. Let $Z_2L_5(3)$ be the loop ring of the loop $L_5(3)$ over $Z_2$. Clearly $x = 1 + g_1 \in Z_2L_5(3)$ is a zero divisor for $y = 1 + g_1 + g_2 + g_3 + g_4 + g_5$ is such that $xy = 0$ and $a = 1 + g_2$ is such that $a.y = 0$ with $a^2 = 0$ so x is a S-pseudo zero divisor and is not a S-zero divisor. It is an obvious zero divisor in $Z_2L_5(3)$.

**THEOREM 2.3.7**: *Let R be a NA-ring a S-pseudo zero divisor in general is not a zero divisor.*

*Proof*: From the above example we see $x = 1 + g_1$ is a S-pseudo zero divisor and is not a S-zero divisor of $Z_2L_5(3)$.

In view of this example we have a nice theorem which guarantees the existence of S-pseudo zero divisors in loop rings using the new class of loops.

**THEOREM 2.3.8**: *Let $Z_2 = \{0, 1\}$ be the prime field of characteristic two and $L_n(m) \in L_n$ be a loop in the class of loops. Then all elements of the form $x = 1 + g_i \in Z_2 L_n(m)$ where $g_i \in L_n(m)$ are S-pseudo zero divisors of $Z_2L_n(m)$.*

*Proof*: Given $x = 1 + g_i \in Z_2L_n(m)$, clearly for $y = 1 + g_1 + \ldots + g_n$ we have $xy = 0$ so x is a zero divisor in $Z_2L_n(m)$. Take $a = (1 + g_j)$, $i \neq j$ with $a.y = 0$ and we have $a^2 = 0$. Hence $x = 1 + g_i$ is a S-pseudo zero divisor of $Z_2L_n(m)$.

We still modify by replacing $Z_2$ by $Z_m$ where $m = 2p$, p an odd prime; in this case we have the following results:

**THEOREM 2.3.9**: *Let $L_n(m) \in L_n$ be a loop $Z_t$ be the ring of integers modulo t where t = 2p where p is an odd prime. $Z_tL_n(m)$ is the loop ring. $x = p + pg_i \in Z_tL_n(m)$ is a S-pseudo zero divisor for all $g_i \in L_n(m)$.*

*Proof*: Let $x = p + pg_i \in Z_tL_n(m)$, consider $y = 1 + g_1 + \ldots + g_n$; we have $x.y = 0$. Take $a = p + pg_j$ ($i \neq j$) then $ay = 0$ with $a^2 = 0$, hence the claim.

The following have been proved about zero divisor in case of power associative and diassociative loop algebras.

**THEOREM [56]**: *Let K be a field, L a simply ordered loop. Then the loop algebra KL has no nontrivial zero-divisors.*

*Proof*: Left for the reader to refer [56].

**THEOREM [56]**: *Let K be a field, L a finite power associative or diassociative loop. Then the loop ring KL has nontrivial divisors of zero.*

*Proof*: KL be the loop algebra. Since L is a diassociative or a power associative loop we have a proper subgroup H; $H \subset L$. Now KH is a group algebra; since $|H| < \infty$ we have for $g \in H$, $g^n = 1$.

Hence $(1 - g)(1 + g + \ldots + g^{n-1}) = 0$. a nontrivial divisor of zero.



In view of these we have the following theorem:

**THEOREM 2.3.10**: *Let KL be a S-loop ring where L is a diassociative loop or a power associative loop of finite order. Then KL has divisors of zero.*

*Proof*: Obvious

Does the loop ring KL have S-divisors of zero?

We define Smarandache weak zero divisor in the following.

**DEFINITION 2.3.5**: *Let R be a NA ring. An element $x \in R \setminus \{0\}$ is a Smarandache weak zero divisor (S-weak zero divisor) if there exist a $y \in R \setminus \{0, x\}$ such that x.y = 0 satisfying the following condition:*

*There exists a, b $\in R \setminus \{x, y, 0\}$ such that*

$$a.x = 0 \text{ or } x.a = 0$$
$$b.y = 0 \text{ or } y.b \ 0$$
$$x.y = 0 \text{ or } y.x = 0$$

*Thus in case of S-weak divisor we have another pair of elements which is a zero divisor.*

**THEOREM 2.3.11**: *Let K be a field and L a diassociative finite loop or a power associative finite loop. The loop ring KL has S-weak divisors of zero.*

*Proof*: We have KL to be a S-loop ring. Also $g \in L$ is such that $g^n = 1$ and $g \in H \subset L$ where H is a group.

$$(1 - g)(1 + g + \ldots + g^{n-1}) = 0$$

Take $a = 1 - g^r$, $r \neq 0$, $r > 1$

$$(1 - g^r)(1 + g + \ldots + g^{n-1}) = 0$$

Let $b = 3 + 3g + \ldots + 3g^{n-1}$, then $b(1 - g) = 0$. Now $a.b = 0$. Hence KL has S-weak divisors of zero.

**THEOREM 2.3.12**: *A S-weak zero divisor of a loop ring in general is not a S-zero divisor but it is a zero divisor.*

*Proof*: Follows from the definition that S-weak zero divisor is a zero divisor. From theorem 2.3.11 it is clearly not a S-zero divisor only a S-weak zero divisor.

It is natural to see that if an element is a zero divisor still it need not be S-weak zero divisor. The reader is advised to construct examples.



**DEFINITION 2.3.6**: *Let KL be a loop ring of a loop L over the field K, an element $\alpha \in$ KL is called a normal element of KL if $\alpha KL = KL\alpha$. If every element in KL is a normal element then we say KL is a normal loop ring.*

**DEFINITION 2.3.7**: *Let KL be a loop ring. S be a proper subring of KL, we say S is a normal subring with respect to a subset T of KL if*

$$xS = Sx \text{ for all } x \in T$$
$$x(yS) = (xy)S$$
$$(Sx)y = S(xy) \text{ for all } x, y \in T.$$

*If in particular T = KL then we say the subring S is a normal subring of KL.*

**THEOREM 2.3.13**: *Let C(L) be the center of the loop L and K any field, KL the loop ring of the loop L over K. KL has a normal subring.*

*Proof*: Obvious form the fact if C(L) is the center of L then KC(L) is a loop ring contained in KL which is a normal subring of KL.

**THEOREM 2.3.14**: *Let L be a loop. S a normal subloop of L. K any field the loop ring KL has KS to be a normal subring.*

*Proof*: Follows from the fact $KS \subset KL$ is such that KS is a normal subring, as S is a normal subloop of L. Hence the claim.

Now we proceed on to define Smarandache normal subring in a non-associative ring.

**DEFINITION 2.3.8**: *Let R be a NA-ring. Suppose A be a S-subring of R. We say an element $x \in R$ to be Smarandache normal element (S-normal element) if xA = Ax. If every $x \in R$ is such that xA = Ax we call R a Smarandache normal ring (S-normal ring) relative to A. Now it is pertinent to mention here that Smarandache normality is a relative concept, for relative to one S-subring A the ring R may be Smarandache normal and relative to some other S-subring B of R, the ring R may fail to be a Smarandache normal.*

**THEOREM 2.3.15**: *Let R be a S-ring which is normal. If R has S-subring then R is S-normal. In fact R is S-normal with respect to every S-subring A of R.*

*Proof*: Follows by the very definition.

**THEOREM 2.3.16**: *Let R be a NA-ring. If R is a S-normal relative to a S-subring A of R then R need not in general be a normal ring.*

*Proof*: The reader is advised to obtain one such example.

Thus we see S-normality of R in general need not imply normality of R or R is to be S-normal relative to every S-subring of R.



**DEFINITION 2.3.9**: *Let R be a ring. if for every S-subring A of R is Smarandache normal relative to A then we call the ring R to be Smarandache strongly normal ring (S-strongly normal ring).*

**THEOREM 2.3.17**: *If R is a S-strongly normal ring then R is a S-normal ring and if R is a S-normal ring then R need not in general be a S-strongly normal ring.*

*Proof*: This is left as an exercise for the reader to prove.

We now proceed on to define the concept of strictly right loop ring and then we will generalize it to the notion of Smarandache strictly right loop ring.

**DEFINITION 2.3.10**: *Let R be a non-associative ring. R is said to be a strictly right ring if the set of all right ideals of R is ordered by inclusion. If all ideals are included by inclusion we call them as strictly ideal rings.*

**DEFINITION 2.3.11**: *Let R be a non-associative ring. If the set of all S-right ideals of R are ordered by inclusion then we call R a Smarandache strongly right ideal ring (S-strongly right ideal ring). In particular if all S-ideals of R are ordered by inclusion then we call the ring R to be a Smarandache strongly ideal ring (S-strongly ideal ring).*

We do not have any form of relations existing between them but we have several problems in this direction suggested in the last chapter.

**THEOREM 2.3.18**: *Let R be a S-strictly right ideal ring then R is a S-ring.*

*Proof*: Obvious by the very definition of these concepts.

*Example 2.3.3*: Let $Z_2 = \{0, 1\}$ be the prime field of characteristic and let L be a loop given by the following table:

| . | e | 1 | 2 | 3 | 4 | 5 |
|---|---|---|---|---|---|---|
| e | e | 1 | 2 | 3 | 4 | 5 |
| 1 | 1 | e | 3 | 5 | 2 | 4 |
| 2 | 2 | 5 | e | 4 | 1 | 3 |
| 3 | 3 | 4 | 1 | e | 5 | 2 |
| 4 | 4 | 3 | 5 | 2 | e | 1 |
| 5 | 5 | 2 | 4 | 1 | 3 | e |

$Z_2L$ be the loop ring of the loop L over the ring $Z_2$. It is left for the reader to find all the S-right ideals of $Z_2L$, right ideals of $Z_2L$ and ideals and S-ideals of $Z_2L$.

The meagerly used concept is the notion of semi-idempotents introduced by [30].

**DEFINITION [30]**: *Let R be a ring. An element $\alpha \neq 0$ in R is said to be a semi-idempotent if and only if $\alpha$ is not in the two sided ideal of R generated by $\alpha^2 - \alpha$, that is $\alpha \notin R (\alpha^2 - \alpha) R$ or $R (\alpha^2 - \alpha) R = R$, where by $R (\alpha^2 - \alpha) R$ we mean the two sided ideal of R generated by $\alpha^2 - \alpha$. So $R (\alpha^2 - \alpha) R$ just denotes the ideal*



*generated by $\alpha^2 - \alpha$. We shall denote the set of all semi-idempotents of a ring R by SI (R).*

*Example 2.3.4*: Let K be any field of characteristic 0. L be a loop given by the table.

| . | 1 | a | b | c | d | e |
|---|---|---|---|---|---|---|
| 1 | 1 | a | b | c | d | e |
| a | a | 1 | d | b | e | c |
| b | b | d | 1 | e | c | a |
| c | c | b | e | 1 | a | d |
| d | d | e | c | a | 1 | b |
| e | e | c | a | d | b | 1 |

The loop ring KL has

$$\frac{3+a}{2}, \frac{3+b}{2}, \frac{3+c}{2}, \frac{3+d}{2} \text{ and } \frac{3+e}{2}$$

as some of its nontrivial idempotents in it.

**THEOREM 2.3.19**: *Let KL be the loop ring of any loop L over the field K. If $\alpha (\neq 0)$ is a semi-idempotent but not a unit in KL, then $\alpha - 1$ is not a unit of KL.*

*Proof*: Suppose $\alpha - 1$ be unit in KL. Then there is an element $\beta$ of KL such that $(\alpha - 1) \beta = 1$. Thus we would have $\alpha = (\alpha^2 - \alpha) \beta \in (\alpha^2 - \alpha)$ KL, but $(\alpha^2 - \alpha)$ KL $\neq$ KL because $\alpha$ is a semi-idempotent, hence $\alpha(\alpha - 1) = \alpha^2 - \alpha$ is not a unit hence the result.

Now we define Smarandache semi-idempotents.

**DEFINITION 2.3.12**: *Let R be a NA-ring. An element $\alpha \in R$ which is a semi-idempotent is called a Smarandache semi-idempotent (S-semi-idempotent) if in R if the S-ideal generated by $\alpha^2 - \alpha$ does not contain $\alpha$. We denote the collection of S-semi-idempotents in R by SSI (R).*

**THEOREM 2.3.20**: *Let R be a NA ring if $\alpha$ is a S-semi-idempotent then $\alpha$ is a semi-idempotent of R.*

*Proof*: Direct by the very definition of S-semi-idempotents.

The proof of the following theorem is left as an exercise for the reader to prove.

**THEOREM 2.3.21**: *Let R be a ring. If $0 \neq \alpha \in R$ is a semi-idempotent then $\alpha$ in general need not be a S-semi-idempotent.*

Clearly if R is not a S-ring then it is impossible for the element $0 \neq \alpha \in R$ which is a semi-idempotent; to be a S-semi-idempotent. Now we proceed on to define Smarandache quasi-regular elements in a loop ring.



**DEFINITION 2.3.13**: *Let R be a non-associative ring. We say an element $x \in R$ is called a Smarandache right quasi regular (S-right quasi regular) if there exists y and z in R such that*

$$x \circ y = x + y - xy = 0$$
$$\text{and} \quad x \circ z = x + z - xz = 0$$
$$\text{but} \quad y \circ z = y + z - yz \neq 0$$
$$\text{and} \quad z \circ y = z + y - zy \neq 0.$$

*Similarly we define Smarandache left quasi regular elements (S-left quasi regular elements). x will be called Smarandache quasi regular (S-quasi regular) if it is simultaneously S-left quasi regular and S-right quasi regular.*

**THEOREM 2.3.22**: *Let R be a NA-ring, if $x \in R$ is S-right (left) quasi regular then x is right (left) quasi regular.*

*Proof*: Follows from the very definitions of these concepts.

**THEOREM 2.3.23**: *Let R be a non-associative ring. If $x \in R$ is right quasi regular then x need not in general be S-right quasi regular.*

*Proof*: Let $Z_2 = \{0, 1\}$ prime field of characteristic two, L loop given by the following table:

| . | e | a | b | c | d |
|---|---|---|---|---|---|
| e | e | a | b | c | d |
| a | a | e | c | d | b |
| b | b | d | a | e | c |
| c | c | b | d | a | e |
| d | d | c | e | b | a |

$Z_2L$ be the loop ring. We have $x = e + a + b + c \in Z_2L$ is right quasi regular for $x \circ (a + b) = 0$ for $x \circ (a + b) = e + a + b + c + a + b - a - e - d - b - b - c - a - d = 0$.

It is easily verified that x is not S-right quasi regular.

**DEFINITION 2.3.14**: *Let R be a NA-ring. We say a S (left) right ideal is Smarandache right (left) quasi regular (S- right (left) quasi regular) if each of its elements are S-right (left) quasi regular. The Smarandache Jacobson radical (S- Jacobson radical) [SJ(R)] of a ring R is defined as follows: $SJ(R) = \{a \in R \,/\, aR$ is a S-quasi regular ideal of R\}.*

**PROBLEMS:**

1. Find the units in $Z_6L_5(2)$ which are SNA-units.
2. Find all quasi regular elements in $Z_6L_5(2)$ which are S-quasi regular.
3. Find all S-units in $Z_7L_9(8)$.



4. Find those idempotents in $Z_2L_{13}(6)$ which are S-idempotents.
5. Find SJ $(Z_3(L_5(3)))$.
6. Find all S-quasi regular elements of $Z_7L_7(3)$.
7. Can $ZL_9(8)$ have S-right quasi regular elements which are not S-left quasi regular elements?
8. Find all quasi regular elements of $Z_6L_5(3)$ which are not S-quasi regular elements.
9. Find all S-semi-idempotents of $Z_7L_9(8)$.
10. Determine SSI $(ZL_5(3))$ and SI $(ZL_5(3))$. Do these sets form any proper algebraic structures?
11. Is the loop ring $Z_8L_9(2)$ a strictly right ring?
12. Find all the S-ideals of $Z_{12}L_5(3)$.

## 2.4 Smarandache substructures in loop rings

In this section we study the substructure properties like Marot loop rings, mod p-envelope of a loop ring, orthogonal ideals, strongly commutative loop ring, inner commutative loop ring, pseudo commutative loop rings and the lattice of substructures of loop rings like ideals, subrings, and their Smarandache analogue. Several interesting results are proved and we see in loop rings several concepts are not fully analyzed. So these Smarandache analogue is still in a dormant state. Thus in this section we introduce several new concepts to loop rings.

The concept of Marot loops rings was introduced and studied in [63]. Here we just recall the concepts and define Smarandache Marot loop rings.

**DEFINITION [63]**: *A commutative loop ring RL with identity is a Marot loop ring if each regular ideal of RL is generated by a regular element of RL, where by a regular element of RL we mean a non-zero divisor of RL.*

**THEOREM [63]**: *Let F be a field. L be the commutative ordered loop without elements of finite order. Then the loop ring FL is a Marot loop ring.*

*Proof*: Since F is a field and L is a loop which is ordered and has no elements of finite order we see the loop ring FL has no zero divisors; hence every element is regular so all ideals in FL will be generated only regular elements. Thus FL is a Marot ring.

In view of this result we have the following theorem that is left as an exercise for the reader to prove.

**THEOREM [63]**: *Let F be a field. FL is a Marot loop ring without divisors of zero if and only if L is a commutative ordered loop without elements of finite order.*

**THEOREM [63]**: *If R is a Marot ring and L is a commutative loop; then the loop ring RL is a Marot ring.*



*Proof*: RL is a commutative loop ring as both R and L are commutative. To prove RL is a Marot ring we need to show every regular ideal generated by regular elements of RL. To prove this we shall prove.

a. If a regular ideal I is generated by a regular element of I then I has no proper divisors of zero.

b. If a regular ideal I is generated by a divisor of zero than I cannot contain regular elements.

*Proof of (a)*: Let I be generated by a regular element $\alpha$. Let $\beta \neq 0$ in I be such that $\beta\gamma = 0$ ($\gamma \neq 0$). Now $\beta = \left(\sum_i \alpha\delta_i\right), \beta\gamma = 0 = \left(\sum \alpha\delta_i\right)\gamma = \alpha\sum_i \delta_i\gamma$ whether associative or not we see $\alpha.\Sigma\delta_i\gamma = 0$, a contradiction to our assumption as $\alpha$ is a regular element of I. Hence I cannot contain divisors of zero.

*Proof of (b)*: Let I be generated by a divisor of zero say $\beta$, ($\beta \neq 0$, $\gamma \neq 0$) with $\beta\gamma = 0$. Let $\alpha \in$ I be a regular element of I. $\alpha = \Sigma\beta\delta_i$, $\alpha\gamma = (\Sigma\beta\delta_i)\gamma = \gamma(\Sigma\beta\delta_i) = \gamma\beta(\Sigma\delta_i) = 0$. That is $\alpha\gamma = 0$ a contradiction, so $\alpha \notin$ I. Hence the theorem.

**DEFINITION 2.4.1**: *Let R be a NA ring. We say R is a Smarandache weak Marot ring (S-weak Marot ring) if R has a commutative S-subring A such that all S-ideals associated with A are regular. We say R is a Smarandache Marot ring (S-marot ring) if for every S-subring A of R is commutative; we have all S-ideals associated with every S-subring A are regular.*

**THEOREM 2.4.1**: *If R is a S-Marot ring then R is a S-weak Marot ring.*

*Proof*: Direct by the very definition of these concepts.

We leave the following theorems as an exercise for the reader to prove.

**THEOREM 2.4.2**: *Every S-weak Marot ring in general need not be a S-Marot ring.*

**THEOREM 2.4.3**: *Let R be a S-Marot ring then R need not in general be a Marot ring.*

To prove the two theorems the reader is advised to construct examples. It is also to be noted even if R is not a commutative ring still R can be a S-Marot ring or S-weak Marot ring.

The concept of mod-p-envelope of G was studied in [38]. The structure of mod-p-envelope was studied in case of loop in [62]. It is proved in [62] the mod-p-envelope of a loop was only a groupoid.

**DEFINITION 2.4.2**: *Let L be a loop and K be any field. We say $L* = 1 + U$ where $U = \{\alpha = \Sigma \alpha_i m_i \in KL / \Sigma\alpha_I = 0\}$ is the mod p envelope of L.*

We study L∗ in our own way analogous to mod p-envelope of a group G.



***Example 2.4.1***: Let L be a loop given by the table below and K = {0, 1} be a field of characteristic two. Then $L^*$, the mod p-envelope of L is a groupoid with $(1 + a + b + c + d)^2 = (1 + b + a + c + d)$ and the number of elements in $L^*$ is equal to 16.

| * | 1 | a | b | c | d |
|---|---|---|---|---|---|
| 1 | 1 | a | b | c | d |
| a | a | d | c | 1 | b |
| b | b | 1 | d | a | c |
| c | c | b | 1 | d | a |
| d | d | c | a | b | 1 |

KL is the loop ring of L over K. U = {0, 1 + a, 1 + b, 1 + c, 1 + d, a + b, a + c, a + d, b + c, b + d, c + d, 1 + a + b + c, 1 + a + b + d, 1 + a + c + d, 1 + b + c + d, a + b + c + d}. Let $L^* = 1 + U$; $\alpha = 1 + a + b + c + d \in L*$ with $\alpha^2 = \alpha$ and $1 \in L^*$, thus $L^*$ is a groupoid with 16 elements in it.

**THEOREM [72]**: *Let L be a commutative finite loop of order 2n in which the square of every element is one and let K = {0, 1} be the prime field of characteristic two. Then $L^*$ is a loop which is commutative such that square of each element is 1 and order of $L^*$ is $2^{2n-1}$.*

*Proof*: Let L = {$a_1 = 1, a_2, ..., a_{2n}$ / $a_i^2 = 1$, $a_i a_j = a_j a_i$; i, j = 1, 2, ..., 2n}. Now U = {0, 1 + $a_2$, ..., $a_{2n-1} + a_{2n}$, ..., sum of 2n elements taken at a time}. Now $L^* = 1 + U$, so $|L^*| = 2^{2n-1}$. Clearly $L^*$ is a commutative loop as L is a commutative loop. So is $L^*$ using the fact $a_i^2 = 1$ for i = 1, 2, ..., 2n and characteristic of K is 2.

**THEOREM 2.4.4**: *Let L be a loop of order 2n + 1 commutative or otherwise; K = {0, 1} be the prime field of characteristic two. $L^*$ is a groupoid of order $2^{2n}$ with a nontrivial idempotent in it.*

*Proof*: The proof is left as an exercise for the reader to prove.

**THEOREM 2.4.5**: *Let L be a finite loop with an element $x \in L$ such that $x^2 = 1$. Let $Z_p$ = {0, 1, 2, ..., p – 1} be a prime field of characteristic p, p > 2. Then $L^*$ is a groupoid with a nontrivial idempotent in it and the order of $L^*$ is $p^{n-1}$ where $|L| = n$.*

*Proof*: Straightforward, hence left for the reader as an exercise.

**THEOREM 2.4.6**: *Let L be a loop of order p + 1 and K = $Z_p$ = {0, 1, 2, ..., p – 1} prime field of characteristic p. Then $L^*$ is a groupoid with a nontrivial idempotent in it.*

*Proof*: Left as an exercise for the reader to prove.

Now we proceed on to define Smarandache mod p-envelope of a loop.



**DEFINITION 2.4.3**: *Let L be a S-loop and K any field, KL the loop ring of K over the loop L. The Smarandache mod p-envelope of L (S-mod p-envelope of L) is the mod p envelope of the S-subloop A of L denoted by $S(L^*)$. Thus unlike in the case of loops where the mod p envelope is unique, in case of S-mod p-envelope of a loop we may have several S-mod p-envelope depending on the number of S-subloops in L. If L has no S-subloop and even if L is a S-loop still the S-mod p envelope of L will only the empty set.*

Thus it is an interesting research problem to find the number of S-mod p-envelopes associated with a loop L.

**THEOREM 2.4.7**: *Let L be a loop having no S-subloops or L is not a S-loop and K any field. Then $S(L^*) = \phi$.*

*Proof*: Straightforward from the very definition of $S(L^*)$.

**THEOREM 2.4.8**: *Let $L_n(m) \in L_n$ where n is a prime. K any field, $S(L^*) = \phi$.*

*Proof*: Clearly $L_n(m) \in L_n$ when n is a prime is a S-loop; but $L_n(m)$ has no proper S-subloop, hence $S((L_n(m))^* = \phi$.

Now we study the loop rings when the loops are unique product or two unique product loops.

**DEFINITION 2.4.4**: *A loop L is called a two unique product loop (t.u.p loop), if given any two non-empty finite subsets A and B of L with $|A| + |B| > 2$, there exists atleast two distinct elements x and y of L that have unique representation in the form $x = ab$, $y = cd$ with $a, c \in A$ and $b, d \in B$. A loop L is called a unique product loop (u.p. loop) if, when A and B are non-empty finite subsets of L, then there always exists atleast one $x \in L$ which has unique representation in the form $x = ab$ with $a \in A$ and $b \in B$.*

[52] proved in case of groups t.u.p and u.p are equivalent for groups. In case of loops it an open question.

**THEOREM 2.4.9**: *Let R be any ring with identity and L be a t.u.p. power associative or diassociative loop. Then the following are equivalent.*

  i. *$U(KS) = \{\Sigma \alpha_g g \,/\, \text{there exists a } \beta_g \in R \text{ with } \Sigma \alpha_g \beta_g^{-1} = 1 \text{ and } \alpha_g \beta_h = 0$ whenever $gh \neq 1\}$, where S is the subgroup generated by single element or a pair of elements.*

  ii. *R has no nonzero nilpotent element.*

*Proof*: To prove (i) implies (ii) follows from the fact that if R has nilpotent elements say $\gamma \in R$ then $1 + \gamma g$ is a unit in RS. For if $\gamma^n = 0$ then

$(1 + \gamma g)(1 - \gamma g + \gamma^2 g^2 - \gamma^3 g^3 + \ldots \pm \gamma^{n-1} g^{n-1}) = 1 - \gamma g + \gamma^2 g^2 - \ldots \pm \gamma^{n-1} g^{n-1} + \gamma g - \gamma^2 g^2 + \ldots \pm \gamma^{n-1} g^{n-1} \pm \gamma^n g^n = 1$ as $\gamma^n = 0$ and all other terms cancel out. To show (ii)



implies (i) the result follows from the fact if R has non-zero nilpotents and p, q ∈ RL where L is a t.u.p loop. If pq = 1, where p = Σ $α_g$ g and q = Σ $β_h$ h then $α_g$ $β_h$ = 0 when gh ≠ 0 with some factors using t.u.p we get the result. Hence the claim.

**DEFINITION 2.4.5**: *Let L be a loop. We say L is a Smarandache unique product loop (S.u.p loop) if there exists a S-subloop A of L such that A is a u.p loop. We call the loop L a Smarandache strongly u.p loop (S-strongly u.p.loop) if every S-subloop A of L is a u.p subloop of L.*

**DEFINITION 2.4.6**: *Let L be a loop; L is said to be Smarandache t.u.p loop (S-t.u.p loop) if L has a S-subloop A such that A is a t.u.p subloop of L. If every S-subloop of L is a t.u.p loop then we say L is a Smarandache strong t.u.p (S-strong t.u.p) loop.*

**THEOREM 2.4.10**: *If L is a S-u.p loop or S.t.u.p loop then L is a S-loop.*

*Proof*: Obvious by the very definitions.

**THEOREM 2.4.11**: *Let L be a S-loop which is a u.p loop or a t.u.p loop If L has no S-subloops then L is not a S.u.p loop or a S.t.u.p loop.*

*Proof*: It is left as an exercise for the reader to prove.

**THEOREM 2.4.12**: *Let $L_p(m) ∈ L_p$, if p is a prime; no $L_p(m)$ is a S.u.p loop or a S.t.u.p loop.*

*Proof*: Follows from the fact if p is a prime then $L_p$(m) has no S-subloop. So $L_p$(m) is not a S.t.u.p loop or a S.u.p loop.

A loop L is Hamiltonian if every subloop is normal. In view of this we define Smarandache Hamiltonian loop as follows:

**DEFINITION 2.4.7**: *Let L be a loop if every S-subloop of L is normal then we say L is a Smarandache Hamiltonian loop (S-Hamiltonian loop).*

**THEOREM 2.4.13**: *Let L be a finite diassociative Hamiltonian loop which is the direct product of A × T × H where A is an abelian group with elements of finite order, T an abelian group of exponent 2 and H satisfies (x, y, z) = 1 and K any field. The loop algebra KL is a direct sum of irreducible KA-modules, KT-modules and KH-modules.*

*Proof*: Since L is a finite diassociative Hamiltonian loop. KA, KT and KH are group algebras of finite groups A, T and H respectively. By Maschke's theorem KA is a direct sum of irreducible KA-modules, KT is the direct sum of irreducible KT-modules and KH is the direct sum of irreducible KH modules. Hence KL is a direct sum of irreducible KA, KT and KH-modules. Hence the theorem.

**DEFINITION 2.4.8**: *Let RL be the loop ring of the loop L over the ring R. Two non-zero ideals J and I of RL are orthogonal if I.J = {i.j = 0 for all i ∈ I and j ∈ J}. An ideal I of RL is self orthogonal if I.I = (0).*



*Example 2.4.2*: Let $Z_2 = \{0, 1\}$ be the prime field of characteristic two and L be the loop given by the following table:

| * | 1 | a | b | c | d | e | f | g |
|---|---|---|---|---|---|---|---|---|
| 1 | 1 | a | b | c | d | e | f | g |
| a | a | 1 | c | e | g | b | d | f |
| b | b | g | 1 | d | f | a | c | e |
| c | c | f | a | 1 | e | g | b | d |
| d | d | e | g | b | 1 | f | a | c |
| e | e | d | f | a | c | 1 | g | b |
| f | f | c | e | g | b | d | 1 | a |
| g | g | b | d | f | a | c | e | 1 |

Consider the loop ring $Z_2 L$. Take two ideals given by $I = \{0, 1 + a + b + c + d + e + f + g\}$ and $J = \{0, \alpha = \Sigma m_i / m_i \in L\}$ such that the sum of the coefficients of $m_i$ in the sum $\alpha$ is zero which are clearly nontrivial and are such that $I.J = \{0\}$. Thus I and J are orthogonal ideals of $Z_2L$. It is easily verified that $I.I = \{0\}$.

*Example 2.4.3*: Let $Z_2 = \{0, 1\}$ be the prime field of characteristic two and L be the loop given by the following table:

| . | 1 | a | b | c | d |
|---|---|---|---|---|---|
| 1 | 1 | a | b | c | d |
| a | a | 1 | c | d | b |
| b | b | d | a | 1 | c |
| c | c | b | d | a | 1 |
| d | d | c | 1 | b | a |

$Z_2L$ be the loop ring of the loop L over $Z_2$.

Consider the ideals

$I = \{0, 1 + a, b + c, b + d, c + d, a + b + c + 1, 1 + a + b + d, 1 + a + d + c\}$ and
$J = \{0, 1 + b, a + d, c + d, a + c, 1 + a + b + d, 1 + b + c + d, 1 + b + c + d, 1 + b + a + c\}$, K = {ideal generated by I and J}. It is easily verified that $I.J \neq \{0\}$, $J.K \neq \{0\}$ and $I.K \neq \{0\}$.

**THEOREM 2.4.14**: *Let $Z_2 = \{0, 1\}$ be the prime field of characteristic two and L be a finite loop of even order. Then the loop ring $Z_2 L$ contains non-zero orthogonal ideals.*

*Proof*: Since $|L| = 2n$, we have $I = \{0, 1 + m_2 + \ldots + m_n\}$ and $J = \{\alpha = \Sigma \alpha_i m_i / \Sigma \alpha_i = 0 \text{ and } m_i \in L\}$ clearly I and J are orthogonal ideals.



**THEOREM 2.4.15**: *Let L be a loop of odd order and $Z_2 = \{0, 1\}$ be the prime field of characteristic two, then the loop ring $Z_2L$ need not in general contain orthogonal ideals.*

*Proof*: Left for the reader to prove.

**THEOREM 2.4.16**: *Let L be a loop of order n with p/n and $Z_p$ be the odd prime field of characteristic p. Then the loop ring $Z_pL$ has nontrivial orthogonal ideals.*

*Proof*: Let $L = \{1, m_2, \ldots, m_n\}$ and $Z_p = \{0, 1, 2, \ldots, p - 1\}$, $p > 2$ and p/n. $Z_pL$ the loop ring of the loop L over $Z_p$. Consider the non-zero ideals given by $I = \{0, \gamma(1 + m_2 + \ldots + m_n) / \gamma = 1, 2, \ldots, p - 1\}$ and $J = \{\alpha = \Sigma \alpha_i m_i / \Sigma \alpha_i = 0\}$; clearly $I.J = \{0\}$.

**THEOREM 2.4.17**: *Let $L = \{1, m_2, \ldots, m_p\}$ and $Z_p = \{0, 1, \ldots, p - 1\}$. $Z_pL$ the loop ring of the loop L over $Z_p$. Clearly $I = \{0, \gamma(1 + m_2 + \ldots + m_p) / \gamma = 1, 2, \ldots, p - 1\}$ is a non-zero ideal of $Z_p L$ such that $I^2 = \{0\}$.*

The study of orthogonal ideals in loop rings in absent. So we have discussed to a possible extent the concept of orthogonal ideals. This notion will find its place in algebraic coding theory. Now we proceed on to define the concept of Smarandache orthogonal ideals.

**DEFINITION 2.4.9**: *Let KL be a loop ring of a loop L over the ring K. Let I and J be two ideals of KL. I and J are said to be Smarandache orthogonal ideals (S-orthogonal ideals) if the following conditions are satisfied.*

  i  *In $I.J = \{i, j / i \in I$ and $j \in J\}$ we have atleast a pair of elements $i \in I$ and $j \in J$ such that $i \neq 0$ and $j \neq 0$ with $i.j = 0$.*

  ii  *Every one of the pairs in I.J which are such that $i.j = 0$ are S-zero divisors.*

**THEOREM 2.4.18**: *Let KL be the loop ring of the loop L over the ring K. If every pair of ideals are orthogonal; then also KL need not have Smarandache orthogonal ideals.*

*Proof*: To prove this it is enough if the reader can show a pair $(x, y) \in I \times J$. I and J ideals of KL with $x \in I \setminus \{0\}$ and $y \in J \setminus \{0\}$ with $x.y = 0$ and $x.y$ is not a S-zero divisor.

**THEOREM 2.4.19**: *Let KL be the loop ring of the loop L over the field K. If KL has orthogonal ideals and no S-divisors of zero then KL has no S-orthogonal ideals.*

*Proof*: Left for the reader to prove.

Now we proceed on to define a new concept called right commutatively and its Smarandache analogue in loop rings.

**DEFINITION 2.4.10**: *Let RL be the loop ring of a loop L over the ring R with unity. If for every triple $\gamma, \beta, \alpha \in RL$ we have $(\alpha\beta) \gamma$ or $\alpha(\beta\gamma) = \alpha(\gamma\beta)$ or $(\alpha\gamma) \beta$, then we call*



*the loop ring RL to be a strongly right commutative, i.e., for every triple $\alpha, \beta, \gamma \in RL$ we must have*

a. $\alpha(\beta\gamma) = \alpha(\gamma\beta)$ or $(\alpha\gamma)\beta$.
b. $(\alpha\beta)\gamma = (\alpha\gamma)\beta$ or $\alpha(\gamma\beta)$.

*then we say RL is strongly right commutative. The term 'or' in general is not used in the mutually exclusive sense.*

**THEOREM 2.4.20**: *Every strongly right commutative loop ring RL is commutative.*

*Proof*: Obvious from the fact for every triple $\alpha, \beta, \gamma \in RL$ we have $\alpha\beta = \beta\alpha$, $\alpha\gamma = \gamma\alpha$ and $\beta\gamma = \gamma\beta$ where one of $\alpha$ or $\beta$ or $\gamma$ is assumed to be the unit of the loop ring RL. Hence every pair is commutative.

**THEOREM 2.4.21**: *Let RL be a strongly right commutative loop ring of the loop L over a ring R. Then the loop L is commutative.*

*Proof*: Obvious from the fact RL is commutative so $L \subset RL$ must also be commutative.

**DEFINITION 2.4.11**: *On similar lines one can define strongly left commutative loop rings as $(\alpha\beta)\gamma = (\beta\alpha)\gamma$ or $\beta(\alpha\gamma)$ for every $\alpha, \beta, \gamma \in L$.*

**DEFINITION 2.4.12**: *Let RL be the loop ring of the loop L over the ring R. If for every pair of elements $\alpha, \beta \in R$ there exists an element $\gamma \in R \setminus \{0, 1\}$ such that $\gamma(\alpha\beta)$ or $(\gamma\alpha) = (\gamma\beta)\alpha$ or $\gamma(\beta\alpha)$; then the loop ring RL is said to be right commutative.*

*Remark*: $\gamma(\alpha\beta) = \gamma(\beta\alpha)$ or $(\gamma\beta)\alpha$ or $(\gamma\alpha)\beta) = \gamma(\beta\alpha)$ or $(\gamma\beta)\alpha$. Similarly one can define left commutative loop rings.

**THEOREM 2.4.22**: *Let R be a commutative ring with unity and L a non-commutative loop. If the loop ring RL is a right commutative loop ring then RL has nontrivial divisors of zero.*

*Proof*: Follows from the fact $\gamma(\alpha\beta) = \gamma(\beta\alpha)$ (when $\gamma \neq 0$ and $\gamma \neq 1$), we have $\gamma(\alpha\beta - \beta\alpha) = 0$. Hence the claim.

**THEOREM 2.4.23**: *Let RL be the loop ring of the loop L over the ring R. If L is a strongly right commutative ring then RL is a right commutative ring.*

*Proof*: Straightforward by the very definition.

**THEOREM 2.4.24**: *Let RL be a right commutative loop ring having no zero divisors. Then RL is commutative.*

*Proof*: Follows from the fact $\gamma(\alpha\beta) = \gamma(\beta\alpha)$; $\gamma \neq 0$, or $1$ $\gamma(\alpha\beta - \gamma\alpha) = 0$. Since RL has no divisors of zero we have $\alpha\beta = \beta\alpha$.



**THEOREM 2.4.25**: *Let RL be a right commutative loop ring then RL is a weakly right commutative ring.*

*Proof*: Straightforward.

**DEFINITION 2.4.13**: *Let RL be the loop ring of the loop L over R. RL is said to be strictly right commutative if for every pair of elements $\alpha, \beta \in RL$ we have a $\gamma \in RL \setminus \{0, 1\}$ such that $\gamma(\alpha\beta) = \gamma(\beta\alpha)$ or $(\gamma\alpha)\beta = (\gamma\beta)\alpha$.*

**THEOREM 2.4.26**: *Every strictly right commutative loop ring is right commutative.*

*Proof*: Straightforward by the definition.

**DEFINITION 2.4.14**: *Let RL be the loop ring of the loop L over R. We say RL is a Smarandache strongly right commutative (S-strongly right commutative) if for every S-subring A of RL; for all $\alpha, \beta, \gamma \in A$ we must have*

$$\alpha(\beta\gamma) = \alpha(\gamma\beta) \text{ or } (\alpha\gamma)\beta \text{ or}$$
$$(\alpha\beta)\gamma = (\alpha\gamma)\beta \text{ or } \alpha(\gamma\beta)$$

*then we say RL is Smarandache strongly right commutative loop ring (S-strongly right commutative loop ring).*

*The term 'or' in general is not used in the mutually exclusive sense. If in a loop ring RL we have atleast one S-subloop A which is strongly right commutative then we call the loop ring RL to be a Smarandache-right commutative (S-right commutative).*

**THEOREM 2.4.27**: *If the loop ring RL is S-strongly right commutative then RL is S-right commutative.*

*Proof*: Direct from the very definition.

We proceed on to define inner commutative loops and its Smarandache analogue.

**DEFINITION 2.4.15**: *Let RL be the loop ring of a loop L over R. We say the loop ring RL is inner commutative if every subring of RL is commutative.*

**THEOREM 2.4.28**: *Let R be a field or a commutative ring with 1. $L_p(m) \in L_p$ be a loop where p is an odd prime. The loop ring $RL_p(m)$ is inner commutative.*

*Proof*: Follows from the fact $L_p(m)$ has only subloops, which are subgroups of order 2. Hence $RL_p(m)$ is inner commutative.

**THEOREM 2.4.29**: *Let R be a field and L be an inner commutative loop. Then the loop ring RL is an inner commutative loop ring.*

*Proof*: Straightforward from the definitions.

**THEOREM 2.4.30**: *Let RL be a commutative loop ring then RL is an inner commutative loop ring.*



*Proof*: Follows from the very definitions.

All loop rings which are inner commutative need not be commutative. In view of this we propose the following.

**THEOREM 2.4.31**: *An inner commutative loop ring RL in general need not be commutative.*

*Proof*: By an example $ZL_{11}(9)$ be the loop ring of the loop $L_{11}(9)$ over Z. Clearly $ZL_{11}(9)$ is an inner commutative loop ring which is not commutative.

Now we proceed on to define Smarandache-analogues.

**DEFINITION 2.4.16**: *Let RL be the loop ring of the loop L over R. We say RL is a Smarandache inner commutative loop ring (S-inner commutative loop ring) if every S-subloop of RL is commutative.*

*Example 2.4.4*: Let $ZL_{13}(5)$ be the loop ring of the loop $L_{13}(5)$ over Z. $Z_{13}(5)$ is a S-inner commutative loop ring.

**DEFINITION 2.4.17:** *Let RL be the loop ring for ab = ba in RL if we have $(a\alpha)b = b(\alpha a)$ (or $= b(\alpha a)$) for all $\alpha \in RL$, then we say the pair a, b is pseudo commutative. [or $a(\alpha b) = (b\alpha)a$ or $b(\alpha a)$] i.e., $(a\alpha)b$ or $a(\alpha b)$ may be taken as the initial point, for both may or may not be equal.*

*We call the loop ring to be a pseudo commutative loop if every pair a, b in L is pseudo commutative.*

**THEOREM 2.4.32**: *Let RL be a commutative loop ring then RL, is a pseudo commutative loop ring.*

*Proof*: Obvious by the very definition.

**DEFINITION 2.4.18**: *Let RL be a loop ring. We say RL is a Smarandache pseudo commutative loop ring (S-pseudo commutative loop ring) if RL has atleast one S-subring A which is pseudo commutative.*

*We say RL is Smarandache strongly pseudo commutative (S-strongly pseudo commutative) if every S-subring A or RL is pseudo commutative.*

**THEOREM 2.4.33**: *Let RL be a S-strongly pseudo commutative loop ring then RL is S-pseudo commutative.*

*Proof*: Straightforward by the very definition.

**THEOREM 2.4.34**: *Let RL be a S-strongly pseudo commutative loop ring; then RL in general need not be a pseudo commutative loop ring.*



*Proof*: Consider the loop ring $Z_7L_{19}(5)$. Clearly $Z_7L_{19}(5)$ is a S-strongly pseudo commutative loop ring.

But $Z_7L_{19}(5)$ is not a commutative loop ring. Further $Z_7L_{19}(5)$ is not a pseudo commutative loop ring.

**DEFINITION 2.4.19**: *Let RL be the loop ring of the loop L over the ring R. The pseudo commutator of RL denoted by $P(RL) = \langle\{p \in RL / a(\alpha b) = p(b\alpha)a; a, b, \alpha \in RL\}\rangle$; where $\langle\rangle$ denotes the subring generated by all p's.*

Now one of the natural question would be if L is a loop having P (L) to be the pseudo commutator; find the relation between P (RL) and R (P (L)).

**DEFINITION 2.4.20**: *Let RL be a loop ring. The Smarandache pseudo commutator (S-pseudo commutator) subring of RL are defined and denoted by $S(P(RL)) = \langle\{p \in RL / a(\alpha b) = p(b\alpha)a$ where $a, b, \alpha \in A$; A a S-subring of RL$\}\rangle$. The subring generated by p. We may have several or one S-pseudo commutator subring of RL depending on varying S-subrings A of RL. It may happen even for varying S-subrings A of RL we may have the same S (P (RL)).*

The reader is expected to find interesting and innovative results in this direction. The major feature about substructures of any algebraic structure is that how does the lattice of the collection of all substructures look like.

For, in case of groups we see the normal subgroups forms a modular lattice the two-sided ideals of a ring form a modular lattice.

Now we study what is the lattice structure of the following substructures.

   a. subrings of a loop ring.
   b. S-subrings of a loop rings.
   c. right (left) ideals of a loop ring.
   d. ideals of a loop ring.
   e. S-right (S-left) ideals of loop ring.
   f. S-ideals of a loop ring.

In our opinion such study for loop rings has not been carried out. Some times the collection of the substructures may be a distributive or a modular lattice.

The more important factor is sometimes the substructures may be a supermodular lattice. (we say a lattice L is super modular if for all x, y, z, a ∈ L. We have $(a \cup x) \cap (a \cup y) \cap (a \cup z) = a \cup \{(x \cap y) \cap (a \cup z)\} \cup \{(x \cap z) \cap (a \cup y)\} \cup \{(y \cap z) \cap (a \cup x)\}$ this will be called as the supermodular identity for simplicity of notation we replace '$\cup$' by '+' and '$\cap$' by '.'; thus the identity reads as $(a + x)(a + y)(a + z) = a + (xy)(a + z) + xz(a + y) + yz(a + x)$ for all a, x, y, z in L)

*Example 2.4.5*: Let $Z_2 = \{0, 1\}$ be the prime field of characteristic two and L be a loop given by the following table:



| . | 1 | a | b | c | d |
|---|---|---|---|---|---|
| 1 | 1 | a | b | c | d |
| a | a | 1 | c | d | b |
| b | b | d | a | 1 | c |
| c | c | b | d | a | 1 |
| d | d | c | 1 | b | a |

The only ideals of $Z_2L$, the loop ring of the loop L over $Z_2$ are

$Z_2L$, $\{0\}$, $I_1 = \{0, (1 + a + b + c + d)\}$, $I_2 = \{\alpha = \Sigma \alpha_i m_i / \Sigma \alpha_i = 0\}$. The set of ideals are $\{Z_2L, \{0\}, I_1, I_2\}$.

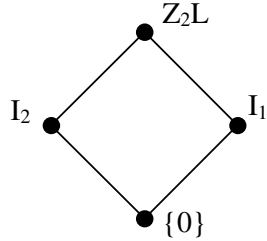

**Figure 2.4.1**

is a modular lattice. In fact we can see that this set also forms a supermodular lattice.

Now we can characterize when will the ideals of a loop ring does not form a super modular lattice. We denote by $A \cup B = A + B$, mean only the ideal generated by the set $A \cup B$.

**THEOREM 2.4.35**: *Let L be a finite loop and K any finite field, KL be the loop ring. S denotes the collection of all ideals of KL. S is a super modular lattice if and only if S does not contain a set of 4 distinct ideals say A, B, C, D such that $A \cup B = A \cup C = A \cup D > A$ and $A > B \cap C$, $A > B \cap D$, $A > D \cap C$.*

*Proof*: By the symbol '>' we mean $A \cup B$ has more elements than A. Similarly $A > B \cap C$ implies A has more elements than $B \cap C$ i.e., not only $B \cap C$, it is properly contained in A.

Suppose S contains a set of 4 distinct elements A, B, C, D such that $A \cup B = A \cup C = A \cup D > A$ and $A > B \cap C$, $A > B \cap D$, $A > D \cap C$ then in the equation.

$(A \cup B) \cap (A \cup C) \cap (A \cup D) = A \cup [(B \cap C) \cap (A \cup D)] \cup [(D \cap B) \cap (A \cup C)] \cup [(C \cap D) \cap (A \cup B)]$ we get right hand side of the equation to be $A \cup B$ where as the left hand side is A. Clearly $A \cup B \neq A$ as A and B distinct. Hence S is not a super modular lattice.

Conversely suppose S is not supermodular to prove S contains a set of 4 distinct ideals A, B, C, D such that $A \cup B = A \cup C = A \cup D > A$ and $A > B \cap C$, $A > B \cap D$, $A > D \cap C$. If S is not supermodular we have $(A \cup B) \cap (A \cup C) \cap (A \cup D) > A \cup \{(B \cap C) \cap (A \cup D)\} \cup \{(D \cap B) \cap (A \cup C)\} \cup \{(D \cap C) \cap (A \cup B)\}$.



Put a = a ∪ [(B ∩ C) ∩ (A ∪ D)] ∪ [(B ∩ D) ∩ (A ∪ C)] ∪ {(D ∩ C) ∩ (A ∪ B)] and b = B, c = C, and d = D. a ∪ b = A ∪ B, a ∪ C = A ∪ C, a ∪ d = A ∪ D. So (a ∪ b) ∩ (a ∪ c) ∩ (a ∪ d) > a. It is easily verified a > b ∩ c, a > c ∩ d and a > b ∩ d. Hence the claim.

Now we proceed on to see what is structure of strongly modular lattice. We at the outset say in a lattice L '+' denotes '∪' and '.' denotes '∩'.

**DEFINITION 2.4.21**: *Let L be a lattice we say an element a ∈ L is a strongly modular element of L if (a + b) (a + c) (a + d) (a + e) = a + bc (a + d) (a + e) + bd (a + c) (a + e) + be (a + b) (a + c) + cd ( a + b) (a + e) + ce (a + b) (a + d) + de (a + b) (a + c) for all b, c, d, e ∈ L.*

*This identity is known as the strongly modular identity. If in a lattice L all elements satisfy the strongly modular identity then we say L is a strongly modular lattice.*

Now our interest is to see when do the collection of all ideals of a loop ring form a strongly modular lattice.

*Example 2.4.6*: Let $Z_2$ = {0, 1} be the prime field of characteristic two and L be a loop given by the following table:

| . | 1 | a | b | c | d | e | f |
|---|---|---|---|---|---|---|---|
| 1 | 1 | a | b | c | d | e | f |
| a | a | f | c | b | e | 1 | d |
| b | b | d | f | 1 | a | c | e |
| c | c | e | d | a | f | b | 1 |
| d | d | 1 | e | f | c | a | b |
| e | e | c | 1 | d | b | f | a |
| f | f | b | a | e | 1 | d | c |

$Z_2L$ be the loop ring of the loop L over $Z_2$. Ideals of $Z_2L$ are {0} = $I_0$, $I_1$ = {0, 1 + a + b + c + d + e + f}, $I_2$ = {0, Σ $α_i$ $m_i$ / Σ $α_i$ = 0} and $Z_2L$. The collection of ideals be denoted by S = {$I_0$, $I_1$, $I_2$, $Z_2L$}.

The lattice diagram of S is

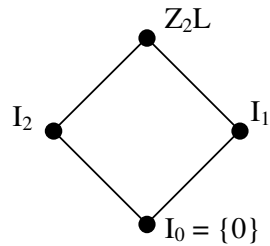

**Figure 2.4.2**

which is easily verified to be a strongly modular lattice.



**THEOREM 2.4.36**: *Let L be a finite loop and K any field. KL the loop ring of the loop L over the field K. S denotes the set of all ideals of KL. S is strongly modular if and only if there does not exists 5 distinct elements A, B, C, D, E in S such that $A \cup B = A \cup C = A \cup D = A \cup E > A$ and $A > B \cap C$, $A > C \cap E$, $A > D \cap E$, $A > B \cap D$, $A > B \cap E$ and $A > C \cap D$.*

*Proof*: The proof follows by using the strongly modular identity and the mentioned relations in the theorem.

Now in case of S-ideals and S-subrings the results given in those theorems hold good. But it is pertinent to mention here that we are not able to find other than the two theorems mentioned in this section about the S-ideals or ideal or subrings or of S-subrings.

**PROBLEMS:**

1. Is the loop ring $Z_7L_9(8)$ a Marot loop ring? Justify your claim.
2. Find all the regular elements of the loop ring $Z_2L_{11}(3)$.
3. Give an example of a loop ring which is a weak Marot loop ring.
4. Find $L^*$ for the loop ring $Z_{12}L_{13}(3)$.
5. Does the loop ring $Z_3L_{15}(8)$ have orthogonal ideals?
6. Can $Z_3L_{15}(8)$ have S-orthogonal ideals?
7. Give an example of a strongly right commutative loop ring.
8. Is $Z_8L_9(5)$ inner commutative? Is it S-inner commutative?
9. Give an example of a pseudo commutative loop ring.
10. Is the loop ring $Z_2L_{19}(3)$ S-pseudo commutative?
11. Find the lattice of S-subring of $Z_4L_5(3)$.
12. Find all S-ideals of $Z_8L_7(3)$. Does it form a modular lattice?
13. Find all ideals of $Z_6L_7(5)$. Do they form a supermodular lattice?
14. Give an example of a loop whose S-ideal forms a strongly modular lattice.
15. Give an example of a loop ring where only S-ideals form a modular lattice but ideals do not form a modular lattice.

## 2.5 General properties of Smarandache loop rings and loop rings

The determination of group of units in a group ring is an interesting one to several researchers. E.G. Goodaire [24 to 29] has studied this problem for loop rings in 1986. The isomorphism problem for group rings was proposed in [32]. He asks if isomorphism of integral group rings ZH ≈ ZG implies that the groups G and H are isomorphic. The problem has been settled only for special cases and remains open. The famous semi simplicity problem of finding necessary and sufficient condition for J(RZ) = 0 which remains still unsettled even over field of characteristic zero this has been studied by [24 to 29] in case of alternative loop rings. Very recently it has been proved that every representation of a semi simple analytic Moufang loop is completely reducible. Here we derive the Smarandache analogue of all these and suggest some problems.



**DEFINITION 2.5.1:** *An alternative ring is one in which the associator $(x,y,z) = (xy)z - x(yz)$ is a skew symmetric function of its arguments.*

It was proved that if in an alternative ring the subring generated by any pair of elements is associative.

**THEOREM [29]:** *Suppose ZL is an alternative ring which is not associative and G is any subgroup of index 2 in L. Then*

    (i)    *The map $* : G \to G$ defined by*

$$g^* = \begin{cases} g & g \in Z(G) \\ eg & g \notin Z(G) \end{cases}$$

        *is an involution on G which extends to an involution on the group ring ZG as follows.*

$$\left( \sum_{g \in G} \alpha_s g \right)^* = \sum \alpha_s g^*.$$

    (ii)    *Every element in ZL can be written in the form $x + yu$ where $x, y \in ZG$. Multiplication in ZL is given by $(x + yu)(z + wu) = (xz + g_o w^* y) + (wx + yz^*)u$ where $x, y, z, w \in ZG$ and $g_o$ is an element in $Z(G)$.*

    (iii)    *The map $* : ZL \to ZL$ defined by $(x + yu)* = x^* + eyu$ is an involution on ZL extending, the one on ZG given in (i). An element $\alpha \in ZL$ is central if and only if $\alpha^* = \alpha$.*

    (iv)    *The centre of the nucleolus of the ring ZL coincides and are equal to $\{x + yu \mid x, y\ Z(RG); ey = y\} = \{x + yu \mid x \in Z(RG)$ and $ey = y\}$.*

The very natural study in the Smarandache direction for the above theorem is for the loop ring ZL where ZL in an alternative ring find Smarandache centre and Smarandache nucleus of the ring ZL.

Do they concide? We see in case of the loops in $L_n$.

**THEOREM 2.5.1:** *Let $L_n(m) \in L_n$, The centre $Z(L_n(m))$ and Nucleus $N(L_n(m))$ are just $\{e\}$.*

*Proof:* Left as an exercise for the reader to prove.

The study of the nucleolus of $RL_n(m)$ and centre of $R L_n(m)$ where R is a ring is an important study. Do they coincide? or Are they just equal to $R \{e\} = R$? and so on.

Find conditions on the units in $ZL_n(m)$ so that they are central units in $ZL_n(m)$. Can S-units be central units in $ZL_n(m)$?



If has been proved by [24].

**THEOREM 2.5.2:** *Suppose r is a central unit is an alternative loop ring ZL such that $\bar{r}$ is trivial. Then ZL is a S-loop ring.*

*Proof:* Given r is central, ey = y. since $\bar{r}$ is trivial either $x = \pm g + (1 - e) x_1$, $y = (1 - e) y_1$ for some $x_1$ and $y_1$ in the group ring [using results from 26] we have ZL is a S-loop ring.

It is interesting to note the Zassenhaus conjecture for group rings which is still open in group rings is also open in loop rings all the more in case of Smarandache loop rings and loop rings which are S-rings.

**Zassenhaus conjecture**: Let $r \in ZG$ be a normalized unit of finite order. Then there exists an invertible element $\alpha \in QG$ and an element $g \in G$ such that $\alpha^{-1} r \alpha = g$.

Here Z denotes the ring of rational integers and Q the field of rationals $\in : ZG \rightarrow Z$ the augmentation map; i.e. map given by $\in (\Sigma \alpha_g g) = \Sigma \alpha_g$. If we denote by U(ZG) the group of units of ZG, then the set $V (ZG) = \{\alpha \in U (ZG) \mid \in (\alpha) = 1\}$ is called the group of normalized units of ZG.

The elements of the form $\pm g$, $g \in G$ are trivial units of ZG. [26] has established the alternative analogue of the Zassenhaus conjecture as follows;

**THEOREM [25]:** *Let $\tau$ be a normalized torsion unit in the integral alternative loop ring ZL of a finite loop L. Then there exists units $\gamma_1$, $\gamma_2 \in Q L$ and $m \in L$ such that $\gamma_2^{-1} (\gamma_1^{-1} r \gamma_1) \gamma_2 = m$.*

He has proved that this extended conjecture is at least true in case of these alternative loop rings, which are not associative. Now we propose a analogue of Zassenhaus conjecture for S-loop rings and SNA-rings. Let KL be a loop ring, which is a SNA-ring. Then the augmentation map $\in : K L \rightarrow K$ where K is the ring of integers given by $\in (\sum \alpha_g g) = \sum \alpha_g$.

If S U (ZG) denotes the set of S-units of ZG then the set $S V (ZG) = \{\alpha \in S U (ZG) \in (\alpha) = 1\}$ is called the group of Smarandache normalized units (S-normalized units) of ZG. Smarandache Zassenhaus conjecture. Let $r \in Z L$ which is Smarandache normalized unit in the integral loop ring of a finite loop L. Then there exists $\alpha_1$, $\alpha_2 \in QL$ and $m \in L$ such that $\alpha_2^{-1} (\alpha_1^{-1} r \alpha_1) \alpha_2 = m$.

It is still an open problem to solve the Smarandache Zassenhous conjecture for

    I. Smarandache loop rings.
    II. S- loop rings.



**PROBLEMS:**

1. Is $ZL_5(3) \cong ZL_5(2)$?
2. Can $ZL_5(3)$ be S-Isomorphic with $Z_7(3)$? Justify your answer
3. Can we say $QL_n(m)$ is S-isomorphic with $Z_2L_n(m)$? Justify your claim.
4. Find all S-units in $Z_5 L_9(8)$.
5. Can $Z_5 L_9(8)$ have central units?
6. Verify Zassenhaus conjecture for the loop ring $ZL_{11}(3)$.
7. Will $ZL_{11}(3)$ satisfy S-Zassenhaus conjecture? Justify.



**Chapter 3**

# GROUPOID RINGS AND SMARANDACHE GROUPOID RINGS

This chapter has 5 sections. In the first section we introduce the concept of groupoids and Smarandache groupoids. Section two is devoted to the study and introduction of groupoid rings and Smarandache groupoid rings. In section three we study special elements in Smarandache groupoid rings and SNA-groupoid rings. Study of substructure is carried out in section four. In the final section several new properties like identities satisfied by Smarandache groupoid rings and other special and new notions about these groupoid rings are introduced and studied, we define level II Smarandache groupoid rings and study them.

## 3.1 Groupoids and Smarandache groupoids

In this section we introduce several of the properties of groupoids and Smarandache groupoids. Clearly as groupoids are the most generalized algebraic structures viz loops and semigroups; all loops are groupoids and groupoids in general are not loops. Likewise all semigroups are groupoids but all groupoids in general are not semigroups. Thus we study these groupoids and Smarandache groupoids with a main motivation to construct a new class of non-associative rings called groupoid rings. Groupoid rings are a generalized class of loop rings. Thus this study of these non-associative rings and mainly their Smarandache NA rings will reveal a lot of interesting results. The study of Smarandache groupoids started in the year 2002 [70]. Now we introduce and study Smarandache groupoid rings and groupoid rings.

**DEFINITION 3.1.1:** *A non-empty set of elements G is said to form a groupoid if in G is defined a binary operation called the product denoted by $*$ such that $a * b \in G$ for all $a, b \in G$. It is important to mention here that the binary operation $*$ defined on the set G need not be associative i.e. $(a * b) * c \neq a * (b * c)$ in general for all $a, b, c \in G$.*

Throughout this chapter by a groupoid (G, $*$) we assume that G is non-associative under the operation $*$.

*A groupoid G is said to be a commutative groupoid if for every $a, b \in G$ we have $a * b = b * a$. A groupoid G is said to have an identity element e in G if $a * e = e * a = a$ for all $a \in G$.*

**DEFINITION 3.1.2:** *Let (G, $*$) be a groupoid a proper subset $H \subset G$ is a subgroupoid if (H, $*$) is itself a groupoid.*

*Example 3.1.1:* Let R be the set of reals (R, –) is a groupoid, where '–' is the usual subtraction of reals in R.



***Example 3.1.2:*** Let G be a groupoid given by the following table:

| * | $a_1$ | $a_2$ | $a_3$ | $a_4$ |
|---|---|---|---|---|
| $a_1$ | $a_1$ | $a_3$ | $a_1$ | $a_3$ |
| $a_2$ | $a_4$ | $a_2$ | $a_4$ | $a_2$ |
| $a_3$ | $a_3$ | $a_1$ | $a_3$ | $a_1$ |
| $a_4$ | $a_2$ | $a_4$ | $a_2$ | $a_4$ |

This has $H = \{a_1, a_3\}$ and $K = \{a_2, a_4\}$ to be subgroupoids of the groupoid (G, *). A groupoid G which has only a finite number of elements in them is called a finite groupoid. A groupoid which has infinite number of elements is called an infinite groupoid or a groupoid of infinite order.

**DEFINITION 3.1.3:** *A groupoid G is said to be a Moufang groupoid if it satisfies the Moufang identity $(xy)(zx) = (x(yz))x$ for all x, y, z in G.*

**DEFINITION 3.1.4:** *A groupoid G is Bol if it satisfies the Bol identity $((xy)z)y = x((yz)y)$ for all x, y, z $\in$ G. A groupoid is a P-groupoid if $(xy)x = x(yx)$ for all x, y $\in$ G. A groupoid G is right alternative if it satisfies the identity $(xy)y = x(yy)$ for all x, y $\in$ G and G is left alternative if $(xx)y = x(xy)$ for all x, y in G. A groupoid G is alternative if it is both right and left alternative. A groupoid is said to be a Jordan groupoid if it satisfies the Jordan identity $x(x^2 y) = x^2 (xy)$ for all x, y $\in$ G.*

**DEFINITION 3.1.5:** *G be a groupoid, P a non-empty proper subset of G. P is said be a left ideal of G if*

   *(a)   P is a subgroupoid of G.*

   *(b)   For all x $\in$ G and a $\in$ P we have x a $\in$ P. One can similarly define a right ideal of the groupoid G. P is called an ideal if P is simultaneously a left and a right ideal of G.*

**DEFINITION 3.1.6:** *Let G be a groupoid. A subgroupoid V of G is said to be a normal subgroupoid of G if*

   *a.   aV = Va*
   *b.   (Vx)y = V(xy)*
   *c.   y(xV) = (yx) V*

*for all x, y, a $\in$ G. A groupoid G is said to be simple if it has no nontrivial normal subgroupoid.*

**DEFINITION 3.1.7:** *A groupoid G is normal if*

   *1.   x G = G x*
   *2.   G(xy) = (Gx)y*
   *3.   y(xG) = (yx)G*

*for all x, y $\in$ G.*



**DEFINITION 3.1.8:** *Let G be a groupoid H and K be two proper subgroupoids of G with H ∩ K = φ or {e}. We say H is conjugate with K if there exists x ∈ H such that H = xK or Kx 'or' in the mutually exclusive sense.*

Direct products are defined as in case of other algebraic structures.

A groupoid G has a zero divisor if 0 ∈ G and for same 0 ≠ a ∈ G, there exists b ≠ 0 ∈ G with a ∗ b = 0. An element a ∈ G is an idempotent of G if a ∗ a = $a^2$ = a. Let G be a groupoid. The center of the groupoid C(G) = {x ∈ G / a x = x a for all a ∈ G}. We say a, b ∈ G is a conjugate pair if a = bx (or xb for some x ∈ G) and b = ay (or ya for some y ∈ G). An element a in G is said to be right conjugate to b in G if we can find x, y ∈ G such that a. x = b and b. y = a (x. a = b and y. b = a). Similarly we define left conjugate elements of G. One of the drawbacks about groupoids was that we did not have natural examples of groupoids built using $Z_n$. So we recall the definition of a new class of groupoids using $Z_n$. [70]

**DEFINITION 3.1.9:** *Let $Z_n$ = {0, 1, 2, …., n-1}; n ≥ 3, n < α. Define a binary operation ∗ on $Z_n$ as follows. For any a, b ∈ $Z_n$ define a ∗ b = (at + bu) (mod n) where (t, u) = 1 where '+' is the usual addition modulo n. {$Z_n$, (t, u), ∗} or $Z_n$ (t, u) is a groupoid. For varying values of t, u ∈ $Z_n$ \ {0} we get a class of groupoids for a fixed n. This class of groupoids is denoted by Z(n) and all groupoids in this class is of order n. i.e. Z(n) = { $Z_n$, (t, u), ∗ | for integers t, u ∈ $Z_n$ \ {0} with (t, u) = 1}.*

Several interesting results can be obtained but the reader is requested to refer [W B V Pad B].

By $Z^*(n)$ we denote the class of groupoids of order n built using $Z_n$ with the binary operation '.' such that for a, b ∈ $Z_n$, a.b = at + bu (mod n) where (t, u) need not always be relatively prime but t ≠ u, t, u ∈ $Z_n$ \ {0}. Clearly the class of groupoids Z(n) is completely contained in $Z^*(n)$. Similarly we define $Z^{**}(n)$ to be the class of groupoids built using $Z_n$ under the binary operation 'o' defined by a o b = ta + ub (mod n) for all a, b ∈ $Z_n$ with t, u ∈ $Z_n$ \ {0}. Thus for varying t and u we get a class of groupoids for any integer n. We denote this class by $Z^{**}(n)$. We have Z(n) ⊆ $Z^*(n)$ ⊆ $Z^{**}(n)$.

Now the most generalized class of groupoids using $Z_n$ denoted by $Z^{***}(n)$ is built using t, u ∈ $Z_n$ with the binary operation '∗' such that a ∗ b = ta + ub (mod n). Z(n) ⊂ $Z^*(n)$ ⊆ $Z^{**}(n)$ ⊆ $Z^{***}(n)$. Thus we get a natural class of groupoids of order n.

Several interesting properties about them can be had from (WBV) P and B).

Now we proceed on to define Smarandache groupoids and study about them. Smarandache groupoids were introduced only in the 2001.

**DEFINITION 3.1.10:** *A Smarandache groupoid (S-groupoid) G is a groupoid which has a proper subset S, S ⊂ G such that S under the operations of G is a semigroup.*

*Example 3.1.3:* The groupoid (G, ∗) given by the following table:



| * | 0 | 1 | 2 | 3 | 4 | 5 |
|---|---|---|---|---|---|---|
| 0 | 0 | 3 | 0 | 3 | 0 | 3 |
| 1 | 1 | 4 | 1 | 4 | 1 | 4 |
| 2 | 2 | 5 | 2 | 5 | 2 | 5 |
| 3 | 3 | 0 | 3 | 0 | 3 | 0 |
| 4 | 4 | 1 | 4 | 1 | 4 | 1 |
| 5 | 5 | 2 | 5 | 2 | 5 | 2 |

Clearly $S_1 = \{0, 3\}$ and $S_2 = \{1, 4\}$ are semigroups of $(G, *)$. A S-groupoid $(G, *)$ is said to be a Smarandache commutative groupoid (S-commutative groupoid) if there is a proper subset which is a commutative semigroup.

*Let $(G, *)$ be a groupoid. A non-empty subset H of G is said to be a Smarandache subgroupoid (S-subgroupoid) if H contains a proper subset $K \subset H$ such that K is a semigroup under the operations of *.*

**THEOREM 3.1.1:** *If G has a S-subgroupoid then G is a itself a S-groupoid.*

*Proof:* Direct by the very definition.

**THEOREM 3.1.2:** *Every subgroupoid of a S-groupoid need not in general be a S-subgroupoid.*

*Proof:* Follows by counter examples.

**DEFINITION 3.1.11:** *A Smarandache left ideal (S-left ideal) A of a S-groupoid G satisfies the following conditions.*

1. *A is a S-subgroupoid of G.*
2. *$x \in G$ and $a \in A$ then $x a \in A$.*

*Similarly, we can define Smarandache right ideal (S-right ideal), if A is both a S-right ideal and S-left ideal simultaneously then we say A is a Smarandache ideal (S-ideal).*

*Let V be a S-subgroupoid of a groupoid G. We say V is a Smarandache seminormal groupoid (S-seminormal groupoid) if*

(i)   $aV = X$ for all $a \in G$
(ii)  $Va = Y$ for all $a \in G$

*where either X or Y is a S-subgroupoid of G but both X and Y are subgroupoids of G.*

*We say V is a Smarandache normal groupoid (S-normal groupoid) if $aV = X$ and $Va = Y$ for all $a \in G$ where both X and Y are S-subgroupoids of G. Let H and P be any two subgroupoids of a groupoid G we say H and P are Smarandache semi conjugate subgroupoids (S-semi conjugate subgroupoids) of G, if*



1. H and P are S-subgroupoids of G.
2. H = xP or Px or
3. P = x H or Hx for some x ∈ G.

**DEFINITION 3.1.12:** *Let H and P be two subgroupoids of G. We say H and P are Smarandache conjugate subgroupoids (S-conjugate subgroupoids) of G if*

1. H and P are S-subgroupoids of G.
2. H = xP or Px and
3. P = xH or Hx.

**DEFINITION 3.1.13:** *Let G be a S-subgroupoid. We say G is a Smarandache inner commutative groupoid (S-inner commutative groupoid) if every semigroup contained in every S-subgroupoid of G is commutative.*

**DEFINITION 3.1.14:** *Let (G, ∗) be a S-groupoids. If every S-subgroupoid H of G satisfies the Moufang identity for all x, y, z in H then G is said to be a Smarandache strong Moufang groupoid (S-strong Moufang groupoid). If there exists at least one S-subgroupoid H of G which satisfies the Moufang identity then we call G a Smarandache Moufang groupoid (S-Moufang groupoid).*

It is easily verified that all S-strong Moufang groupoids are S-Moufang groupoid, whereas S-Moufang groupoids are not in general S-strong Moufang groupoids.

On similar lines one can define Smarandache Bol groupoid, Smarandache strong Bol groupoid, Smarandache alternative (right/left) groupoid, Smarandache strong alternative (right/left) groupoid and Smarandache strong P-groupoid and Smarandache P-groupoid. For more about groupoids and S-groupoids please refer [4, 5, 67, 70].

**PROBLEMS:**

1. Find an example of a groupoid of order 7, which is not a S-groupoid.
2. Is $Z_7(3,2)$ a S-groupoid?
3. Give an example of a S-normal groupoid.
4. Give an example of a normal groupoid which is not a S-groupoid.
5. Does there exist a S-strong Bol groupoid of order 10?
6. Give an example of a Bol groupoid, which is not a S-Bol groupoid.
7. Will a S-Moufang groupoid be a Moufang groupoid?
8. Show by an example that a Bol groupoid need not be a S-Bol groupoid.
9. Find a S-groupoid which is simple and of order 8.
10. Show by example all simple groupoids need not be S-simple groupoids.

## 3.2 Groupoid rings and Smarandache groupoid rings

In this section we introduce the concept of groupoid rings and Smarandache groupoid rings. The study of groupoid rings is very recent. Groupoid rings are constructed analogous to group rings where groups are replaced by groupoids. Rings



are always assumed to be commutative with unit. Further groupoid rings leads to a collection of non-associative rings.

These non-associative rings may or may not contain identity. Using different types of groupoids we can have non-associative rings which are commutative or non-commutative ring of infinite or finite order. Thus using groupoid rings we get a new class of non-associative rings which are a generalized class of loop rings for every loop ring is a groupoid ring. Study of loop rings are carried out in chapter 2. Throughout this section, G will denote a groupoid may or may not contain 1. R is a commutative ring with 1 or a field.

**DEFINITION 3.2.1:** *Let G be a groupoid and R a ring or a field. The groupoid ring RG consists of all finite formal sums of the form $\Sigma \, r_i \, g_i$ (i – running over a finite number) where $r_i \in R$ and $g_i \in G$ satisfying the following conditions.*

1. $\sum_{l=1}^{n} r_i \, g_i = \sum_{l=1}^{n} s_i \, g_i \Leftrightarrow r_i = s_i$ *for i = 1, 2, 3,...., n.*

2. $\left( \sum_{i=1}^{n} \alpha_i \, m_i \right) + \left( \sum_{i=1}^{n} \beta_i \, m_i \right) = \sum_{i=1}^{n} (\alpha_i + \beta_i) \, m_i$

3. $\left( \sum_{i=1}^{n} (\alpha_i \, m_i) \right) \left( \sum_{i=1}^{n} (\beta_i \, m_i) \right) = \sum \gamma_k m_k$ *where* $m_k = m_i m_j$ *where* $\gamma_k = \Sigma \, \alpha_i \, \beta_j$

4. $r_i \, m_i = m_i \, r_i$ *for all* $r_i \in R$ *and* $m_i \in G$.

5. $r \sum_{l=1}^{n} r_i \, m_i = \sum_{l=1}^{n} (r \, r_i) \, m_i$ *for* $r_i, r \in R$. *Since* $I \in R$ *and* $m_i \in G$, *we have* $G = 1$. $G \subseteq RG$ *and* $R \subseteq RG$ *if and only if G has identity otherwise* $R \not\subseteq RG$.

*Clearly RG is a non-associative ring with $0 \in R$ as its additive identity. The groupoid ring RG is an alternative ring if $(xx) y = x(xy)$ and $x(yy) = (xy)y$ for all $x, y \in RG$.*

*Let RG be a ring and let $x, y \in RG$. A binary operation known as the circle operation denoted by 'o' is defined by $x \circ y = x + y - xy$. The concepts of right (left) quasi regular elements, the Jacobson radical J(RG) of RG are defined as in the case of any non-associative ring.*

***Example 3.2.1:*** Let $Z_2 = \{0,1\}$ be the prime field of characteristic two and G be a groupoid given by the following table:

| . | $a_0$ | $a_1$ | $a_2$ |
|---|---|---|---|
| $a_0$ | $a_0$ | $a_2$ | $a_1$ |
| $a_1$ | $a_1$ | $a_0$ | $a_2$ |
| $a_2$ | $a_2$ | $a_1$ | $a_0$ |



$Z_2 G$ is the groupoid ring given by $Z_2 G = \{0, a_0, a_1, a_2, a_0 + a_1, a_0 + a_2, a_1 + a_2, a_0 + a_1 + a_2\}$. $|Z_2G| = 8$ and $Z_2G$ is a non-commutative and non-associative ring. We have $\alpha = a_0 + a_1 + a_2$ to be an idempotent for $\alpha^2 = (a_0 + a_1 + a_2)^2 = a_0 + a_1 + a_2$. Clearly $Z_2G$ has no unit element. This is the smallest non-associative ring known to us without having unit element. We define the concepts of subring and ideals.

**DEFINITION 3.2.2:** *Let RG be a groupoid ring. We call a subring V of RG to be a normal subring if*

*(1) aV = Va*
*(2) (Vx)y = V(xy)*
*(3) y(xV) = (yx)V*

*for all $x, y, a \in G$.*

*The groupoid ring is simple if it has no normal subrings. We call RG a normal groupoid ring if*

*(a) x(RG) = (RG)x*
*(b) RG(xy) = ((RG)x)y*
*(c) y(x(RG)) = (yx)(RG)*

*for all $x, y \in RG$.*

**DEFINITION 3.2.3:** *Let RG be a groupoid ring of a groupoid G over the ring R. The centre of $RG = C(RG) = \{x \in RG / xa = ax$ for all $a \in RG\}$.*

**DEFINITION 3.2.4:** *Let RG be a groupoid ring. We say $a, b \in RG$ is a conjugate pair if $a = bx$ (or $xb$) for some $x \in RG$ and $b = ay$ (or $ya$ for some $y \in RG$).*

*We define an element $a \in RG$ to be right conjugate in RG if we can find $x, y \in RG$ such that $ax = b$ and $by = a$ ($xa = b$ and $yb = a$). Similarly we define left conjugate.*

*We define a groupoid ring RG to be a Moufang ring if the Moufang identity $(\alpha\beta)(\gamma\alpha) = (\alpha(\beta\gamma))\alpha$ is satisfied for all $\alpha, \beta, \gamma \in RG$.*

Similarly, we define P-ring, Bol ring and right (left) alternative ring in case of groupoid rings RG.

The study of these rings mainly depends on the structure of the groupoid G but not totally on the groupoid G. We can only prove if RG is a Moufang (Bol, P or alternative) groupoid ring then G is a Moufang groupoid (Bol groupoid, P-groupoid or alternative groupoid). Now we proceed on to define Smarandache groupoid rings and Smarandache non-associative rings.

**DEFINITION 3.2.5:** *Let RG be a groupoid ring; RG is a Smarandache groupoid ring (S-groupoid ring) if and only if G is a S-groupoid.*



**DEFINITION 3.2.6:** *Let RG be a groupoid ring; RG is a SNA-ring if and only if RG has a proper subset which is an associative ring under the operations of RG.*

*Remark:* If RG is a groupoid ring and G is a groupoid with unit then this groupoid ring is always a SNA-ring. It is inevitable to note that it is not very easy to construct groupoids with identity. For when we define groupoids with identity we do so by only adjuring an element with G, but this adjoined element usually is not well adjusted with G so we in all cases do not get all groupoid rings to be SNA-rings as most of the groupoids in the new classes of groupoids do not have identity 1 in them. Thus in all cases we have groupoid rings, which are neither SNA-rings nor S-groupoid rings. All S-groupoid rings are SNA-rings, we prove this in the following theorem.

**THEOREM 3.2.1:** *Let RG be a groupoid ring which is a S-groupoid ring then RG is a SNA-ring.*

*Proof:* By the very definition of the S-groupoid ring we see G is a S-groupoid so G has a proper subset H which is a semigroup. Now consider the semigroup ring RH. Clearly RH $\subset$ RG and RH is an associative ring so RG is a SNA-ring.

**THEOREM 3.2.2:** *A groupoid ring RG which is a SNA-ring need not in general be a S-groupoid ring.*

*Proof:* Suppose G is a groupoid which is not a S-groupoid with unit 1. Then in the groupoid ring R $\subset$ RG is an associative subring but as G is not a S-groupoid, RG is not a S-groupoid ring.

Thus we see the class of S-groupoid rings contains completely the class of groupoid rings which are SNA-rings.

**THEOREM 3.2.3:** *Let RG be the groupoid ring. If g $\in$ G is an idempotent of G and G has unit, then the groupoid ring has idempotents.*

*Proof:* Obvious from the fact $g^2$ = g for g $\in$ G and 1 $\in$ RG as 1 $\in$ G so g $\in$ RG which is an idempotent

*Remark*: If $g^2$ = g in RG then we have g (g – 1) = 0 is a divisor of zero. But the more interesting part of the question is that can every RG have S-idempotents for S-idempotents are defined in non-associative rings only as in the case of associative rings. For the existence of S-idempotents does not demand the groupoid ring to a S-groupoid ring or a SNA-ring, what we need is the existence of idempotents in the groupoid ring.

**THEOREM 3.2.4:** *Let RG be a groupoid ring. If the groupoid G has S-subgroupoids then RG has S-subrings.*

*Proof:* Follows from the very definition of S-subgroupoids and S-subrings. In fact we can say RG has S-groupoid subrings.

*Example 3.2.2:* Let G be a groupoid with identity given by the following table:



| . | 1 | $a_0$ | $a_1$ | $a_2$ |
|---|---|---|---|---|
| 1 | 1 | $a_0$ | $a_1$ | $a_2$ |
| $a_0$ | $a_0$ | 1 | $a_2$ | $a_1$ |
| $a_1$ | $a_1$ | $a_2$ | 1 | $a_0$ |
| $a_2$ | $a_2$ | $a_1$ | $a_0$ | 1 |

$Z_2 = \{0,1\}$ the prime field of characteristics two. $Z_2G$ be the groupoid ring.

$Z_2G = \{0, 1, a_0, a_1, a_2, 1+ a_0, 1+ a_1, 1 + a_2, a_0 + a_1, a_0 + a_2, a_1 + a_2, 1 + a_0 + a_1, 1 + a_0 + a_2, 1 + a_2 + a_1, a_0 + a_1 + a_2, 1 + a_0 + a_1 + a_2\}$.

$|Z_zG| = 16$. We have zero divisors in $Z_2G$, for $(1 + a_0)^2 = 0 = (1 + a_1)^2 = (1 + a_2)^2 = (1 + a_0 + a_1 + a_2)^2 = 0$. $(a_0 + a_1)^2 = 0$, $(a_0 + a_2)^2 = 0$ and $(a_1 + a_2)^2 = 0$ and $(1 + a_0 + a_1)^2 = 1 = (1 + a_0 + a_2)^2 = (1 + a_1 + a_2)^2 = (a_0 + a_1 + a_2)^2$.

In view of this we have the following theorem.

**THEOREM 3.2.5:** *Let $Z_2 = \{0, 1\}$ be the prime field of characteristic two and G be a commutative groupoid with unit of even order in which every $g \in G$ is such that $g^2 = 1$. Then $Z_2 G = K \cup H$ where K consist of nilpotent element of order two and H units of a special form $x^2 = 1$ for all $x \in H$.*

*Proof:* $Z_2H = K \cup H$ where $K = \{\Sigma \alpha_i g_i \mid \alpha_i \in Z_2$ and $g_i \in G$ with $\Sigma \alpha_i = 0\}$ and $H = \{\Sigma \alpha_i g_i / \alpha_i \in Z_2$ and $g_i \in G$ with $\Sigma \alpha_i = 1\}$. It is given for all $g \in G$, $g^2 = 1$ and G is commutative and order of G is even. Hence for any $h \in H$ we have

$$h = \sum_{1=1}^{n} \alpha_i\, g_i$$

where n is odd i.e. $h = g_1 + \ldots + g_n$; $h^2 = 1$ using $g_i^2 = 1$ and $g_i g_j = g_j g_i$. Similarly if $k \in K$ then we have $k^2 = 0$. Hence the claim. If G is of odd order we have such a division only when G is a commutative groupoid with 1 and for every $g \in G$ we have $g^2 = 1$. In $Z_2 G$ the groupoid ring we have $Z_2 G = H \cup K$ as mentioned in the theorem and $H \cap K = \phi$.

It is nice to note that H is not even closed under addition where as K is a subring under the operations of $Z_2G$.

**DEFINITION 3.2.7:** *We say the groupoid ring RG is a Smarandache commutative (S-commutative) if RG has a proper Smarandache subring (S-subring) which is commutative.*

*We say RG is Smarandache strongly commutative (S-strongly commutative) if every proper S-subring of RG is commutative.*

**<u>Remark:</u>** We see even if RG is S-strongly commutative ring still RG need not be commutative. Further if RG is S-strongly commutative then it is S-commutative.



Suppose RG is a commutative groupoid ring which has no S-subrings then RG is not S-strongly commutative or S-commutative.

So for S-commutative to exist RG should posses S-subrings which are commutative even if RG is a commutative groupoid ring still RG need not be S-commutative. We as in case of S-seminormal groupoid define S-seminormal groupoid ring and S-normal groupoid rings.

**DEFINITION 3.2.8:** *Let RG be a groupoid ring. We call a subset W of RG to be a Smarandache seminormal subring (S-seminormal subring). If*

1. *W is a S-subring of RG.*
2. *aW = X for all a $\in$ RG.*
3. *Va = Y for all a $\in$ RG*

*where either X or Y is a S-subring of RG and X and Y are both subrings. We say W is Smarandache normal subring (S-normal subring) if aV = X and Va = Y for all a $\in$ RG where both X and Y are S-subrings of RG.*

**Result:** Every S-normal subring is a S-seminormal subring of RG. Here we define the concept of n-ideal groupoid rings.

**DEFINITION 3.2.9:** *Let RG be a groupoid ring we say RG is a n-ideal groupoid ring if for n-ideals $S_1, S_2, \ldots, S_n$ of RG there exists exactly n-distinct elements $x_1, x_2, \ldots, x_n \in RG \setminus \{S_1 \cup S_2 \cup \ldots \cup S_n\}$ such that $\langle x_1 \cup S_1 \cup S_2 \cup \ldots \cup S_n \rangle = \langle x_2 \cup S_1 \cup S_2 \cup \ldots \cup S_n \rangle = \ldots = \langle x_n \cup S_1 \cup S_2 \cup \ldots \cup S_n \rangle$ where $\langle \rangle$ denotes the ideal generated by $x_i \cup S_i \cup S_2 \cup \ldots \cup S_n$ for $1 \leq i \leq n$.*

**DEFINITION 3.2.10:** *Let RG be a groupoid ring. $S_1, S_2, \ldots, S_n$ be any n S-ideals of RG. We say RG is a Smarandache n-ideal ring (S-n-ideal ring) if there exists n distinct elements $x_1, \ldots x_n$ in $RG \setminus \{S_1 \cup S_2 \cup \ldots \cup S_n\}$ such that $\langle x_1 \cup S_1 \cup S_2 \cup \ldots \cup S_n \rangle = \langle x_2 \cup S_1 \cup S_2 \ldots \cup S_n \rangle = \ldots = \langle x_n \cup S_1 \cup S_2 \cup \ldots \cup S_n \rangle$ where $\langle \rangle$ denotes the S-ideal generated by $x_i \cup S_1 \cup S_2 \cup \ldots \cup S_n$ for $1 \leq I \leq n$.*

It is pertinent to note that n can be any integer value 2 or 3 or 4 or ... Depending on the groupoid G and on the ring R we find whether groupoid ring RG is a S-n-ideal ring or n-ideal ring for some suitable integer n.

***Example 3.2.3:*** Let $Z_2 = \{0, 1\}$ be the prime field of characteristic two and G be a groupoid given by the following table:

| *   | e   | $a_0$ | $a_1$ | $a_2$ | $a_3$ | $a_4$ |
|-----|-----|-------|-------|-------|-------|-------|
| e   | e   | $a_0$ | $a_1$ | $a_2$ | $a_3$ | $a_4$ |
| $a_0$ | $a_0$ | e   | $a_2$ | $a_4$ | $a_1$ | $a_3$ |
| $a_1$ | $a_1$ | $a_2$ | e   | $a_1$ | $a_3$ | $a_0$ |
| $a_2$ | $a_2$ | $a_4$ | $a_1$ | e   | $a_0$ | $a_2$ |
| $a_3$ | $a_3$ | $a_1$ | $a_3$ | $a_0$ | e   | $a_4$ |
| $a_4$ | $a_4$ | $a_3$ | $a_0$ | $a_2$ | $a_4$ | e   |



$Z_2G$ be the groupoid ring, $K_1 = \{0, 1 + a_0, 1 + a_1, 1 + a_2, 1 + a_3, 1 + a_4, a_0 + a_1, a_0 + a_2 + a_0 + a_3, a_0 + a_4, a_1 + a_2, a_1 + a_3, a_1 + a_4, a_2 + a_3, a_2 + a_4, a_3 + a_4, 1 + a_0 + a_1 + a_2, 1 + a_0 + a_1 + a_3, 1 + a_0 + a_1 + a_4, 1 + a_0 + a_2 + a_3, 1 + a_0 + a_2 + a_4, 1 + a_0 + a_3 + a_4, 1 + a_1 + a_2 + a_3, 1 + a_1 + a_2 + a_4, 1 + a_1 + a_3 + a_4, a_0 + a_1 + a_2 + a_3, a_0 + a_1 + a_2 + a_4, a_0 + a_1 + a_3 + a_4, a_0 + a_2 + a_3 + a_4, a_1 + a_2 + a_3 + a_0, 1 + a_1 + a_3 + a_4, 1 + a_2 + a_3 + a_4, 1 + a_0 + a_1 + a_2 + a_3 + a_4\}$, is an ideal of $Z_2G$.

$K_2 = \{0, 1 + a_0 + a_1 + a_2 + a_3 + a_4\}$ is also an ideal of $Z_2G$. Though there are 2 ideals in $Z_2G$ still it is not an n-ideal ring for any n as $K_2 \subset K_1$.

Using the notion of semi-ideals in groupoid we define Smarandache semi-ideals and semi-ideals and Smarandache semi-ideals in groupoid rings.

**DEFINITION 3.2.11:** *Let G be a groupoid. A non-empty subset S of G is said to be a left semi-ideal of G if*

*(1) S is closed under the operations of G*
*(2) $x^2 s \in S$ for all $x \in G$ and $s \in S$.*

*Similarly we can define right semi-ideals of G. We call S a semi-ideal if it is both right and left semi-ideal of G.*

*Example 3.2.4:* G be a groupoid given by the following table:

| * | a | b | c |
|---|---|---|---|
| a | b | c | a |
| b | a | b | c |
| c | c | c | c |

$S = \{b, c\}$ is a semi-ideal of G. We define the concept called the strong left (right) semi-ideal of a groupoid G.

**DEFINITION 3.2.12:** *Let G be a groupoid. A closed subset S of G is said to be strong left (right) semi-ideal of G if $(xy)s \in S$ $(s(xy) \in S)$ for every pair $x, y \in G$ and $s \in S$.*

**THEOREM 3.2.6:** *If G is a groupoid with S a strong semi-ideal of G, then S is a semi-ideal of G.*

*Proof:* Straightforward by the definition.

**DEFINITION 3.2.13:** *Let G be a groupoid. A non-empty subset S of G is said to be a Smarandache left semi-ideal (S-left semi-ideal) of G if*

*(1) S is S-subgroupoid of G*
*(2) $x^2 s \in S$ for all $x \in G$ and $s \in S$*



*Similarly we define Smarandache right semi-ideal (S-right semi-ideal) of G. If S is simultaneously both a Smarandache left and right semi-ideal of G then we call S a Smarandache semi-ideal (S-semi-ideal) of G.*

**DEFINITION 3.2.14:** *Let G be a groupoid, let A be a non-empty subset of G;*

  *(1) A is a S-subgroupoid of G*
  *(2) For every pair x, y $\in$ G and s $\in$ G, (xy) s $\in$ S.*

*Then S is said to be a Smarandache strong left semi-ideal (S-strong left semi-ideal) of G. If s(xy) $\in$ S for every s $\in$ S and x, y $\in$ S we say S is a Smarandache strong right semi-ideal (S-strong right semi-ideal) of G.*

*If S is simultaneously both a S-strong right semi-ideal and S-strong left semi-ideal then we call S a Smarandache strong semi-ideal (S-strong semi-ideal) of G.*

Several interesting results about groupoids and its Smarandache properties can be derived. As the main aim of this book is the study of Smarandache non-associative rings we do not recall or obtain properties about groupoids or S-groupoids in this book. Now we proceed on to define Smarandache generalized semi-ideals in groupoid rings and generalized semi-ideals of groupoid rings.

**DEFINITION 3.2.15:** *Let G be a groupoid, R any field or a commutative ring with 1. The groupoid ring RG of the groupoid G over R. Let V be a non-empty subset of RG. V is said to be a generalized left semi-ideal of RG if V is closed under addition and $x^2 s \in V$ for every $s \in V$ and $x \in RG$. Similarly we can define generalized right semi-ideals of RG. If V is simultaneously a generalized right semi-ideal and left semi-ideal of RG then we call RG a generalized semi-ideal.*

**DEFINITION 3.2.16:** *Let RG be a groupoid ring. We say a non-empty subset V of RG is a Smarandache generalized left semi-ideal (S-generalized left semi-ideal) of RG if*

  *(1) V is a S-semigroup under '+'.*
  *(2) For all $s \in V$ and $x \in RG$ we have $x^2 s \in RG$.*

*We define similarly Smarandache generalized right semi-ideal (S-generalized right semi-ideal) of RG. If V is simultaneously both a S-generalized right and left semi-ideal of RG then we call V a Smarandache generalized ideal (S-generalized ideal) of RG.*

**THEOREM 3.2.7:** *Every generalized semi-ideal of a groupoid ring in general need not be a S-generalized semi-ideal of the groupoid ring of RG.*

*Proof:* We prove this by an example. Consider the groupoid ring $Z_2 G$ where G is of finite order say G = $\{g_1, \ldots, g_n\}$ such that $g_i (g_1 + \ldots + g_n) = g_1 + \ldots + g_n$.

$$V = \left\{0, \sum_{l=1}^{n} g_i \right\}$$



is a generalized semi-ideal of $Z_2G$, but clearly V is not a S-generalized semi-ideal of $Z_2G$.

**THEOREM 3.2.8:** *Let $Z_2G$ be the groupoid ring where G is a commutative groupoid with 1 such that $g^2 = 1$ for every of $g \in G$. Then $Z_2G$ has*

> *(1) S-generalized semi-ideal.*
> *(2) Generalized semi-ideals which is not a S-generalized semi-ideal.*

*Proof:* Consider the groupoid ring $Z_2G$. Let $V = \{\Sigma\alpha_I g_i \mid \Sigma\alpha_I = 0\}$. Clearly it is verified that V is a S-generalized semi-ideal of $Z_2G$. Now $K = \{0, 1 + g_1 + \ldots + g_n\}$ is a generalized semi-ideal which is not a S-generalized semi-ideal of $Z_2G$.

**PROBLEMS:**

1. Find normal subrings and S-normal subrings if any in the groupoid ring ZG where G is the groupoid from $Z^{***}(n)$.
2. Find S-seminormal subrings and seminormal subrings in the groupoid ring $Z_2G$ where $G \in Z^{**}(n)$.
3. Find all units and idempotents in $Z_5G$ where $G \in Z^{**}(7)$.
4. Find a conjugate pair in $Z_5G$, $G \in Z^{*}(8)$.
5. Give an example of a groupoid ring RG which is S-Moufang.
6. Find an example of a S-Bol groupoid ring.
7. Give an example of a S-groupoid ring.
8. Give an example of a groupoid ring which is not a S-groupoid ring.
9. Give an example of a groupoid ring which has S-generalized semi-ideals.
10. Give an example of a groupoid ring, which has no S-generalized semi-ideals.
11. Give an example of a groupoid ring, which has S-normal subring.
12. Find for the groupoid ring which has
    (1) S-5-ideal ring.
    (2) 5-ideal ring.

### 3.3 Smarandache special elements in Groupoid rings

In this section we specially study the properties of not only elements like units idempotents etc. but also special propertied rings like E-rings, p-rings, pre J-rings, zero square rings, quasi commutative rings, pseudo commutative rings etc and specially their Smarandache analogue. Such study in case of groupoid rings in totally absent.

The concept of S-units, S-zero divisors and S-idempotents for any non-associative ring is also the same as that of the associative rings. So we request the reader to refer [68]. Study of such concepts for loop rings have been carried out we in this section apart from study the existence of units, S-units, zero divisors, S-zero divisor, idempotents and S-idempotents in groupoid rings and S-groupoid rings we study and introduce concepts like normal elements, quasi regular elements, semi idempotents, pseudo commutative elements and their Smarandache analogue.



*Example 3.3.1*: Let $Z_2 = \{0, 1\}$ be the prime field of characteristic two. G be a groupoid given by the following table:

| . | e | $a_0$ | $a_1$ | $a_2$ |
|---|---|---|---|---|
| e | e | $a_0$ | $a_1$ | $a_2$ |
| $a_0$ | $a_0$ | e | $a_2$ | $a_1$ |
| $a_1$ | $a_1$ | $a_2$ | e | $a_0$ |
| $a_2$ | $a_2$ | $a_1$ | $a_0$ | e |

Take $a = (1+a_0 + a_1 + a_2)$ and $b = (1+a_0)$. Clearly ab = 0. Choose x, y, x ≠ y such that $x = a_2 + a_0$ and $y = a_2 + a_1$.

a . x = 0 and by = 0 but

$$\begin{aligned} x \cdot y &= (a_0 + a_2)(a_1 + a_2) \\ &= a_0 a_1 + a_2 a_1 + a_0 a_2 + a_2^2 \\ &= a_2 + a_0 + a_1 + 1 \\ &\neq 0. \end{aligned}$$

Thus we see $Z_2G$ has S-zero divisors. In view of this example we have the following theorem.

**THEOREM 3.3.1:** *Let $Z_2 = \{0, 1\}$ be the field and G be a finite commutative groupoid with 1 of even order such that $g^2 = 1$, for every $g \in G$. The groupoid ring $Z_2G$ has S-zero divisors.*

*Proof*: Let $G = \{1, g_1, \ldots, g_n\}$, where n is odd so |G| is even. Let $a = 1 + g_1 + \ldots + g_n$ and $b = 1 + g_1$ be in $Z_2G$. Clearly ab = 0. Choose $x = g_1 + g_2$ and $y = g_2 + g_3$ such that ax = 0 and by = 0 with xy ≠ 0. Thus the groupoid ring has S-zero divisors.

*Example 3.3.2*: Let $Z_2 = \{0,1\}$ and G be a groupoid given by the following table:

| * | e | $a_1$ | $a_2$ | $a_3$ | $a_4$ | $a_5$ | $a_6$ | $a_7$ |
|---|---|---|---|---|---|---|---|---|
| e | e | $a_1$ | $a_2$ | $a_3$ | $a_4$ | $a_5$ | $a_6$ | $a_7$ |
| $a_1$ | $a_1$ | e | $a_5$ | $a_2$ | $a_6$ | $a_3$ | $a_7$ | $a_4$ |
| $a_2$ | $a_2$ | $a_5$ | e | $a_6$ | $a_3$ | $a_7$ | $a_4$ | $a_1$ |
| $a_3$ | $a_3$ | $a_2$ | $a_6$ | e | $a_7$ | $a_4$ | $a_1$ | $a_5$ |
| $a_4$ | $a_5$ | $a_6$ | $a_3$ | $a_7$ | e | $a_1$ | $a_5$ | $a_2$ |
| $a_5$ | $a_5$ | $a_3$ | $a_7$ | $a_4$ | $a_1$ | e | $a_2$ | $a_6$ |
| $a_6$ | $a_6$ | $a_7$ | $a_4$ | $a_1$ | $a_5$ | $a_2$ | e | $a_3$ |
| $a_7$ | $a_7$ | $a_4$ | $a_1$ | $a_5$ | $a_2$ | $a_6$ | $a_3$ | e |

$Z_2G$ is the groupoid ring. $(1 + a_i + a_j)^2 = 1$. For $a_c, a_i \in G$ ($a_i \neq a_j$) and $(a_i + a_j + a_k)^2 = 1$ for $a_i, a_j, a_k \in G$, $i \neq j$, $j \neq k$ and $i \neq k$. Similarly $(1 + a_i + a_j + a_k + a_m) = 1$. $a_i, a_j, a_k$,



$a_m$ are distinct elements of G, also $(a_i + a_j + a_k + a_l + a_m)^2 = 1$ where $a_i, a_j, a_k, a_l, a_m$ are 5 distinct elements of G. Now,

$$\left(\sum_{l=1}^{7} a_i\right)^2 = 1$$

Thus we have units. The study of how many of these units happens to be S-units is an interesting one. Further we see every S-unit is a unit but all units in general are not S-units.

Now we proceed on to define pseudo commutative elements.

**DEFINITION 3.3.1:** *Let RG be non-commutative groupoid ring that is, the groupoid G is non-commutative. A pair of distinct elements $x, y \in RG$ different from the identity of RG which are such that $xy = yx$ is said to be a pseudo commutative pair of R if $(xa)y$ or $((ya)x)$ for all $a \in RG$. If in the groupoid ring RG every commutative pair happens to be a pseudo commutative pair of RG then we say RG is a pseudo commutative ring. Clearly if G is a commutative groupoid the groupoid ring RG is a pseudo commutative ring.*

**Note:** For us to define a pseudo commutative pair in the groupoid ring what we need is that if $xy = yx$. We can have at least one of the equalities to be true for all $a \in RG$ is $(xa)y = y(ax)$ or $= (ya)x$ or $x(ay) = (ya)x$ or $y(ax)$.

Now we define Smarandache pseudo commutative groupoid ring.

**DEFINITION 3.3.2:** *Let RG be a groupoid ring. Let A be a S-subgroupoid ring of RG. A pair of elements $x, y \in A$ which are such that $xy = yx$ is said to be a Smarandache pseudo commutative pair (S-pseudo commutative pair) of RG if $x(ay) = y(ax)$ or $(ya)x$ or $(xa)y = y(ya)x)$ or $y(ax)$ for all $a \in A$. If in the S-subgroupoid ring A every commutative pair happens to be a S-pseudo commutative pair of A, then A is said to be a Smarandache pseudo commutative groupoid ring (S-pseudo commutative groupoid ring).*

**THEOREM 3.3.2:** *Let RG be a groupoid ring. If RG is a S-pseudo groupoid ring then RG is a SNA-ring.*

*Proof*: Follows from the fact that RG is a S-pseudo commutative ring then by the very definition RG is a SNA-ring.

**DEFINITION 3.3.3:** *The groupoid ring RG is strongly regular if for every $x, y \in RG$ we have $(xy)^n = xy$ for $n = n(x,y) > 1$, then R is called a strongly regular ring. We say the groupoid ring RG is a Smarandache strongly regular ring (S-strongly regular ring) if RG contains a SNA-subring A such that for every $x, y \in A$ we have $(xy)^n = xy$ for some integer $n = n(x, y) > 1$.*

It is left as an exercise for the reader to prove.

**THEOREM 3.3.3:** *If RG is a S-strong regular ring then RG is a SNA-ring.*



Such type of study has not been carried out by any researcher in case of any of the non-associative rings be it loop rings, groupoid rings, Lie rings or Jordan rings. Now we proceed on to define quasi commutatively in case of groupoid rings RG and their Smarandache analogue.

**DEFINITION 3.3.4:** *Let RG be a groupoid ring we say RG is quasi commutative if $ab = b^t a$ for every pair of elements $a, b \in R$ and $t > 1$.*

**THEOREM 3.3.4:** *Let RG be a groupoid ring, if RG is quasi commutative then every pair of elements in the groupoid G satisfies $ab = b^t a$, $t > 1$.*

*Proof*: Clearly given in the groupoid ring RG we have $ab = b^t a$, $t > 1$ for every $a, b \in$ RG, further we know $I.G \subseteq RG$ so in G we have $ab = b^t a$ for every pair $a, b \in G$; $t > 1$.

*Remark*: If G is a groupoid with 1 and if we have RG to be quasi commutative then, R is quasi commutative. If G has no unit element 1 then we may not be in a position to say even if RG in quasi commutative so is R.

**THEOREM 3.3.5:** *Let RG be a quasi commutative groupoid ring. Then for every pair of elements $a, b \in RG$ there exists $s \in RG$ such that $a^2 b = bs^2$ provided $(b.b^{\gamma-1}) a = b(b^{\gamma-1} a)$.*

*Proof*: RG is a quasi commutative groupoid ring, so $ab = b^\gamma a$ for every pair of elements $a, b \in RG$, $\gamma > 1$. Now $ab = b^\gamma a$ so $a^2 b = ab^\gamma a = (ab) b^{\gamma-1} a = (b^\gamma a)(b^{\gamma-1} a) = b (b^{\gamma-1} a)(b^{\gamma-1} a) = bs^2$ (using $b(b^{\gamma-1} a) = (b. b^{\gamma-1}) a$ where $s^2 = (b^{\gamma-1} a)^2$).

**DEFINITION 3.3.5:** *Let RG be a groupoid ring. We say RG is a Smarandache quasi-commutative ring (S-quasi-commutative ring) if for any S-subring A of RG we have $ab = b^\gamma a$ for every $a, b \in A$, $\gamma > 1$.*

**THEOREM 3.3.6:** *If RG is a S-quasi commutative ring then RG is a SNA-ring.*

*Proof*: Obvious by the very definitions.

The reader is expected to prove the following theorem.

**THEOREM 3.3.7:** *Every SNA-ring in general need not be S-quasi commutative.*

**DEFINITION 3.3.6:** *An element x of the groupoid ring RG is called semi nilpotent if $x^n - x$ is a nilpotent element of RG. If $x^n - x = 0$ we say x is a trivial semi nilpotent.*

**THEOREM 3.3.8:** *If x is a nilpotent element of the groupoid ring RG then x is a semi nilpotent element of RG.*

*Proof*: Straightforward from the fact $x^n = 0$ so $x^n - x = -x$. so $(-x)^n = 0$, hence our claim.



**DEFINITION 3.3.7:** *Let RG be a groupoid ring. An element $x \in RG$ is a Smarandache semi nilpotent (S-semi nilpotent) if $x^n - x$ is S-nilpotent. An element $0 \neq x \in RG$ is said to be Smarandache nilpotent if $x^n = 0$ and there exists a $y \in R \setminus \{0, x\}$ such that $x^r y = 0$ or $yx^s = 0$, $r, s > 0$ and $y^m \neq 0$ for any integer $m > 1$.*

**DEFINITION 3.3.8:** *A groupoid ring RG is reduced if R has no non-zero nilpotent elements.*
*We say RG is Smarandache-reduced (S-reduced) if RG has no S-nilpotents.*

The following theorem is straightforward and hence left for the reader to prove.

**THEOREM 3.3.9:** *Let RG be a reduced groupoid ring then RG is a SNA-reduced ring. If R is a SNA-reduced ring then RG need not in general be a reduced ring.*

**DEFINITION 3.3.9:** *A groupoid ring RG is a zero square ring if $x^2 = 0$ for all $x \in RG$.*

**DEFINITION 3.3.10:** *A groupoid ring RG is a Smarandache zero square ring (S-zero square ring) if RG has a S-subring which is a zero square ring.*

***Example* 3.3.3:** Let $Z_2 = \{0, 1\}$ be the prime field of characteristic two and G be a commutative groupoid with unit given by the following table:

| . | 1 | $a_0$ | $a_1$ | $a_2$ | $a_3$ | $a_4$ |
|---|---|---|---|---|---|---|
| 1 | 1 | $a_0$ | $a_1$ | $a_2$ | $a_3$ | $a_4$ |
| $a_0$ | $a_0$ | 1 | $a_3$ | $a_1$ | $a_4$ | $a_2$ |
| $a_1$ | $a_1$ | $a_3$ | 1 | $a_4$ | $a_2$ | $a_0$ |
| $a_2$ | $a_2$ | $a_1$ | $a_4$ | 1 | $a_0$ | $a_3$ |
| $a_3$ | $a_3$ | $a_4$ | $a_2$ | $a_0$ | 1 | $a_1$ |
| $a_4$ | $a_4$ | $a_2$ | $a_0$ | $a_3$ | $a_1$ | 1 |

Clearly
$$A = \left\{ \sum \alpha_i a_j \mid a_i \in G, a_i \text{ can also be 1 with } \sum \alpha_i = 0 \right\}$$

is a S-subgroupoid ring which is a S-zero square ring as $x^2 = 0$ for all $x \in A$.

Thus we have the following nice theorem.

**THEOREM 3.3.10:** *Let $Z_2 = \{0,1\}$ be the prime field of characteristic two and G be a commutative groupoid of even order with 1 such that $g_i^2 = 1$ for every $g_i \in G$. Then the groupoid ring $Z_2 G$ is a S-zero square ring.*

*Proof:* Set $U = \left\{ \sum \alpha_i g_i \mid g_i \in G \text{ with } \sum \alpha_i = 0 \right\}$; U is such that every element x is of the form $x^2 = 0$. U is a S-subgroupoid ring. Hence $Z_2 G$ is a S-zero square ring.

**THEOREM 3.3.11:** *Let RG be a groupoid ring such RG is a S-zero square ring then RG need not be a zero square ring.*



*Proof*: We prove this by the following example. Consider the groupoid ring RG given in the previous theorem; clearly in RG we have elements x ∈ RG say x = g ∈ RG, g ∈ G such that $g^2 = 1$. Hence RG is a zero square ring.

**DEFINITION 3.3.11:** *Let RG be a groupoid ring. RG is called an inner zero square ring if every proper subring of RG is a zero square ring.*

*We call the groupoid ring RG to be Smarandache inner zero square ring (S-inner zero square ring) if every proper S-subring of RG is a zero square ring.*

Several interesting results can be got from these definitions.

**DEFINITION 3.3.12:** *We call the groupoid ring RG to be a Smarandache weak inner zero square ring (S-weak inner zero square ring) if RG has at least a S-subring A ⊆ R G such that a subring B of A is a zero square ring.*

**DEFINITION 3.3.13:** *A groupoid ring RG is said to be p-ring if $x^p = x$ and px = 0 for every x ∈ R. We say the groupoid ring RG is said to be Smarandache p-ring (S-p-ring) if RG is a S-ring and RG has a subring P such that $x^P = x$ and px = 0 for every x ∈ P.*

**THEOREM 3.3.12:** *Let G be a groupoid with nilpotent elements and R any ring. The groupoid ring RG is not a p-ring.*

*Proof*: Straightforward from the very definition.

**DEFINITION 3.3.14:** *Let RG be a groupoid ring. R G is said to be a S-p-ring if RG is a SNA-ring and has a subring A such that A is a p-ring.*

Thus the following theorem is direct hence left for the reader to solve.

**THEOREM 3.3.13:** *If RG is a groupoid ring which is a p-ring then RG is a S-p-ring.*

**THEOREM 3.3.14:** *Every Smarandache p-ring RG need not in general be a p-ring.*

*Proof*: We prove by an example. Take RG to be a groupoid ring where R is a p-ring and G is a groupoid with nontrivial nilpotents. Then the groupoid ring RG is a S-p-ring and not a p-ring. Hence the claim.

Now we proceed on to define a particular case of p-rings in case of associative rings they are termed as E-ring, So far they have not extended this notion to non-associative rings. Further this notion cannot be defined on Lie rings which also forms a class of non-associative rings. We can define the concept of E-rings for loop rings and groupoid rings that too only when we take the rings R to be of characteristic two or the prime field of characteristic 2.

**DEFINITION 3.3.15:** *Let RG be a groupoid ring. We call RG an E-ring if $x^{2n} = x$ and 2x = 0 for every x in RG and n a positive integer. The minimal such n is called the degree of the E-ring.*



We have the following nice theorem from the definition which we leave as an exercise to the reader.

**THEOREM 3.3.15:** *Let RG be a E-ring then we have characteristic of R is 2.*

**DEFINITION 3.3.16:** *Let RG be a groupoid ring. RG is said to be a Smarandache E-ring (S-E-ring) if RG has a proper subset A which is SNA-subring of RG and A is a E-ring.*

Thus we have the following.

**THEOREM 3.3.16:** *Let RG be a groupoid ring which is a E-ring. If RG has a proper SNA-subring. A then RG is a S-E-ring.*

*Proof*: Straightforward.

**THEOREM 3.3.17:** *Let RG be a groupoid ring. Every S-E-ring need not in general be an E-ring.*

*Proof*: By an example. Consider the groupoid ring $Z_2G$ where G is given by the following table:

| . | 1 | $a_0$ | $a_1$ | $a_2$ |
|---|---|---|---|---|
| 1 | e | $a_0$ | $a_1$ | $a_2$ |
| $a_0$ | $a_0$ | e | $a_2$ | $a_1$ |
| $a_1$ | $a_1$ | e | $a_2$ | $a_1$ |
| $a_2$ | $a_2$ | $a_1$ | $a_0$ | e |

$P = \{1, a_0, a_1, a_2, 1 + a_0 + a_1, 1 + a_0 + a_2, 1 + a_1 + a_2, a_0 + a_1 + a_2\}$ is a SNA-ring which is a E-ring.

So RG is a Smarandache E-ring which is not a E-ring. Now we proceed on to define pre-J-ring in groupoid rings.

**DEFINITION 3.3.17:** *Let RG be a groupoid ring. We say RG is a pre-J-ring if $a^n b = ab^n$ for any pair $a, b \in R G$ and n a positive integer.*

**DEFINITION 3.3.18:** *Let RG be a groupoid ring we say RG is a Smarandache pre-J-ring (S-pre-J-ring) if for every pair of elements $a, b \in P$, (where P is a proper subset of RG which is a S-subring of RG) we have $a^n b = ab^n$ for some positive integer n.*

**THEOREM 3.3.18:** *Suppose R is a pre-J-ring and G any groupoid, the groupoid ring RG is a S-pre-J-ring but is general not a pre-J-ring.*

*Proof*: Given R is a pre-J-ring, so for any groupoid, G we have RG to be a S-pre-J-ring as $R \subset RG$. Now to prove RG in general is not a pre-J-ring even if RS is a S pre-J-ring. Suppose G is a groupoid having elements like say $a^n b$ to exists and different



from $ab^n = 0$ then we have the result to be true consider $(1 + a_0) = b$ and $a = a_0$ we see $a_0^n (1 + a_0) = 1+a_0$ but $a.b^n = a_0 (1+a_0) = 1 + a_0$ but $a.b^n = a_0$, $(1+ a_0)^n = 0$. Hence the claim.

Several such properties studied as in case of associative rings can also be found and introduced in case of the special class of non-associative rings viz loop rings and groupoid rings.

**PROBLEMS:**

1. Does the groupoid ring $Z_5G$ where $G = Z_7$ (3,2) have pseudo commutative pair 2?
2. Does the groupoid ring $Z_5G$ given in example 1 have a S-pseudo commutative pair 2.
3. Find all the strongly regular elements of the groupoid ring $Z_4G$ where $G = Z_{10}(3, 5)$
4. Does the groupoid ring given in problem 3 have S-strongly regular element?
5. Find a quasi commutative pair in $Z_9(3, 7)$.
6. Give an example of a groupoid ring which has no semi nilpotent elements.
7. Find groupoid rings, which has no semi nilpotent elements.
8. Give an example of a S-reduced groupoid ring.
9. Can $Z_2G$ where $G = Z_{11}$ (3,2) be a S-zero square ring? Justify your answer.
10. Give an example of each of a groupoid ring which is a
    a. E-ring.
    b. S-E-ring.
    c. p-ring.
    d. S-p-ring.
    e. pre-J-ring.
    f. S-pre J-ring.
    g. inner zero square ring.
    h. S-weak inner zero square ring.
    i. S-quasi commutative ring.

### 3.4. Smarandache substructures in Groupoid rings

This section is completely devoted to the introduction of Smarandache substructures in groupoid rings. Study of groupoid rings is totally absent in literature, so the study of Smarandache substructures in groupoid rings seems not only interesting but also very new. These new class of non-associative rings is distinct from Lie algebras or Jordan algebras or loop rings. First we introduce the concept of mod-p-envelope of a groupoid ring. The notion of mod p envelope for group rings was first started [38] in the year 1978, and later introduced in [55] in the year 1986 and [62]. For more about group rings and above all finite loop rings refer [62].

**DEFINITION 3.4.1:** *Let G be a finite groupoid and $Z_p$ be the prime field of characteristic p. $Z_pG$ be the groupoid ring of the groupoid G over $Z_p$. The mod p-envelope of G denoted by $G^* = 1+U$ where*

$$U = \{\alpha = \Sigma \alpha_i g_i \in Z_p G | \Sigma \alpha_i = 0\}.$$



**Example 3.4.1:** Let $Z_2 = \{0, 1\}$ be the prime field of characteristic two. G be the groupoid given by the following table:

|     | 0   | 1   | $a_1$ | $a_2$ | $a_3$ | $a_4$ |
| --- | --- | --- | ----- | ----- | ----- | ----- |
| 1   | 1   |     | $a_1$ | $a_2$ | $a_3$ | $a_4$ |
| $a_1$ | $a_1$ | 1   |       | $a_3$ | $a_2$ | $a_1$ |
| $a_2$ | $a_2$ | $a_1$ | 1     |       | $a_3$ | $a_2$ |
| $a_3$ | $a_3$ | $a_2$ | $a_1$ | 1     |       | $a_3$ |
| $a_4$ | $a_4$ | $a_3$ | $a_2$ | $a_1$ | 1     |       |

$Z_2$ G is the groupoid ring of G over $Z_2$. Now $U = \{0, 1 + a_1, 1 + a_2, 1 + a_3, 1 + a_4, a_1 + a_2, a_1 + a_3, a_1 + a_4, \ldots, a_3 + a_4, 1 + a_1 + a_2 + a_3, \ldots, a_1 + a_2 + a_3 + a_4\}$. $G^* = 1 + U = \{1, a_1, a_2, a_3, a_4, 1 + a_1 + a_2, 1 + a_1 + a_3, \ldots, 1 + a_1 + a_2 + a_3 + a_4\}$. It is easily verified $G^*$ is a groupoid under multiplication and $G \subset G^*$.

We have nice theorems about the mod p envelope of the groupoid $G^*$.

**THEOREM 3.4.1:** *Let G be a commutative groupoid of order 2n in which square of every element is one and $Z_2 = \{0,1\}$ be the prime field of characteristic two. Then $G^*$ is a commutative groupoid such that the square of every element is 1 and the order of $G^*$ is $2^{2n-1}$.*

*Proof:* Let G be the groupoid with elements, $G = \{1, g_1, g_2, \ldots g_{2n-1}, g_i^2 = 1$ and $g_i g_j = g_j g_i, i = 1, 2, \ldots, 2n - 1\}$. $Z_2G$ be the groupoid ring the mod p-envelop of G, $G^* = 1 + U$ where $U = \{0, 1 + g_1, 1 + g_2, 1 + g_3, \ldots, 1 + g_{2n-1}, 1 + g_1 + g_2 + \ldots + g_{2n-1}\}$. $G^* = \{1, g_1, g_2, \ldots, g_{2n-1}, \ldots, g_1 + g_2 + \ldots + g_{2n-1}\}$. Now $\alpha = \sum \alpha_i g_i \in G^*$ is such that $\alpha^2 = 1$ as $g_i^2 = 1$ and $g_i, g_j = g_j g_i$ and characteristic of $Z_2$ is two. Thus $G^*$ is a groupoid in which square of every element is 1 and G is strictly contained in $G^*$ and $|G^*| = 2^{2n-1}$. Clearly if G is non-commutative we may not be in a position to arrive at this conclusion even if each $g_i^2 = 1$ in G.

Now we study the case when $Z_2$ is replaced by any other prime field $Z_p$.

**Example 3.4.2:** Let $Z_3$ be the prime field of characteristic three. G be a groupoid given by the following table:

| *   | e   | $a_0$ | $a_1$ | $a_2$ |
| --- | --- | ----- | ----- | ----- |
| e   | e   | $a_0$ | $a_1$ | $a_2$ |
| $a_0$ | $a_0$ | e     | $a_2$ | $a_1$ |
| $a_1$ | $a_1$ | $a_2$ | e     | $a_0$ |
| $a_2$ | $a_2$ | $a_1$ | $a_0$ | e     |

We see G is commutative with $a_i^2 = e$ for each $a_i \in G$.



Now let $Z_3G$ be the groupoid ring with $1.e = e.1 = e=1$ acting as the identity of $Z_3G$.
$G^* = 1 + U = 1 + \{0, 1 + 2a_0, 1 + 2a_1, 1 + 2a_2, 2 + a_0, 2 + a_1, 2 + a_2, a_1 + 2a_0, a_0 + 2a_1$
$a_0 + 2a_2, a_1 + 2a_2, a_2 + 2a_0, a_2 + 2a_1, 1 + a_0 + a_1, 1 + a_0 + a_2, 1+a_1+a_2, a_0 + a_1+ a_2, 2+ 2a_0 + 2 a_1, 2 + 2a_0 + 2a_2, 2 + 2a_2 + 2a_1, 2a_0 + 2a_1 + 2a_2, 1+ a_0 + 2a_1 + 2a_2, 1 + a_1 + 2a_0 + 2a_2, 1 + a_2 + 2a_0 +2a_1, 2 + a_0 + 2a_1 + a_2, 2+ a_0 + a_1 + 2a_2, 2+2 a_0 + a_1 + a_2\} =$
$\{1, 2 + 2a_0, 2 + 2a_1, 2 + 2a_2, a_0, a_1, a_2, 2 + a_0 + a_1, 2 + a_0 + a_2, 2 + a_1 + a_2, 1 + a_0 + a_1 + a_2, 2a_0 + 2a_1, 2a_0 + 2a_2, 2a_2 + 2a_1, 1 + 2a_0 + 2a_1 + 2a_2, 2 + a_0 + 2a_1 + 2a_2, 2 + a_1 + 2a_0 + 2a_2, 2 + a_2 + 2a_0 + 2a_1, 1 + a_1 + 2a_0, 1 + a_0 + 2a_1, 1 + a_0 + 2a_2, 1 + a_1 + 2a_2, a_0 + 2a_1 + a_2, a_0 + a_1 + 2a_2, 1 + a_2 + 2a_0, 1 + a_2 + 2a_1, 2a_0 + a_1 + a_2\}$. Clearly $G^*$ is closed under the operation product and $G \subset G^*$. Further $|G^*| = 27 = 3^3 = 3^{4-1}$

We see only some of the elements are idempotents some of them units and so on. In view of this we have the following theorem:

**THEOREM 3.4.2:** *Let G be a commutative finite groupoid of even order with unit in which every $g_i \in G$ in such that $g_i^2 = 1$, $Z_p$ ($p > 2$) be the prime field of characteristic p. $Z_pG$ be the groupoid ring of G over $Z_p$. Then $G^* = 1+ U$ the mod p envelope of G is a groupoid of order $p^{2n-1}$.*

*Proof*: Given $|G| = 2n$ and G is commutative with $g^2 = 1$ with unit. $Z_pG$ the groupoid ring of G over $Z_p$. $G^* = 1+ U = 1+ \{0, \sum \alpha_i g_i \in Z_pG | \sum \alpha_i = 0\}$. It is easily verified $G^*$ is a groupoid of order $p^{2n-1}$. Further $G \subset G^*$.

If p / n then the extra properties enjoyed by $G^*$ is an interesting one.

Now we proceed on to define the Smarandache-mod-p-envelop of the groupoid G.

**DEFINITION 3.4.2:** *Let G be a finite groupoid. $Z_p$ the prime field of characteristic p. $Z_pG$ the groupoid ring of the groupoid G over $Z_p$. The Smarandache mod p-envelope (S-mod-p envelope) of G denoted by $SG^* = 1 + U$ where $U = \{0, \Sigma\alpha_i g_i \in Z_pB \mid g_i \in B \subset G;$ B a S-subgroupoid of G with $\Sigma\alpha_i = 0\}$.*

*Thus unlike the mod p-envelope of G we have for a given groupoid G several mod p-envelopes. The importance is to find whether $SG^*$ is a groupoid or a S-groupoid. If G has no S-subgroupoid but G is a S-groupoid then we assume the mod-p-envelope of G to be $G^* = SG^*$.*

Thus the study is interesting only when the groupoids G has S-sub groupoids.

*Example 3.4.3*: Let $Z_2 = \{0, 1\}$ be the prime field of characteristic two and G be the groupoid given by the following table:

| * | e | $a_0$ | $a_1$ | $a_2$ |
|---|---|---|---|---|
| e | e | $a_0$ | $a_1$ | $a_2$ |
| $a_0$ | $a_0$ | e | $a_2$ | $a_1$ |
| $a_1$ | $a_1$ | $a_2$ | e | $a_0$ |
| $a_2$ | $a_2$ | $a_1$ | $a_0$ | e |



Let $Z_2$ G be the groupoid ring. G is a S-groupoid, but G does not contain S-sub groupoid. So $G^* = SG^*$.

***Example 3.4.4:*** Let $Z_2 = \{0, 1\}$ be the prime field of characteristic two and G be a groupoid given by the following table.

| * | e | $a_1$ | $a_2$ | $a_3$ | $a_4$ |
|---|---|---|---|---|---|
| e | e | $a_1$ | $a_2$ | $a_3$ | $a_4$ |
| $a_1$ | $a_1$ | $a_1$ | $a_3$ | $a_1$ | $a_3$ |
| $a_2$ | $a_2$ | $a_3$ | $a_1$ | $a_3$ | $a_1$ |
| $a_3$ | $a_3$ | $a_1$ | $a_3$ | $a_1$ | $a_3$ |
| $a_4$ | $a_4$ | $a_3$ | $a_1$ | $a_3$ | $a_1$ |

$Z_2G$ be the groupoid ring. $Z_2G$ has S-subgroupoids given by $H = \{e, a_1, a_3\}$. Identifying e with 1 we get. $SG^* = 1 + \{0, 1 + a_1, 1 + a_3, a_1 + a_3\} = \{1, a_1, a_3, 1 + a_1 + a_3\}$. Clearly $H \subseteq SG^*$, $SG^*$ has unit for $(1 + a_1 + a_3) = 1$. $H_1 = \{a_1, a_3\}$ is a S-sub groupoid. Then $SG^* = 1 + \{a_1 + a_3, 0\} = \{1, 1 + a_1 + a_3\}$. Clearly $H_1 \not\subseteq SG^*$. Further these S mod-p-envelopes are related.

Thus we have the following.

**THEOREM 3.4.3:** *Let G be a groupoid with 1, $Z_2$ be the prime field of characteristic two. The S-mod p envelope of G associated with proper S-sub groupoids H of G need not in general have $H \subset SG^*$.*

*Proof*: By an example. In the above example we have for $H_1 = \{a_1, a_3\}$ which is a S-subgroupoid of G. The $SG^*$ associated with the S-subgroupoid H of G is such that $SG^* = \{1, 1 + a_1 + a_3\}$. Clearly $H \not\subseteq G^*$. Hence the claim.

Now we proceed on to define the concept of Smarandache pseudo right ideals in groupoid rings.

**DEFINITION 3.4.3:** *Let RG be a groupoid ring. B a proper subset of RG which is an associative subring of RG. A non-empty subset X of RG is said to be a Smarandache pseudo right ideal (S-pseudo right ideal) of RG related to B if*

*(1) (X, +) is an additive abelian group.*
*(2) For $b \in B$ and $s \in S$ we have $sb \in S$.*

*On similar lines we define Smarandache pseudo left ideal (S-pseudo left ideal). A nonempty subset X of RG is said to be Smarandache pseudo ideal (S-pseudo ideal) of RG if X is both a S-pseudo right ideal and S-pseudo left ideal.*

***Example 3.4.5:*** Let ZG be the groupoid ring of a groupoid G given by the following table.



| . | e | $g_0$ | $g_1$ | $g_2$ | $g_3$ | $g_4$ | $g_5$ |
|---|---|---|---|---|---|---|---|
| e | e | $g_0$ | $g_1$ | $g_2$ | $g_3$ | $g_4$ | $g_5$ |
| $g_0$ | $g_0$ | $g_0$ | $g_4$ | $g_2$ | $g_0$ | $g_4$ | $g_2$ |
| $g_1$ | $g_1$ | $g_2$ | $g_0$ | $g_4$ | $g_2$ | $g_0$ | $g_4$ |
| $g_2$ | $g_2$ | $g_4$ | $g_2$ | $g_0$ | $g_4$ | $g_2$ | $g_0$ |
| $g_3$ | $g_3$ | $g_0$ | $g_4$ | $g_2$ | $g_0$ | $g_4$ | $g_2$ |
| $g_4$ | $g_4$ | $g_2$ | $g_0$ | $g_4$ | $g_2$ | $g_0$ | $g_4$ |
| $g_5$ | $g_5$ | $g_4$ | $g_2$ | $g_0$ | $g_4$ | $g_2$ | $g_0$ |

G is a S-groupoid, G has S-subgroupoid. $Z \subset ZG$ is a S-groupoid ring for Z is the associative subring. Let $\{e, g_0, g_2, g_4\} = H$. $ZH \subset ZG$ and $ZH$ is a S-pseudo ideal of ZG.

**DEFINITION 3.4.4:** *Let RG be a groupoid ring which is a S-ring i.e. RG has a proper subset A which is an associative ring. Let X be the S-pseudo ideal related to A. X is said to be Smarandache minimal pseudo ideal (S-minimal pseudo ideal) of RG if $X_1$ is another S-pseudo ideal related to A and $(0) \subset X_1 \subset X$ implies $X_1 = X$ or $X_1 = \{0\}$.*

*The minimality may vary with different related associative rings. On similar lines we may define Smarandache maximal pseudo ideals (S-maximal pseudo ideals).*

*Let RG be a groupoid ring we say X is a Smarandache cyclic pseudo ideal (S-cyclic pseudo ideal) related to the associative ring A if X is generated by a single element. We say the S-pseudo ideal X of a groupoid ring RG related to the associative ring A, $A \subset RG$ is Smarandache prime pseudo ideal (S-prime pseudo ideal) related to A if $xy \in I$ implies $x \in I$ or $y \in I$.*

The study of groupoid rings in mathematical literature is completely absent so the study of these ideals S-ideals, pseudo ideals and S-pseudo ideals will be very innovative and interesting.

The first open problem is, will the set of S-ideals of a finite groupoid ring of any groupoid G over a ring R form a modular lattice. Atleast study groupoids in $Z^{***}(G)$ over $Z_p$. $Z_p$ may be a prime field of characteristic p or a ring with p-elements that is p is a composite number.

The study of A.C.C. and D.C.C conditions on groupoid rings is yet another interesting feature for in case of group ring such study have been made. Several unsolved problems exist even in case of group rings to be Noetherian. It has been proved that KG is Artinian if and only if G is finite. Is it true, KG is Noetherian if and only if G has a series of subgroups?

$G = G_0 \supset G_1 \supset \ldots \supset G_n = \{1\}$ such that $G_{i+1}$ is normal in $G_i$ and $G_i /G_{i+1}$ is either a finite group or infinite cyclic.

Now the analogous study for groupoid rings remains untouched. Analogous to group rings which satisfy S-ACC or S-D CC we can define for the groupoid rings RG to satisfy S-ACC and S-D.C.C on S-ideals.



**DEFINITION 3.4.5:** *Let RG be a groupoid ring. If the set of all S-right ideals of RG is totally ordered by inclusion we say RG is Smarandache right chain groupoid ring (S-right chain groupoid ring).*

**DEFINITION 3.4.6:** *Let RG be a groupoid ring, A a S-subring of RG. If the set of all S-ideals of A in RG is totally ordered by inclusion then the groupoid ring RG is said to be a Smarandache weakly chain groupoid ring (S-weakly chain groupoid ring).*

**DEFINITION 3.4.7:** *Let RG be a groupoid ring. $I \neq 0$ be an ideal of RG . If for any nontrivial ideals X and Y of RG, $X \neq Y$ we have $\langle X \cap I, Y \cap I \rangle = \langle X, Y \rangle \cap I$ then I is called the obedient ideal of RG, $\langle \rangle$ denotes the ideal generated by the collection of elements in $X \cap I$ and $Y \cap I$.*

**DEFINITION 3.4.8:** *Let RG be a groupoid ring. If for any S-ideal I of RG we can find obedient S-ideals of RG say X, Y in RG, $X \neq Y$ such that $\langle X \cap I, Y \cap I \rangle = \langle X, Y \rangle \cap I$. Then we say I is a Smarandache obedient ideal of RG (S-obedient ideal). If in a groupoid ring RG we have every S-ideal I of RG to be a S-obedient ideal of RG then we say RG is a Smarandache ideally obedient ring (S-ideally obedient ring).*

Now we give a nice result about groupoid rings, which are not true in case loop rings or group rings.

**THEOREM 3.4.4:** *Let RG be a groupoid ring. If G has no unit and no subgroupoids and R is a S-ideally obedient ring, then RG is an S-ideally obedient ring.*

*Proof***:** We know every $I \neq 0$ in R is a S-obedient ideal, so is every $I G \neq 0$ in RG as G has no sub groupoids. Hence the theorem.

***Example 3.4.6*:** Let $Z_{12} = \{0, 1, 2, \ldots, 11\}$ be the ring of integers modulo 12. Let G be a groupoid given by the following table:

| *   | $a_0$ | $a_1$ | $a_2$ | $a_3$ | $a_4$ |
|-----|-------|-------|-------|-------|-------|
| $a_0$ | $a_0$ | $a_4$ | $a_3$ | $a_2$ | $a_1$ |
| $a_1$ | $a_2$ | $a_1$ | $a_0$ | $a_4$ | $a_3$ |
| $a_2$ | $a_4$ | $a_3$ | $a_2$ | $a_1$ | $a_6$ |
| $a_3$ | $a_1$ | $a_0$ | $a_4$ | $a_3$ | $a_2$ |
| $a_4$ | $a_3$ | $a_2$ | $a_1$ | $a_0$ | $a_4$ |

Clearly the groupoid ring $Z_{12}$ G has S-obedient ideals as $Z_{12}$ has I={0,6} to be an obedient ideal of $Z_{12}$. When we try to find non-associative Lin rings we will not be in a position to say nontrivially all Lie algebras are Lin rings as $(xy - yx)^n = xy - yx$ or $(xy + yx)^n = xy + yx$.

**DEFINITION 3.4.9:** *Let RG be a groupoid ring. We say RG is a Lin ring if $(xy - yx)^n = xy - yx$ or $(xy + yx)^n = xy + yx$ for every x , y in R and for some $n = n(x, y) > 1$.*



**THEOREM 3.4.5:** *Let $Z_2G$ be a groupoid ring where $Z_2 = \{0,1\}$ and $G$ a commutative groupoid with unit in which every $g \in G$ is such that $g^2 = 1$ is a Lin Ring.*

*Proof*: Obvious by the very definition of Lin ring.

**DEFINITION 3.4.10:** *Let $RG$ be a groupoid ring. We say $RG$ is an ideally strong ring if every subring of $RG$ not containing identity is an ideal of $RG$. We call $RG$ a Smarandache ideally strong ring (S-ideally strong ring) if every S-subring of $RG$ is an S-ideal of $RG$.*

We define $I^*$ rings in case of groupoid rings.

**DEFINITION 3.4.11:** *Let $RG$ be a groupoid ring. $\{I_j\}$ be the collection of all ideals of $RG$. $RG$ is said to be a $I^*$ ring if every pair of ideals $I_1, I_2 \in \{I_j\}$ in $RG$ and for every $a \in RG \setminus (I_1 \cup I_2)$ we have $\langle a \cup I_1 \rangle = \langle a \cup I_2 \rangle$ where $\langle \rangle$ denotes the ideal generated by $a$ and $I_j$ ; $j = 1, 2$.*

*Now we define Smarandache $I^*$ ring (S-$I^*$-ring) if for all S-ideals $\{A_i\}$ of $RG$; we have for every pair of ideals $A_1, A_2 \in \{A_i\}$; and for every $x \in RG \setminus \{A_i \cup A_j\}$, $\langle A_1 \cup x \rangle = \langle A_2 \cup x \rangle$ and they generate S-ideals of $RG$. We call the ring to be a Smarandache weakly $I^*$ ring (S-weakly $I^*$-ring) if for some $x \in R \setminus \{A_1 \cup A_2\}$ we have $\langle A_1 \cup x \rangle = \langle A_2 \cup x \rangle$.*

**THEOREM 3.4.6:** *All S-$I^*$ rings are S-weakly $I^*$ rings.*

*Proof:* Follows from the very definition.

The concept of Q-rings in case of groupoid rings and all the more in case of non-associative rings remains practically absent. So now we introduce Q-ring only in case of groupoid rings.

**DEFINITION 3.4.12:** *Let $RG$ and $R_1G_1$ be any two non-isomorphic finite groupoid rings, if there exists non-maximal ideals $I$ if $RG$ and $J$ of $R_1 G_1$ such that $RG/I$ is isomorphic to $R_1 G_1/ J$. Then the finite groupoid rings are said to be Q-groupoid rings.*

*Example 3.4.7:* $Z_2 G$ and $Z_4 G_1$ be two finite groupoid rings where $G \in Z_6$ (2,3) and $G_1 \in Z_4$ (1,3). Can $Z_2G$ and $Z_4G_1$ be Q-groupoid rings. We define further the concept of weakly Q-ring.

**DEFINITION 3.4.13:** *Let $RG$ be a groupoid ring such that all of its ideals are maximal and if we have $RG/\{0\}$ is isomorphic to some groupoid ring $R_1G_1$; then we call $RG$ a weakly Q-ring.*

**DEFINITION 3.4.14:** *Let $RG$ be a groupoid ring and $A$, a S-ideal of $RG$. $RG/A$ is defined as the Smarandache quotient ring (S-quotient ring) related to the S-ideal $A$ of $RG$.*



**DEFINITION 3.4.15:** *Let RG and $R_1G_1$ be any two groupoid rings, if we have S-ideals A and B of RG and $R_1G_1$ respectively such that the S quotient ring RG/A is S-isomorphic with the S-quotient ring $R_1G_1$/B then we say RG is a Smarandache Q-ring (S-Q-ring).*

Several examples and results can be got using the groupoids from the new class of groupoids $Z^{***}(n)$ and using the modulo ring of integers $Z_n$.

**DEFINITION 3.4.16:** *Let RG be a groupoid ring. RG is called an F-ring if there exists a finite set X of non-zero elements in RG such that a RG $\cap$ X $\neq \phi$ for any non-zero a in RG.*

*If in addition X is contained in the center of RG, RG is called a FZ-ring.*

*For the groupoid ring RG we say a S-subring A of RG is a Smarandache F-ring (S-F-ring) if we have a subset X in RG and a non-zero b $\in$ RG such that bA $\cap$ X $\neq \phi$.*

*Here it is pertinent to mention that we need not take X as a subset of A but nothing is lost even if we take X to be a subset of A. Similarly b can be in A or in RG. In case of associative rings one defines the concept of $\gamma_n$-rings.*

**DEFINITION 3.4.17:** *Let RG be a groupoid ring. RG is said to be a $\gamma_n$-ring if for n >1, n an integer $x^n - x$ is an idempotent for all $x \in R$.*

**DEFINITION 3.4.18:** *Let RG be a groupoid ring. RG is a Smarandache $\gamma_n$-ring (S-$\gamma_n$-ring) if for every $x \in$ RG, $x^n - x$ is a S-idempotent for some integer n >1.*

**THEOREM 3.4.7:** *If RG is a Smarandache $\gamma_n$-ring then RG is a $\gamma_n$-ring.*

*Proof:* By the very definitions the result follows.

**THEOREM 3.4.8:** *Let $Z_2G$ be the groupoid ring where G is a commutative groupoid with unit in which every $g \in$ G is such that $g^2 = 1$. Then $Z_2G$ is not a $\gamma_n$-ring.*

*Proof:* Follows from the fact that RG has elements $\alpha \neq 0$ such that $\alpha^2 = 0$ so $\alpha^n - \alpha$ can never be an idempotent. Hence the claim.

**DEFINITION 3.4.19:** *Let RG be a groupoid ring. A non-empty subset S of RG is called a closed net of RG if S is a closed set of RG under '.' and is generated by a single element. That is S is a groupoid.*

*If the groupoid ring RG is contained in a finite union of closed nets of RG then we say RG has a closed net.*

*If RG is a groupoid ring RG is a CN-ring if RG = $\cup S_i$ where $S_i$ 's are closed nets such that $S_i \cap S_J$'s are closed nets such that $S_i \cap S_J = \phi$ or {1} or {0} . if $i \neq j$ and 1 $\in$ RG and $S_i \cap S_J = S_i$ if i = j and each $S_i$ is a nontrivial closed net of R.*

*RG is said to be a Smarandache closed net if*



1. S is a groupoid
2. S is a S-groupoid

**PROBLEMS:**

1. Find all SG* for G = $Z_9$ (3,4) over the prime field $Z_7$.
2. Find G* for G = $Z_6$ (2,5) over the prime field $Z_3$. Is SG* = G*?
3. Whether there exists a relation between G* and atleast one of the SG*?
4. Let G = $Z_{11}$ (5,7) be a groupoid. Does ZG have S-pseudo ideals?
5. Give an example of a groupoid ring, which has no S-pseudo ideals.
6. Does the groupoid ring ZG where G = $Z_{12}$ (5,7) satisfy A.C.C. or D.C.C or S.A.C.C or S.D.C.C.
7. Give an example of a groupoid ring which satifies A.C.C and not S.A.C.C.
8. Give an example of a groupoid ring which satisfies both S.D.C.C and D.C.C on S-ideals.
9. Find the set of all S-ideals of the groupoid ring $Z_{12}$G where G = $Z_{16}$ (5,12).
10. Give an example of a groupoid ring, which is a $\gamma_n$-ring.
11. Is the groupoid ring $Z_6$G where G = $Z_{10}$ (5, 3) a $\gamma_n$- ring? Justify your answer.
12. Give an example of a groupoid ring, which is a F-ring.
13. Can the groupoid ring $Z_3$G where G = $Z_6$(2,3) be a F-ring? Justify your claim.
14. Give an example of a groupoid ring which is a CN-ring.
15. Is the groupoid ring $Z_3$G where G = $Z_7$ (3,6) a CN –ring?
16. Can the groupoid ring given in example (15) be a S-CN –ring?
17. Does the groupoid ring $Z_8$G where G = $Z_{12}$ (4,8) have a S-obedient ideal?
18. Will the groupoid ring $Z_6$G where G = $Z_2$ (3,9) have a obedient ideal?

## 3.5 Special Properties in Groupoid rings and Smarandache Groupoid rings

This section is solely devoted to the introduction of several interesting properties in groupoid rings. The main concepts introduced in this section are the n-capacitor sub groups, essential and S-essential groupoid rings, exponential groupoid rings and Smarandache exponential rings and groupoid rings satisfying special identity (xy) x = x(yx). Finally we introduce the concept of subring link relation in groupoid rings. Radicals and S-radicals are introduced and studied in case of groupoid ring RG. These studies are carried out both in case of finite and infinite groupoid rings.

Here we introduce the concept of a new class of groupoids using [-∞, ∞].

**DEFINITION 3.5.1:** *Let L = [-∞, ∞]. Define for a , b ∈ L, a ∗ b = ma + nb where m, n ∈ L. Clearly {L, ∗, (m, n)} is a groupoid of infinite order. For varying m, n in L we get an infinite class of infinite groupoids.*

**Result 1:** [L, ∗, (1, 1)] is a semigroup.

**Result 2:** The groupoids [L, ∗, (m, n)] ≠ [L, ∗, (n, m)] when m ≠ n.



**Result 3:** The groupoid [L, ∗, (m, n)] satisfies the identity $(x * y) * (z * y) = x * z$ if and only if m = 1 and n = 0 or m = 1 and n = -1.

In view of these results we can study the groupoid ring RL where L is a groupoid built using [-∞, ∞].

**THEOREM 3.5.1:** *Let $L = \{[-\infty, \infty], *, (m, n) \mid m, n \in [-\infty, \infty]\}$ be an infinite groupoid. $Z_2 = \{0, 1\}$ be the prime field of characteristic two. $Z_2L$ is an infinite non-commutative groupoid ring if $m \neq n$ of characteristic two and has every element to be of order 2 if $m = -n$.*

*Proof:* Straightforward by the very definition.

It is left for the author to built interesting and innovative results in this direction by incorporating the properties of Smarandache substructures introduced in section 4 and Smarandache special elements in section 3.

**_Remark_:** Using these groupoids L and choosing any $Z_p$ ($p < \infty$), p a prime. $Z_p L$ will be an infinite groupoid ring of characteristic p. If both m and n are chosen to be positive we see every element in L is torsion free and the groupoid L will be a torsion free non-abelian groupoid.

**DEFINITION 3.5.2:** *Let RG be a groupoid ring. We say $a \in RG$ ($a \neq 0$) is said to be a normal element of RG if a RG = RG a. If every element $a \in R G$ is a normal element then we call RG a normal groupoid ring.*

We give nice condition for a groupoid ring to be a normal groupoid ring.

**THEOREM 3.5.2:** *Let K be any field and G a groupoid such that G is a normal groupoid then the groupoid ring KG is a normal groupoid ring.*

*Proof:* Straightforward by the very definition.

**DEFINITION 3.5.3:** *Let KG be a groupoid ring We say KG is a Smarandache normal groupoid (S-normal groupoid) ring if KG has a proper subset A, $A \subset KG$ and A a SNA-subring of KG, such that for all $\alpha \in K G$, $\alpha A = A \alpha$.*

Thus we have the immediate consequence.

**THEOREM 3.5.3:** *Suppose KG is a S-normal groupoid then KG need not in general be a normal groupoid.*

*Proof:* The reader is requested to construct an example to prove this result.

The natural question would be if KG is a normal groupoid will KG be a S-normal groupoid, even this question is not that easy to be settled.



**THEOREM 3.5.4:** *Let G be a groupoid. K any field, KG the groupoid ring. $g \in G$ is a normal element of KG if and only if $gG = Gg$.*

*Proof:* One way is obvious if $gG = Gg$ then $gKG = KGg$.

Suppose $g \in G$ is such that $gKG = KGg$ then we must have $gG = Gg$ for otherwise we will not have $gKG = KGg$ as $G \subset KG$. It is interesting to note all elements of G are not normal elements of KG for if $gG \neq Gg$ then $gKG \neq KGg$. This is clear by the example. For take, $Z_2 = \{0,1\}$ and G a groupoid given by the following table:

| . | $a_1$ | $a_2$ | $a_3$ | $a_4$ |
|---|---|---|---|---|
| $a_1$ | $a_2$ | $a_2$ | $a_2$ | $a_2$ |
| $a_2$ | $a_4$ | $a_4$ | $a_4$ | $a_4$ |
| $a_3$ | $a_2$ | $a_2$ | $a_2$ | $a_2$ |
| $a_4$ | $a_4$ | $a_4$ | $a_4$ | $a_4$ |

Clearly in the groupoid ring $Z_2G$ we have $a_1 Z_2G \neq Z_2G\, a_1$ as $a_1 G \neq Ga_1$. Hence the claim.

**DEFINITION 3.5.4:** *Let KG be a groupoid ring $N(KG) = \{$ set of all elements in KG such that $\alpha KG = KG\alpha$ where $\alpha \in KG\}$ that is $N(KG)$ is the set of normal elements in KG.*

**THEOREM 3.5.5:** *Let KG be a groupoid ring. $N(KG)$ is not even closed under,'.' .*

*Proof:* By an example. Let K be any field and G be a groupoid given by the following table:

| * | a | b | c |
|---|---|---|---|
| a | b | c | a |
| b | a | b | c |
| c | c | c | c |

$aKG = KGa$ for $a \in G \subset KG$ but $a.b = b \notin N(KG)$ as $b.KG \neq KGb$. Hence the claim.

Now we study the condition for a groupoid ring RG to satisfy the identify $(\alpha \beta)\alpha = \alpha(\beta\alpha)$ for all $\alpha, \beta \in RG$.

*Example 3.5.1:* Consider the groupoid ring $Z_2G$ where G is given by the following table:

| . | $a_0$ | $a_1$ | $a_1$ | $a_3$ | $a_4$ |
|---|---|---|---|---|---|
| $a_0$ | $a_0$ | $a_3$ | $a_1$ | $a_4$ | $a_2$ |
| $a_1$ | $a_3$ | $a_1$ | $a_4$ | $a_2$ | $a_0$ |
| $a_2$ | $a_1$ | $a_4$ | $a_2$ | $a_0$ | $a_3$ |
| $a_3$ | $a_4$ | $a_2$ | $a_0$ | $a_3$ | $a_1$ |
| $a_4$ | $a_2$ | $a_0$ | $a_3$ | $a_1$ | $a_4$ |



It is easily verified, the groupoid ring $Z_2G$ satisfies the identify $(\alpha \beta)\alpha = \alpha(\alpha\beta)$ for all $\alpha, \beta \in Z_2G$.

*Example 3.5.2:* Let $Z_2G$ be a groupoid ring where G is a groupoid given by the following table:

| . | $a_0$ | $a_1$ | $a_2$ | $a_3$ | $a_4$ |
|---|---|---|---|---|---|
| $a_0$ | $a_0$ | $a_2$ | $a_4$ | $a_1$ | $a_3$ |
| $a_1$ | $a_3$ | $a_0$ | $a_2$ | $a_4$ | $a_1$ |
| $a_2$ | $a_1$ | $a_3$ | $a_0$ | $a_2$ | $a_4$ |
| $a_3$ | $a_4$ | $a_1$ | $a_3$ | $a_0$ | $a_2$ |
| $a_4$ | $a_2$ | $a_4$ | $a_1$ | $a_3$ | $a_0$ |

Clearly the groupoid ring does not satisfy the identity for $(a_2 \ a_4) \ a_2 \neq a_2 \ (a_a \ a_2)$ for $a_2$, $a_4 \in G$; so for all $\alpha, \beta \in Z_2G$ we see $(\alpha\beta) \ \alpha \neq (\beta\alpha)$. Hence all groupoid rings in general do not satisfy the identity.

[46] has introduced the concept of exponentiation in a ring R. He says the triple (R, B, E) is an exponential ring if for the ring R, B is the multiplicative subsemigroup of (R, .) (B the basis) which does not contain zero and E (the exponents) be a semiring with unit element 1; that is (E, +) is a semigroup and (E,.) is a semigroup with 1 and '.' distributes over '+' both from the left and from the right with a binary operation from

$B \times E \to B \subseteq R$ (that is $(b, e) \to b^e$) such that for all $b, d \in B$ and $e, k \in E$, $b^e d^e = (bd)^e$. $b^{ek} = (b^e)^k$, $b^{e+k} = b^e \ b^k$ and $b^1 = b$.

Using this concept we build non-associative rings by replacing the associative rings which we call as the first type of exponentiation. The second type of exponentiation would be by replacing in the triple (R, B, E); E the exponents by a groupoid under + and a groupoid under '.' with 1 and + and '.' distributors over each other.

For more about exponentiation please refer [46]

**DEFINITION 3.5.5:** *Let A = RG be a groupoid ring of a groupoid G with 1. Let P be a multiplicative groupoid of (RG,.) which does not contain zero and let E (the exponents) be a semi ring with unit 1.*

*A binary operation '$\times$' from $P \times E \to P \subseteq A = RG$ is $((p,e) \to p^e)$ makes (RG, P, E) a non-associative ring if it has the usual properties of exponentiation viz. for all $b, d \in P$ and $e, k \in E$, $b^e d^e = (bd)^e \ b^{ek} = (b^e)^k$, $b^{e+k} = b^e \ b^k$ and $b^1 = b$. Thus the exponentiation in a non-associative ring is just defined from the exponentiation of a ring by replacing the semigroup by a groupoid.*

*Example 3.5.3:* Let R be a commutative ring with 1 or a field. L a loop so that L is trivially a groupoid. RL be the loop ring. Take P a multiplicative groupoid from RL



where P does not contain 0. Take E exponents and (RL, P, E) is the exponential ring of the non-associative ring RL.

Clearly (R, S, E) where R is a ring and S is a subsemigroup of (R, .) without zero is the exponential subring of (RL, P, E). Thus we can build exponential non-associative rings from exponential associative rings. If we replace the loop by a groupoid with 1 still one gets a non-associative ring on which exponentiation can be defined.

**DEFINITION 3.5.6:** *Let (RG, P, E) be a non-associative exponentiation ring. We say (RG, P, E) is a Smarandache non-associative exponentiation ring (S-non-associative exponentiation ring) if this has a sub triple (R, S, E), which is an associative ring with exponentiation.*

Thus we have got the following theorem the proof of which is straightforward.

**THEOREM 3.5.6:** *Let RG be a groupoid ring (where G is a groupoid with 1) or a loop ring. Then the triple (RG, P, E) is a S-non-associative exponentiation ring.*

We define n-capacitor for commutative groupoid rings.

**DEFINITION 3.5.7:** *Let RG be a commutative groupoid ring P an additive subgroup of RG, P is called the n-capacitor group of RG if $x^n P \subset P$ for every $x \in R$ and $n > 1$ and n a positive integer.*

*Example 3.5.4:* Let $Z_2G$ be a commutative groupoid ring where the groupoid G is given by the following table:

| *     | $a_0$ | $a_1$ | $a_2$ | $a_3$ | $a_4$ | $a_5$ |
|-------|-------|-------|-------|-------|-------|-------|
| $a_0$ | $a_0$ | $a_2$ | $a_4$ | $a_0$ | $a_2$ | $a_4$ |
| $a_1$ | $a_2$ | $a_4$ | $a_0$ | $a_2$ | $a_4$ | $a_0$ |
| $a_2$ | $a_4$ | $a_0$ | $a_2$ | $a_4$ | $a_0$ | $a_2$ |
| $a_3$ | $a_0$ | $a_2$ | $a_4$ | $a_0$ | $a_2$ | $a_4$ |
| $a_4$ | $a_2$ | $a_4$ | $a_0$ | $a_2$ | $a_4$ | $a_0$ |
| $a_5$ | $a_4$ | $a_0$ | $a_2$ | $a_4$ | $a_0$ | $a_2$ |

$Z_2G$ is a commutative groupoid ring. Let $P = Z_2 \{a_0, a_2, a_4\}$. It is easily verified P is a n capacitor group of $Z_2G$.

**DEFINITION 3.5.8:** *Let RG be a commutative groupoid ring and P an additive S-semigroup of RG. P is called a Smarandache n-capacitor group (S-n-capacitor group) of RG if $x^n P \subseteq P$ for every $x \in RG$ and $n \geq 1$; n a positive integer.*

We now proceed on to define radical for groupoid rings.

**DEFINITION 3.5.9:** *Let RG be a groupoid ring. A subset P of RG is a radical if*

       *(1). P is an ideal of RG*
       *(2). P is a nil ideal*
       *(3). R/P has no non-zero nilpotent right ideals.*



*The sum of all ideals in RG satisfying (1) and (2) is the upper radical of RG and is denoted by U (RG). The intersection of all these ideals in RG satisfying (1) and (3) is the lower radical of RG denoted L(RG).*

**DEFINITION 3.5.10:** *Let RG be a groupoid ring. The Smarandache radical ideal (S-radical ideal) P of RG is defined as follows:*

1. *P is a S-ideal of RG*
2. *$S \subset P$ where S is a subideal of P is a nil ideal*
3. *RG/P has no non-zero nilpotent right ideals.*

*The sum of all S-ideals of RG satisfying (1) and (2) is called the Smarandache upper radical (S-upper radical) of RG and is denoted by S(U (RG)). The intersection of those S-ideals in RG satisfying (1) and (3) is the Smarandache lower radical (S-lower radical) of RG denoted by S(L(RG)).*

Now we proceed on to define the subring link relation in groupoid rings.

**DEFINITION 3.5.11:** *Let RG be a groupoid ring. A pair of elements $x, y \in RG$ is said to have a subring right link relation if there exists a subring M of RG in RG\ {x, y} that is $M \subseteq RG \setminus \{x, y\}$ such that $x \in My$ and $y \in Mx$. Similarly subring left link relation if $x \in yM$ and $y \in xM$. If it has both a left and a right link relation for the same subring M then we have x and y to have a subring link relation and is denoted by xMy.*

**DEFINITION 3.5.12**: *Let RG be a groupoid ring. We say a pair x, y in RG has a weakly subring link with a subring P in RG\ {x , y} if either $y \in Px$ or $x \in Py$ 'or' in the strictly mutually exclusive sense and we have subring Q, $Q \neq P$ such that $y \in Q x$ (or $x \in Q y$).*

*We say a pair x, y in RG is one way weakly subring link related if we have a subring $P \subset RG \setminus \{x, y\}$ such that $x \in Py$ and for no subring $S \subset RG \setminus \{x, y\}$ we have $y \in Sx$. Let RG be a groupoid ring a pair $x, y \in R G$ is said to have a Smarandache subring link relation (S-subring left link relation) if there exists a S-subring P in RG\{x , y} such that $x \in Py$ and $y \in xP$. If it has both a Smarandache left and right link relation for the same S-subring P then we say x and y have a Smarandache subring link (S-subring link).*

*We say x, y $\in$ RG is a Smarandache weak subring link (S-weak subring link) with a S-subring P in RG \ {x, y} if either $x \in Py$ or $y \in Px$ (or in strictly mutually exclusive sense) we have a S-subring $Q \neq P$ such that $y \in Q x$ (or $x \in Qy$).*

*We say pair x, y $\in$ RG is said to be Smarandache one way weakly subring link related (S-one way weakly link related) if we have a S-subring $P \subset RG \setminus\{x, y\}$ such that $x \in Py$ and for no subring $Q \subseteq RG \setminus\{x, y\}$ we have $y \in Qx$.*

Let us now define essential subring of the groupoid ring RG.



**DEFINITION 3.5.13:** *Let RG be a groupoid ring; A S-subring x of RG is said to be a Smarandache essential subring (S-essential subring) of RG if the intersection of X with every other S-subring is zero. If every S-subring of RG is an S-essential subring of RG then we call RG an S-essential ring.*

Study of essential subrings, S-essential subrings and groupoid rings, which are essential rings, happens to be an interesting field.

*Example 3.5.5:* $RG = Z_2G$ be the groupoid ring where G is the groupoid given by the following table:

| . | $a_0$ | $a_1$ | $a_2$ | $a_3$ |
|---|---|---|---|---|
| $a_0$ | $a_0$ | $a_2$ | $a_0$ | $a_2$ |
| $a_1$ | $a_2$ | $a_0$ | $a_2$ | $a_0$ |
| $a_2$ | $a_0$ | $a_2$ | $a_0$ | $a_2$ |
| $a_3$ | $a_2$ | $a_0$ | $a_2$ | $a_0$ |

The groupoid ring $Z_2$ G has $Z_2[a_0]$, $Z_2[a_0\ a_2]$, $Z_2[a_3, a_2, a_0]$ and $Z_2[a_0\ a_1\ a_2]$ to be subgroupoid rings. We see $Z_2G$ is not an essential ring or does not even have a single essential subring.

**PROBLEMS:**

1. Construct an infinite groupoid ring of characteristic two using $Z^+ = \{1, 2, 3, …\}$.
2. Give an example of a normal groupoid ring.
3. Find an example of a S-normal groupoid ring.
4. Characterize a class of groupoid rings RG which satisfy the identify $(\alpha\beta)\alpha = \alpha(\beta\alpha)$ for $\alpha, \beta \in RG$.
5. Is the groupoid ring $Z_7G$ be made into an exponential ring where $G \in Z***(9)$.
6. Construct rings, which are not S-exponential rings.
7. Characterize those groupoid rings, RG which are n-capacitor rings G taken from $Z***(p)$ and $R = Z_2$.
8. Find the radical and S- radical for the groupoid rings $Z_mG$ where $G \in Z***(n)$ where
    a. $(m, n) = 1$.
    b. $(m, n) = d \neq 1$ or $n$
    c. $(m, n) = n$.
    d. $m = n$.
9. Classify those groupoid rings R G which are
    a. essential or
    b. S-essential.
   when $G \in Z**(7)$ and $R = Z_{12}$.
10. Characterize those groupoid rings RG which are never essential or S-essential when $R= Z_2$ and $G \in Z*(6)$.



**Chapter 4**

# LIE ALGEBRAS AND SMARANDACHE LIE ALGEBRAS

In this chapter we recall the definition of Lie algebras and its properties and introduce the concept of Smarandache Lie algebras. This Chapter has five sections. In the first section we introduce the basic notions of Lie algebras. Section two introduces Smarandache Lie algebras and basic properties of Lie algebras. Section three is devoted to defining several new notions in Smarandache Lie algebras and Lie algebras. Section four is devoted to the study of certain special properties and in the final section we introduce some new concepts in Lie algebras and S-Lie algebras and define Smarandache mixed direct product of Lie algebras of type A and B. As every section has an introduction we do not recall the contents of each section here.

## 4.1 Basic Properties of Lie algebras

In this section we just recall the definition of real and complex Lie algebras from [36, 74] and some of the basic properties enjoyed by them: but we do not claim that we have fully exhausted recalling all the basic properties of the Lie algebras. We have enlisted only some of the properties, which has interested us. We leave it for the reader to refer books on Lie algebras for more literature.

**DEFINITION 4.1.1:** *A real Lie algebra L is a vector space over the real field furnished with a bilinear product $[\lambda, \mu]$ satisfying the identities.*

$$[\lambda, \mu] = - [\lambda, \mu]$$
$$\text{and } [ [\lambda, \mu], \nu] + [ [\mu, \nu], \lambda] + [[\nu, \lambda], \mu] = 0$$

*[Jacobi identity]*

In what follows we will always suppose unless the contrary is explicitly specified, that a Lie algebra is finite dimensional. Now we just recall the definition of complex Lie algebra.

**DEFINITION 4.1.2:** *Let L be a real Lie algebra. Let $\overline{L}$ be the set of all formal sums $\lambda + i\mu$; $\lambda, \mu \in L$. We define addition and multiplication by complex numbers and a product $[\lambda + i\mu, \lambda' + i\mu']$ by the formula $(\lambda + i\mu) + (\lambda' + i\mu') = \lambda + \lambda' + i(\mu + \mu')$ $(x + iy) (\lambda + i\mu) = (x\lambda - y\mu) + i(y\lambda + x\mu), [\lambda + i\mu, \lambda' + i\mu'] = ([\lambda'\lambda'] - [\mu' \mu']) + i ([\lambda, \mu'] + [\mu, \lambda']).$*

It is easily verified that with these definition L becomes a complex Lie algebra. The complex Lie algebra $\overline{L}$ is called the complexification of the real Lie algebra L.

Thus we can also say non-associative algebra L is said to be associative if its multiplication satisfies the associative law that is (xy) z = x (yz).



A non-associative algebra L is a Lie algebra, if its multiplication satisfies the Lie conditions i.e., $x^2 = 0$, $(xy)z + (yz)x + (zx)y = 0$, that is the Jacobi identify is true.

Clearly if L is a Lie algebra and x, y ∈ L then $0 = (x + y)^2 = x^2 + xy + yx + y^2 = 0$ so that xy = - yx. Conversely if this condition holds $2x^2 = 0$ so that if the characteristic of L is not two then $x^2 = 0$. Thus for characteristic not two we can assume instead of $x^2 = 0$; xy + yx = 0 or xy = - yx. In general if two elements x,y ∈ L is such that [x, y] = 0 then we say x and y commute . The Lie algebra L is abelian if [ x , y ] = 0 for all x , y ∈ L.

*Example 4.1.1*: Let $C^3$ be a three dimensional vector space over C and $e_1$, $e_2$, $e_3$ be a basis of C. For any two elements $x = x_1e_1 + x_2e_2 + x_3e_3$, and $y = y_1e_1 + y_2e_2 + y_3e_3$, define $[x, y] = (x_2y_3 – x_3y_2) e_1 + (x_3y_1 – y_3x_1) e_2 + (x_1y_2 – x_2y_1)e_3$. Then $C^3$ becomes a Lie algebra.

*Example 4.1.2:* Let S denote the collection of all 3 × 3 matrices with entries from C. For any X, Y ∈ S define [X, Y] = XY – YX, then S is a Lie algebra.

**DEFINITION 4.1.3:** *Let L and $L_1$ be any two Lie algebras. A one to one map from L onto $L_1$ is called a Lie algebra isomorphism if it satisfies:*

*(1) if $X_1 \to Y_1$ , $X_2 \to Y_2$ then for any $\lambda$ , $\mu \in$ C (where L and $L_1$ are Lie algebras over C) $\lambda X_1 + \mu X_2 \to \lambda Y_1 + \mu Y_2$.*

*(2) if $X_1 \to Y_1$ , $X_2 \to Y_2$ then
$[X_1 , X_2 ] \to [Y_1 , Y_2]$.*

*We say L and $L_1$ are isomorphic and denote it by $L_1 \cong L_2$. In particular an isomorphism from L onto itself is called an automorphism.*

One of the basic fundamental problems in Lie algebras is to determine all non-isomorphic Lie algebras. The following is a very important and commonly used example.

*Example. 4.1.3:* Let gl (n, C) be the set of all n × n matrices. It is well known with respect to matrix addition and scalar multiplication gl(n, C) forms $n^2$ dimensional vector space.

Now for any X , Y ∈ gl(n, C), define [X, Y] = XY – YX, then gl (n, C) forms a Lie algebra. gl(n, C) can be considered as the set of all linear transformation of an n-dimensional vector space V, then in this case it is denoted by gl (V).

**DEFINITION 4.1.4:** *Let L be a Lie algebra. M and N denote subsets of L. Denote by M + N the linear subspace spanned by elements of the form m + n (m ∈ m, n ∈ N ) and by [M, N] the subspace spanned by elements of the form [m , n] (m ∈ M and n ∈ N) . If M, $M_1$, $M_2$, $N_1$P are subspaces of L then*

*(1) $[M_1 + M_2, N] \subseteq [M_1, N] + [M_2, N]$.*
*(2) [M, N] = [N, M].*



$$(3)\ [M, [N, P]] \subseteq [N, [P, M]] + [P, [M,N]].$$

*A subspace H of L is said to be a subalgebra if $[H,H] \subseteq H$, i.e. $X, Y \in H$ implies that $[X, Y] \in H$. A subspace H of L is said to be an ideal if $[L, H] \subseteq L$ i.e. $X \in L$ and $Y \in H$ implies $[X, Y] \in H$.*

*An ideal is a subalgebra if $H_1$ and $H_2$ are ideals then $H_1 + H_2$ and $H_1 \cap H_2$ are also ideals of L. Subalgebras of gl(n, C) are called matrix Lie algebras or linear Lie algebras. If H is an ideal of L, then the quotient L/H which consists of all cosets (congruence classes mod H) is defined. For $X \in L$ denote the congruence class containing X by $\overline{X}$. Define $[X, Y] = [\overline{X}, \overline{Y}]$ the definition is independent of X and Y.*

*The quotient space L/H thus becomes a Lie algebra, this algebra is called the quotient algebra of L with respect to H. If L is a Lie algebra and H be an ideal of L, then it can be proved that the mapping $X \to \overline{X}$ satisfies the conditions.*

*(1) if $X \to \overline{X}$, $Y \to \overline{Y}$ then for any $\lambda, \mu \in C$, $\lambda X + \mu X \to \lambda \overline{X} + \mu \overline{Y}$.*
*(2) if $X \to \overline{X}$, $Y \to \overline{Y}$ then $[X, Y] \to [\overline{X}, \overline{Y}]$.*

*In general a mapping $X \to X_1$ from a Lie algebra L to Lie algebra $L_1$ is said to be a homomorphism if it satisfies;*

*(1) if $X \to X_1$, $Y \to Y_1$ then for any $\lambda, \mu \in C$, $\lambda X + \mu Y \to \lambda X_1 + \mu Y_1$.*
*(2) if $X \to X_1$, $Y \to Y_1$ then $[X, Y] \to [X_1, Y_1]$.*

*If this mapping is onto then $L_1$ is said to be a homomorphism image of L.*

***Example 4.1.4***: All trace zero matrices of gl(n, C) form a subalgebra, denote it by $A_{n-1}$. In fact $A_{n-1}$ is an ideal of gl(n, C) for if $X, Y \in$ gl(n, C) then Tr $[X, Y]$ = Tr $[XY - YX] = 0$, therefore $[X, Y] \in A_{n-1}$. All scalar matrices of gl(n, C) form a one dimensional subalgebra which is also an ideal of gl(n, C) for if $\lambda I$ is a scalar matrix then for any $x \in$ gl(n, C) we have $[X, \lambda I] = X\lambda I - \lambda IX = 0$. All diagonal matrices of gl(n, C) form an n-dimensional abelian subalgebra; denote it by d (n, C). The set of trace zero diagonal matrices form an (n – 1) dimensional abelian subalgebra of $A_{n-1}$.

**DEFINITION 4.1.5**: *Let L be a Lie algebra. Obviously L and {0} are ideals of L. If L does not have any other ideals, then it is said to be a simple Lie algebra.*

**DEFINITION 4.1.6**: *Let L be a Lie algebra and $M_1, \ldots, M_n$ be ideals of L, if any element X in L can be uniquely written as $X = X_1 + X_2 + \ldots + X_n$ where $X_i \in L$, i = 1, 2,…, n, then L is said to be the direct sum of $M_1, \ldots, M_n$. We also denote direct sum by $L = M_1 + \ldots + M_n$. If L is a direct sum of $M_1, M_2, \ldots, M_n$ then for $i \neq j$, $M_i \cap M_j = \{0\}$, for if $X \in M_i \cap M_j$, then $X = 0 + \ldots + X_i + \ldots + 0 = 0 + \ldots + X_j + \ldots + 0$ are two expressions of X and this forces X = 0. From $M_i \cap M_j = \{0\}$ ($i \neq j$), it follows that $[M_i, M_j] = \{0\}$ ($i \neq j$) are ideals.*



**DEFINITION 4.1.7**: *Let L be a Lie algebra, [L, L] is called the derived algebra of L and is denoted by DL. If H is an ideal of L then so is DH. In fact $[L, DH] = [L, [H, H]] \subseteq [H, [H, L]] + [H, [L, H]] \subseteq [H, H] + [H, H] = DH$. $D^0 L = L$, $D^1 L = DL$, $D^{n+1} L = D(D^n L)$. In this way, we obtain a series of ideals satisfying $D^0 L \supseteq D^1 L \supseteq \ldots \supseteq D^n L \supseteq \ldots$*

*This series is called the derived series of L. If there exists a positive integer n such that $D^n L = \{0\}$ then L is called a solvable Lie algebra. The following are true about solvable Lie algebras.*

*Subalgebras of solvable algebras are solvable, homomorphic images of solvable algebras are solvable, quotient algebras of solvable algebras are solvable. If L is a Lie algebras H an ideal and if H and L / H are solvable then L is solvable. Direct sum of solvable Lie algebras are solvable.*

*If an ideal H of a Lie algebra L is solvable then it is called a solvable ideal.*

*If L does not contain any solvable ideal except {0} then L is said to be semi simple.*

*An equivalent definition is that L is semi simple if it does not contain any non-zero abelian ideal. In fact non-zero abelian ideals are certainly solvable.*

*A semi simple Lie algebra is the direct sum of all the minimal ideals.*

*The one-dimensional Lie algebra is the only Lie algebra which is simple but not semi simple. So in this text the one-dimensional Lie algebra will be excluded from the simple Lie algebras, hence all simple Lie algebras are semi simple.*

*If L has a unique maximal solvable ideal $\tau$, $\tau$ is called the radical of L.*

**DEFINITION 4.1.8**: *A Lie algebra L is said to be nilpotent if there exists a positive integer n such that $C^n L = \{0\}$, where $C^0 L = L$, $C^1 L = [L, C^0 L]$, …, $C^{(n-1)} L = [L, C^n L] = \ldots$ If $C^{n-1} L$ is an ideal of L then $[L, C^n L] \subseteq [L, C^{(n-1)} L] = C^n L$, hence $C^n L$ is also an ideal of L. In this way we obtain a series of subalgebras satisfying $C^0 L \supseteq C^1 L \supseteq \ldots \supseteq C^{(n)} L \supseteq \ldots$ and each $C^{(n)} L$ is an ideal of L. This series is called the descending central series. If there exists a positive integer n such that $C^{(n)} L = \{0\}$, then L is said to be nilpotent, that is L is nilpotent for some positive integer n, $[X_n, [\ldots [X_2 X_1] \ldots]] = 0$ for any $X_1, \ldots, X_n \in L$. It can be proved $D^{(n)} L \subseteq C^{(n)} L$. In fact from $D^n L \subseteq C^n L$, it follows that $D^{(n+1)} L = [D^{(n)} L, D^{(n)} L] \subseteq [L, C^{(n)} L] = C^{(n+1)} L$.*

*Some of the interesting properties on nilpotent Lie algebras are:*

1. *If L is nilpotent then it is solvable.*

2. *Subalgebras of nilpotent Lie algebras are nilpotent. Homomorphic images of nilpotent Lie algebras are nilpotent in particular quotient algebras of nilpotent Lie algebras are nilpotent.*

*Direct sums of nilpotent algebras are nilpotent.*



***Example 4.1.5***: If t (n, C) denote the set of all upper triangular matrices of gl(n, C) i.e. matrices of the form

$$\begin{bmatrix} x_{11} & \cdots & & x_{1n} \\ 0 & x_{22} & \cdots & x_{2n} \\ \cdots & & & \\ 0 & \cdots & 0 & x_{nn} \end{bmatrix}$$

is a Lie algebra t(n, C) . Elements in t(n, C) that have equal diagonal entries form a Lie algebra m(n, C); t(n, C) is solvable and m(n, C) is nilpotent.

The image of the mapping A → adA is denoted by adL; the kernel of this mapping is the set of all A ∈ L satisfying [A, X] = 0 for all X ∈ L. These elements are called central elements.

For a linear transformation A in a finite dimensional vector space V the fitting lemma asserts that V = $V_{0A}$ ⊕ $V_{iA}$ where $V_{iA}$ are invariant relative to A and the induced transformation in $V_{0A}$ is nilpotent and in $V_{iA}$ it is an isomorphism , [36] calls $V_0A$ the fitting null component of V relative to A and $V_{iA}$ fitting one component of V relative to A. The space $V_{oA}$ = {x / x$A^i$ = 0 for some i} and

$$V_{iA} = \bigcap_{i=1}^{\infty} VA^i.$$

Let L be a Lie algebra of linear transformations in a finite dimensional vector space over a field of characteristic 0. An element h ∈ L is called regular if the dimensionality of the Fitting null component of L relative to adh is minimal. A Lie algebra L is a linear transformations of a finite dimensional vector space over a field of characteristic zero is called almost algebraic if it contains the nilpotent and semi simple components for every X ∈ L.

For more about Lie algebras refer [36, 74].

## 4.2 Smarandache Lie Algebras and its Basic Properties

In this section we introduce the concept of Smarandache Lie algebras Till date Smarandache Lie algebras were not defined. To define this concept we basically need the notion of Smarandache vector space, Smarandache linear transformation etc, but except for the definition of Smarandache vector space nothing has been developed about Smarandache linear transformation etc.

So in chapter one some basic definitions about Smarandache vector spaces have been introduced mainly keeping in mind the essential concepts used in the definition of Smarandache Lie algebras. Further we do not even try to define Lie groups for we felt the very study of Smarandache Lie groups can be taken as a separate one, hence we have defined only Smarandache Lie algebras and introduce some basic Smarandache properties about them. We do not promise to exhaust all notions about



Lie algebras and introduce them to Smarandache Lie algebras. Only those concepts which has interested us and which we felt would be developable to Smarandache concepts has been defined in this book.

**DEFINITION 4.2.1:** *A Smarandache real Lie algebra (S-real Lie algebra) L is a vector space over the real field which has a proper subset X ($X \subset L$, $X \neq \phi$) such that X is a Smarandache K- vectorial subspace of L furnished with a bilinear product [ $\lambda$ , $\mu$] satisfying the identities [ $\lambda$ , $\mu$] = - [$\mu$ ,$\lambda$ ], [[ $\lambda$ , $\mu$] , $\nu$] + [[,$\mu$, $\nu$] , $\lambda$] + [[$\nu$,$\lambda$],$\mu$] = 0 for all $\lambda$ , $\mu$ , $\nu$ $\in$ X. It is important to note that this bilinear product need not be defined on the totallity of L.*

*We always assume L to be a finite dimensional vector space over the field of reals.*

*This type of Smarandache Lie algebras will be called as Smarandache Lie algebras of type I, but we will for short say Smarandache Lie algebras. We define type II and other types in what follows. Similarly we define Smarandache complex Lie algebras (S-complex Lie algebra). All real or complex Lie algebras need not in general be S-real or complex Lie algebras.*

**THEOREM 4.2.1:** *Let L be a complex or a real Lie algebra; L is a S-complex or a real Lie algebra if and only if L has a S-vectorial subspace X.*

*Proof*: Straightforward from the very definitions and the fact that L is a complex or a real Lie algebra so satisfies the identities.

Now we define the conditions under which a non-associative ring is a Lie algebra.

**DEFINITION 4.2.2:** *Let V be a non-associative algebra. V is said to be a Smarandache Lie algebra (S-Lie algebra) if V has a proper subset $W \subset V$, $W \neq \phi$ such that W is a S-non-associative ring and every element of W satisfies the Lie conditions $x^2 = 0$ and (xy)z + (yz) x + (zx) y =0 for all x, y, z $\in$ W.*

**THEOREM 4.2.2:** *If V is a S-Lie algebra where V is a non-associative ring then V is a SNA-ring.*

*Proof*: Follows from the fact that if a non-associative ring has proper subset W, where W is a SNA-ring then obviously by the very definition of SNA-ring, V is a SNA-ring.

The following theorem is given as an exercise for the reader to prove.

**THEOREM 4.2.3:** *If V is a SNA-ring then V in general need not be a S-Lie algebra even if V is a Lie algebra.*

**Note:** The reader must construct an example of a SNA-ring V which is a Lie ring but V must not contain any SNA-subring.

***Example 4.2.1*:** Let $Z_3$ = {0, 1, 2} be the prime field of characteristic three. G be a groupoid given by the following table:



| * | $a_0$ | $a_1$ | $a_2$ | $a_3$ | $a_4$ | $a_5$ |
|---|---|---|---|---|---|---|
| $a_0$ | $a_0$ | $a_0$ | $a_0$ | $a_0$ | $a_0$ | $a_0$ |
| $a_1$ | $a_3$ | $a_3$ | $a_3$ | $a_3$ | $a_3$ | $a_3$ |
| $a_2$ | $a_0$ | $a_0$ | $a_0$ | $a_0$ | $a_0$ | $a_0$ |
| $a_3$ | $a_3$ | $a_3$ | $a_3$ | $a_3$ | $a_3$ | $a_3$ |
| $a_4$ | $a_0$ | $a_0$ | $a_0$ | $a_0$ | $a_0$ | $a_0$ |
| $a_5$ | $a_3$ | $a_3$ | $a_3$ | $a_3$ | $a_3$ | $a_3$ |

$Z_3G$ be the groupoid ring of G over $Z_3$. $W = \{0, 2a_0 + a_3, 2a_3 + a_0\}$. Clearly W satisfies the Lie identities and W is trivially a S-non-associative ring. Hence W is a S-Lie ring but not a Lie ring.

***Example 4.2.2:*** Let $Z_4 = \{0, 1, 2, 3\}$ be the ring of characteristic 4 and G be the groupoid given as in example 4.2.1. Clearly the groupoid ring is a S-Lie ring and not a Lie ring. For take $X = \{0, 2a_0 + 2a_3, a_0 + 3a_3, 3a_0 + a_3\}$.

Clearly $Z_2G$ is a S-Lie ring. But $Z_2G$ is never a Lie ring. In consequence of this we have the following theorems, we in these examples have taken $G \in Z^{***}(n)$.

**THEOREM 4.2.44:** *Let F be any field or a commutative ring. FG the groupoid ring where $G \in Z***(n)$. No groupoid ring in this class of groupoid rings is a Lie algebra.*

*Proof*: F is a field or a commutative ring. Let $G \in Z^{***}(n)$ for some n. FG the groupoid ring of G over F. If FG is to a Lie algebra then it should satisfy the following two identities.

(1) $x^2 = 0$ for all $x \in FG$.
(2) $(xy)z + (yz)x + (zx)y = 0$.

Now let $g = Z_n(t, u)$; atleast one of t or $u \neq 0$. Now $x^2 = 0 \Rightarrow tx + ux = 0$ so that $(t + u)x = 0$ for all $x \in G$ is possible only if $t + u = n$.

Now $(xy)z + (yz)x + (zx)y = (tx + uy)z + (ty + uz)x + (tz + ux)y = t^2x + uty + uz + t^2y + utz + ux + t^2z + tux + uy = 0$.
$t^2 (x + y + z) + ut (y + z + x) + u (x + y + z) = 0$ so that $(t^2 + ut + u) (x + y + z) = 0$.

Thus for this equation to be true we must have $t^2 + ut + u = 0$, $t(u + t) + u = 0$ and $u + t = n$ are true if and only if $u = 0$ forcing $t = 0$. Hence no groupoid in the class of groupoid $Z^{***}(n)$ satisfies the Jacobi identity together with $x^2 = 0$. So no groupoid ring using the class of groupoid from $Z^{***}(n)$ is a Lie algebra.

**THEOREM 4.2.5:** *Let R be a commutative ring with unit or a field, $L_n(m)$ be a loop from the new class of loops $L_n$. The loop ring $RL_n(m)$ is never a Lie algebra.*

*Proof*: Follows from the fact that if $RL_n(m)$ is to be Lie algebra then $x^2 = 0$ and $(xy)z + (yz)x + (zx)y = 0$ for all $x_1y_1z \in RL_n(m)$. If $x \in L_n(m) \subseteq RL_n(m)$ we know $x^2$



= 0 is impossible as in $L_n(m)$, $x^2 = e$ that is every element is self inversed or of order two where e is the identity in $L_n(m)$.

Now if $(xy)z + (yz)x + (zx)y = 0$ for some x, y, z $L_n(m)$, then we have
$[my – (m – 1)x ]z + [mz – (m – 1) y]x + [mx – (m - 1)z]y$

$$= mz – (m – 1)[my – (m – 1) x] + mx – (m – 1)[mz – (m – 1)y]$$
$$+ my – (m – 1)[mx – (m – 1)z]$$

$$= mz – m (m – 1) y + (m – 1)^2 x + mx – (m – 1) mz + (m – 1)^2 y$$
$$+ my – (m – 1) mx + (m – 1)^2 z.$$

$$= m (x + y + z) + (m – 1)^2 (x + y + z) – m (m – 1) (x + y + z)$$

$$= [m + (m – 1)^2 – m (m – 1)] (x + y + z)$$

$$= [ m + m^2 – 2m + 1 – m^2 + m] [x + y + z]$$

$$= 1. (x + y + z ) = 0.$$

This may not be true for all x, y, z $\in RL_n(m)$. Thus no loop ring $RL_n(m)$ from the new class of loops $L_n$ will lead to a Lie algebra. Thus from these theorems we see that the loop rings and groupoid rings constructed by us from the new class of loops or groupoids does not become Lie algebras. Thus they form a distinctly interesting class of non associative rings. But now our study will be "Can ever these class of non-associative rings have Smarandache Lie rings". For this we proceed to study some examples.

These examples show that we can have some of the loop rings or groupoid rings where the loop are from the class of new loops and groupoid rings which are from the new class of groupoids can be Smarandache Lie rings when we take appropriate fields $Z_p$; p can be prime or a composite number.

**THEOREM 4.2.6:** *Let $Z_mG$ be a groupoid ring of the groupoid G, $G \in Z^{***}(n)$ and $Z_m$ be the field or ring of integers modulo m. $Z_mG$ can be a S-Lie ring but never a Lie ring.*

*Proof*: We have already shown that $Z_mG$ is never a Lie ring for all $G \in Z^{***}(n)$. To prove ZnG is a Lie ring it is sufficient if we illustrate by example.

Now we proceed on to show that in the groupoid ring $Z_mG$ where $G \in Z^{***}(n)$ the set $U = \{\alpha = \Sigma\alpha_ig_i / \Sigma\alpha_i = 0\}$; is a ring which is non-associative. Clearly no element in U can be from G. Since m is a composite number or a prime different from 2 we can have in some cases U to be a S-Lie ring, if elements of G satisfy the two identities $x^2 = 0$ and $(xy) z + (yz) x + (zx) y = 0$ for x, y, z $\in Z_mG$.

Thus we have classes of non-associative rings which are never Lie rings but can be S-Lie rings.



**PROBLEMS:**

1. Is $Z_{12}G$ a S-Lie ring where $G = Z_6(3,0)$?
2. Can $Z_{12}L_7(3)$ be a S-Lie ring?
3. Let $Z_{15}L_9(5)$ be a loop ring. $U = \{\alpha = \sum \alpha_i g_i \mid \sum \alpha_i = 0 \text{ and } \alpha \in Z_{15}L_9^{(5)}\}$. Is U a Smarandache NA ring? Does U satisfy Lie identity?
4. Find whether $Z_2G$ where $G = Z_{10}(5,4)$ be a S-Lie ring.
5. Show $Z_2L_n(m)$ for $n = 15$, $m = 14$ is not a Lie ring.
6. Give an example of a S-Lie ring using $Z_6$ and a groupoid from $Z^{***}(8)$.
7. Can $Z_2G$ for any $G \in Z^{***}(n)$ be a S-Lie ring? Justify your answer.
8. Show $Z_4G$ for $G \in Z^{***}(n)$ is a S-Lie ring.
9. Prove $Z_4L_{19}(7)$ is never a Lie ring. Is it a S-Lie ring?
10. Show $Z_5L_{19}(8)$ is not a S-Lie ring.

## 4.3 Substructures in Smarandache Lie algebras

In this section we introduce the concept of substructures like subalgebras, ideals and quotient algebras in Lie algebras and also define their Smarandache analogue. Here we introduce several new concepts in Lie algebras and Smarandache Lie algebras which are not discussed or defined till date. Some of the concepts introduced in this section are, Chinese Lie algebras and its Smarandache analogue, $SI^*$ and $I*$ Lie algebras, ideally obedient Lie algebras and Smarandache ideally obedient Lie algebras CN-Lie algebras and Smarandache CN-Lie algebras. Finally the concepts in Lie algebras like pseudo ideal, semiprime etc are introduced and their Smarandache analogue also discussed.

**DEFINITION 4.3.1:** *Let L be S-Lie algebra. Suppose L has a subalgebra say X and X is a S-Lie algebra then we call X the Smarandache Lie subalgebra of L (S-Lie subalgebra).*

**THEOREM 4.3.1:** *Let L be a Lie algebra if L has a S-Lie subalgebra then L is itself a S-Lie algebra.*

*Proof:* Straightforward from the very definitions.

**DEFINITION 4.3.2:** *Let R and S be two non-associative rings, which are S-Lie algebras. We say $\phi: R \to S$ is a Smarandache Lie algebra homomorphism (S-Lie algebra homomorphism) if $\phi$ restricted to the subsets W of R and V of S respectively where W and V are SNA-subrings which satisfy the Lie identity and $\phi: W \to V$ is a Lie algebra homomorphism.*

*Note:* It is important to note that the mapping need not even be defined on the whole of R it is sufficient if it is defined on W, W a proper subset which is a Lie algebra in R. The concept of Smarandache Lie isomorphism (S-Lie isomorphism) is defined in a similar way.

**DEFINITION 4.3.3:** *Let L be a Lie algebra. A subspace S of L is said to be a subalgebra if $[S, S] \subset S$ i.e. for $x, y \in S$, $[x, y] \in S$. A subspace S of L is said to be an*



*ideal if $[L, S] \subseteq S$ i.e. $x \in L$ and $y \in S$ then $[x, y] \in S$. An ideal is a subalgebra. We say $X \subset L$ is a Smarandache Lie ideal (S-Lie ideal) of L if (1) $X \subset Y$ is a Smarandache Lie Algebra (S-Lie Algebra) where Y is a proper subspace of L which satisfies the Lie identify for all $x \in X$ and $y \in Y$, we have $[x, y] \in X$. If a Lie algebra L has no Smarandache Lie ideals then we call L a Smarandache simple Lie algebra (S-simple Lie algebra).*

**DEFINITION 4.3.4:** *Let L be a real or complex vector space over a field K. A S-Lie ideal X of L is said to be Smarandache Lie solvable (S-Lie solvable) if*

$$X \supset X' = [X, X]$$
$$\supseteq X'' = [X', X'] \ldots$$
$$\supseteq X^k = [X^{k-1}, X^{k-1}] \supseteq \ldots$$

*are Smarandache Lie ideals and $X^k = 0$ for some positive integer k. Similarly we say a S-Lie algebra L is Smarandache Lie solvable if, X the Smarandache K-vectorial subspace of L is such that*

$$X \supseteq X' = [X, X]$$
$$\supseteq X'' = [X^1, X^1] \ldots$$
$$\supseteq X^k = [X^{k-1}, X^{k-1}] \supseteq \ldots$$

*are Smarandache K-vectorial subspaces and $X^k = 0$ for some positive integer K.*

**DEFINITION 4.3.5:** *Let L be a S-Lie algebra, X a S-Lie ideal of L; we say X is a Smarandache minimal Lie ideal (S-minmal Lie ideal) of L if $Y \subset X$ and Y is a S-Lie ideal of X then $X = Y$ or $Y = \phi$.*

**THEOREM 4.3.2:** *Let L be a S-Lie algebra X a S-minimal ideal of L then $[X, X] = 0$ or if $X' = [X X] \subset X$ then X' is not a S-Lie ideal of L or any $[X^k, X^k]$ is not a S-Lie ideal.*

*Proof:* Straightforward by the very definitions of these concepts.

**DEFINITION 4.3.6**: *Let L be a S-Lie algebra. Let $B \subset L$ be a S-Lie ideal of L, we say B is a Smarandache maximal Lie ideal (S-maximal Lie ideal) of L if $L \supset A \supset B$, A any S-Lie ideal of L then $B = A$ or $A = L$.*

**THEOREM 4.3.3:** *Let L be a Smarandache Lie algebra. If B is a S-maximal Lie ideal of L and if $[B B] \subset C$ then either $C = L$ or C is not a S-Lie ideal of L.*

*Proof:* As in the case of minimality, one can easily prove this result.

**DEFINITION 4.3.7:** *Let L be a S-Lie algebra we say the Smarandache D.C.C (S.D.C.C) conditions on S-Lie ideals is satisfied by L if every descending chains $X_1 \supset X_2 \supset \ldots$ of S-Lie ideals of L is sationary. We say the S-Lie algebra satisfies Smarandache A.C.C (S.A.C.C) condition on S-Lie ideals if every ascending chain $Y_1 \subset Y_2 \subset \ldots$ of S-Lie ideals of L is stationary.*



**THEOREM 4.3.4:** *Let L be a S-finite dimensional Lie algebra. L satisfies S.A.C.C on S- Lie ideals of L.*

*Proof:* Straightforward by using the fact that

1. S-basis is finite.
2. All ideals X are basically subspaces of L i.e. $[X\ X] \subset L$.

**DEFINITION 4.3.8:** *Let L be a vector space over the field K. L be a S-Lie algebra. X be a subspace of L i.e. X is only a subspace and not a Lie subspace of L. We call X a Smarandache pseudo Lie ideal (S-pseudo Lie ideal) of L if for every $y \in Y$ where Y is a subset of L and a K-algebra of L, which satisfies the Lie identity and we have $[xy] \in X$ for all $x \in X$.*

**Note***:* $[X, X] \not\subset X$ i.e. $[x, x'] \notin [X, X]$ for some $x, x' \in X$.

Clearly S-pseudo Lie ideals are in general not S-Lie ideals of L. The study of characterization of Smarandache pseudo Lie ideals is left for the reader as an open problem. Any researcher can define Smarandache pseudo minimal Lie ideals and Smarandache pseudo maximal Lie ideals using the usual notions of minimality and maximality. We say a S-Lie ring L is S-pseudo Lie simple if L has no proper Smarandache pseudo Lie ideals (S-pseudo Lie ideals).

**DEFINITION 4.3.9:** *Let L be a S-Lie ring, X a S-pseudo Lie ideal of L then L/X = {a + X / a $\in$ L} is defined as the Smarandache pseudo quotient ring (S-pseudo quotient ring).*

The structure possessed by L / X is an interesting problem for it depends both on the properties enjoyed by L and that of X.

**DEFINITION 4.3.10:** *Let L be a Lie algebra. X be a proper subset which is a K-algebra and satisfies the Lie identity. We say the Lie algebra L is Smarandache semi prime (S-semi prime) if and only if L has no non-zero S-Lie ideal A with [A A] = 0.*

**THEOREM 4.3.5:** *Let L be a Lie ring. If L is S-semi prime then L has no S-minimal ideals.*

*Proof:* Straightforward from the definitions and earlier results.

**DEFINITION 4.3.11:** *Let L be a Lie ring. X a proper S-Lie subring of L; i.e. X is a proper subset which is S-Lie subring. A Smarandache subsemi Lie ideal (S-subsemi Lie ideal) I related to X i.e. I is a S-Lie ideal related to X. and I is not a S-Lie ideal related to any other S-subring containing X.*

*If L has a S-subsemiideal then we say L is a Smarandache subsemi Lie ideal ring (S-subsemi Lie ideal ring). We call the Lie ring a Smarandache filial Lie ring (S-filial Lie ring) if the relation S-Lie ideal in L is transitive that is a S-subring J is a S-Lie ideal in S-subring I and I a S-Lie ideal in L then J is an S-ideal of L.*

Such study has not been imitated till date.



**DEFINITION 4.3.12:** *Let L be a Lie algebra. We say the Lie algebra is a Chinese Lie algebra, if given elements a, b ∈ L and ideals I and J in L such that [J + I, a] = [I + J, b], there exists c ∈ L such that [I, a] = [I, c] and [J, c] = [J, b]; L will be called Smarandache Chinese Lie ring (S-Chinese Lie ring) if we replace the ideals I and J of L by S-Lie ideals.*

**DEFINITION 9.3.13:** *Let L be a Lie algebra. An ideal I of L is said to be a Lie multiplication ideal J ⊂ I if J = [I, C] for some ideal C. We call the multiplication as Smarandache Lie multiplication ideal (S-Lie multiplication ideal) if for S-Lie ideal I of L we have J ⊂ I (J also a S- Lie ideal) such that J = [I,C] where C is a S-Lie ideal.*

**DEFINITION 4.3.14:** *Let L be a Lie algebra. We call I a Lie obedient ideal of L if we have two ideals X, Y in L, X ≠ Y such that [X ∩ I, Y ∩ I] = [X, Y] ∩ I. We say I is a Smarandache Lie obedient ideal (S-Lie obedient ideal) of L, if I is a S-Lie ideal and X and Y X ≠ Y are also S-Lie ideals satisfying the relation [X ∩ I, Y ∩ I] = [X, Y] ∩ I. If every ideal is Lie obedient we call the Lie ring an ideally obedient Lie ring. If every S-Lie ideal is Smarandache Lie obedient we call L a Smarandache ideally obedient Lie ring (S-ideally obedient Lie ring).*

**DEFINITION 4.3.15:** *Let L be a Lie algebra. L is said to be a I∗ - ring if for every pair of ideals $I_1, I_2 \in \{I_J\}$, where $I_J$ denotes the collection of all ideals of $L_v$ and for every $a \in L \setminus (I_1 \cup I_2)$ we have $[a, I_1] = [a, I_2]$ is an ideal of $\{I_f\}$. Suppose for the Lie algebra L if $\{I_j^s\}$ denotes the collection of S-Lie ideals of L we call L a SI∗ ring if for every $a \in L \setminus (I_1 \cup I_2)$ we have $[a, I_1] = [a, I_2]$ is an ideal of $\{I_j^s\}$.*

**DEFINITION 4.3.16:** *Let L be a Lie algebra we say L is locally unitary if for each x ∈ L there exists an idempotent e ∈ L such that xe = ex = x. We call the Lie algebra L Smarandache locally semi unitary (S-locally semi unitary) if for every x ∈ L there exists a S-idempotent e in L such that x e = e x = x.*

**DEFINITION 4.3.17:** *Let L be a Lie algebra. A subset P of L is called a closed net if*

1) *P is a semi group under '+' and '.'.*
2) *aP ⊂ P for all a ∈ K.*

*(K is the field over which L is defined).* We call P a Smarandache closed net (S-closed net) if P is a S-semi group under '.' and a P ⊂ P for all a∈ K.

**DEFINITION 4.3.18:** *We call the Lie algebra L a CN-ring (closed net ring) if L = ∪$P_i$ where $P_i$ ' s are closed nets in L such that $P_i \cap P_j$ = φ or {1} or {0} if i ≠ j and $P_i \cap P_j = P_i$ if i = j (i.e. each $P_i$ is nontrivial).*

*We call the Lie algebra a weakly CN ring or algebra if L ⊂ ∪$P_j$ with $P_i \cap P_j \neq$ φ or {1} or {0} even if i ≠ j.*

**DEFINITION 4.3.19**: *Let L be a Lie algebra {$SP_i$} denote the collection of all S-closed nets of L. We say L is a Smarandache CN-ring (S-CN-ring) if L = ∪ $SP_i$ and $SP_i \cap$*



$SP_j = \phi$ or $\{0\}$ or $\{1\}$ if $i \neq j$ and $SP_i \cap SP_j = SP_i$ if $i = j$. We say L is a Smarandache weakly CN-ring (S-weakly CN-ring) if $L \subset \cup \{SP_i\}$ and $SP_i \cap SP_j \neq \phi$ or $\{0\}$ or $\{1\}$.

**THEOREM 4.3.6:** *Let L be a S-CN ring then L is a S-weakly CN-ring.*

*Proof:* Follows from the very definition.

The proof of the following theorem is left for the reader.

**THEOREM 4.3.7:** *Every S-weakly CN ring in general is not a S-CN-ring.*

<u>PROBLEMS:</u>

1. Give an example of a Lie algebra, which is not a S-Lie algebra.
2. Find S-ideals I in any Lie algebra and show I is not an ideal.
3. Find a Lie ring, which has no S-pseudo ideals.
4. Give an example of a Lie algebra, which has S-pseudo ideals.
5. Give an example of a S-Chinese Lie ring.
6. Prove all algebras are not S-Chinese Lie algebras by an example.
7. Give an example of a S-weakly CN-ring which is not a S-CN-ring
8. Give an example of a Lie algebra, which has no S-obedient Lie ideals.
9. Give an example of a S-ideally obedient ring.
10. Find examples of semi prime Lie rings and S-semi prime Lie rings.
11. Find S-ideals in a Lie ring L, which are not ideals of L.
12. Illustrate by examples: Lie multiplication ideal and S-Lie multiplication ideal of a Lie algebra L.

## 4.4 Special Properties in Smarandache Lie Algebras

In this section we introduce several properties about elements in Lie algebras and their Smarandache analogous. All the properties introduced in this section are very new and interesting. A deep study of these definitions and notions will certainly lead to more innovative and unraveled properties about Lie algebras and in particular of S-Lie algebras. The main new concepts defined in this section are Smarandache Lie related element, Smarandache Lie links, essential Lie algebra and its Smarandache analogue. Stabilizer and stable Lie algebra and its Smarandache notion Lie super ore conditions, conjugate Lie subalgebras and Smarandache conjugate Lie subalgebras and finally the concept of Smarandache generalized stable Lie algebras.

These notions are defined using mainly the Lie product. To the best of our knowledge so far no one has explored or defined these new concepts in Lie algebras, hence cannot be found in S-Lie algebras.

**DEFINITION 4.4.1:** *Let L be a Lie algebra. We say any two elements $a, b \in L$ are Lie related if there exists a subspace A of L such that $a, b \notin A$ but $a \in [A, b]$ and $b \in [A, a]$. We say $a, b \in L$ are Smarandache Lie related (S-Lie related) if there exists a S-subspace A of L such that $a, b \notin A$ but $a \in [A, b]$ and $b \in [A, a]$. We say a and b are Smarandache Lie link related elements (S-Lie link related elements) by the link A.*



**DEFINITION 4.4.2:** *Let L be a Lie algebra. A subalgebra B of L is said to be an essential algebra of L if the intersection of B with every other subalgebra is zero. (By subalgebras we mean only proper subalgebras).*

*A S-subalgebra A of L is said to be a Smarandache essential algebra (S-essential algebra) if the intersection of A with every other S-subalgebra of L is zero. If every subalgebra of a Lie algebra is essential then we call L a essential Lie ring. If every S-subalgebra of a Lie algebra is S-essential then L is a Smarandache essential Lie algebra (S-essential Lie algebra).*

**DEFINITION 4.4.3:** *Let L be a Lie algebra. If for a pair of subalgebras P and Q of L there exists a subalgebra T of L ($T \neq L$) such that [P, T] and [T, Q] are subrings with [P, T] = [T, Q]; then we say the pair of subalgebras are stabilized subalgebras and T is called the stabilizer subalgebra of P and Q.*

*We say a pair of S- subalgebras P and Q of L are Smarandache stabilized algebras (S-stabilized algebras) if there exists a S-subalgebra T such that [P, T] = [T, Q] are S-subalgebras and T is called the Smarandache stabilized subalgebra (S-stabilized subalgebra) of P and Q.*

It is left for the reader to prove the following theporem by constructing an example:

**THEOREM 4.4.1:** *Let L be a Lie ring. If the subalgebras A, B of L is a stable pair then A, B is a stabilized pair and not conversely. The same result holds good for Smarandache stable pair of S-subalgebras.*

**DEFINITION 4.4.4:** *Let L be a Lie algebras if every pair of subalgebras of L is a stable pair, then we say the Lie algebra is a stable algebra. The Lie algebra is a Smarandache stable algebra (S-stable algebra) if every pair of S-subalgebras of L is a Smarandache stable pair (S-stable pair) then we say the Lie algebra is a Smarandache stable algebra (S-stable algebra).*

**THEOREM 4.4.2:** *Every S-stable Lie algebra is a S-stabilized Lie algebra.*

*Proof*: Follows from the very definition hence left for the reader to prove.

**DEFINITION 4.4.5:** *Let L be a Lie algebra. We say a pair of elements $r, s \in L$ satisfies the Lie super ore condition if there is some $r' \in L$ such that [r, s] = [s, r']. If every pair of elements satisfies the Lie super one condition then we say L satisfies the Lie super ore condition.*

**DEFINITION 4.4.6:** *Let L be a Lie algebra such that L is a S-Lie algebra i.e. $X \subset L$ is such that X is a S-K-vectorial subspace satisfying the Lie identity. We say a pair of elements $x, y \in X$ satisfies Smarandache super ore condition (S-super ore condition) if there is some $z \in L$ such That [x, y] = [y, z]. If every pair of elements in X satisfies S-super ore condition then we say L satisfies the Smarandache super ore condition (S-super ore condition).*

**THEOREM 4.4.3:** *Let L be a Lie algebra if L satisfies the S-super ore condition then L need not in general be a Lie algebra satisfying super ore condition.*



*Proof*: The reader is requested to construct an example to prove the theorem which is direct.

**THEOREM 4.4.4:** *If L is a Lie algebra which satisfies the super ore condition and if L is a S-Lie algebra then L satisfies S-super ore condition.*

*Proof*: Straightforward by the very definition.

**DEFINITION 4.4.7:** *Let L be a Lie ring. Two subalgebras A and B of same dimension are conjugate if there exists some $x \in L$ such that $Ax = xB$ ($xA = By$).*

*We call L a conjugate Lie algebra if every distinct pair of subalgebras of same dimension are conjugate.*

Now we define it to the case of Smarandache conjugates.

**DEFINITION 4.4.8:** *Let L be a Lie algebra. Two S-subalgebras A and B of same Smarandache dimension are said to be Smarandache conjugate (S-conjugate) if there exists $x \in L$ such that $Ax = xB$ (or $xA = Bx$). If every pair of S-subalgebras of same Smarandache dimension are S-conjugate then we call L a Smarandache conjugate Lie ring (algebra) (S-co-Lie-ring).*

*We call L a Smarandache weak Lie conjugate ring (S-weak co-Lie ring) if there exists at least one pair of distinct S-subalgebras with same Smarandache dimension which are conjugate to each other.*

**THEOREM 4.4.5:** *Every S-co-Lie-ring is a S-weak co-ring.*

*Proof*: Direct by the very definitions.

**DEFINITION 4.4.9:** *Let L be a Lie algebra and P be a additive subgroup of L. P is called the n capacitor group of L if $[xP] \subset P$ for every $x \in L$. If in case of P an S-semigroup and L a S-Lie algebra having $X \subset L$ where X is a K-algebra of L then we call P a Smarandache n-capacitor group (S-n-capacitor group) of L if $[x, P] \subset P$ for every $x \in X$.*

**THEOREM 4.4.6:** *Let L be a Lie algebra. Suppose I is an ideal of L then I is a n-capacitor group of L.*

*Proof*: Straightforward.

**THEOREM 4.4.7:** *Let L be a Lie algebra if I is a S-ideal of L then I is a S-n-capacitor group of L.*

*Proof*: Direct by the very definitions.

**DEFINITION 4.4.10:** *Let L be a Lie algebra. A subset I of L which is closed under '+' is called a n-ring ideal of L if $[r, i] \subset I$ for all $i \in I$ and $r \in L$. We say the subset I of L which is closed under '+' and has a subgroup under '+' i.e. I is a S-semigroup to*



*be a Smarandache ring n-ideal (S-ring n-ideal) of L if L is a S-Lie algebra and r ∈ X ⊂ L where X is a proper subset which is a K-algebra.*

**DEFINITION 4.4.11:** *Let L be a Lie ring or a Lie algebra. We define P (L) = {x ∈ L / there exist s, t in L such that s(xt) = 1 or (sx)t = 1 }. L is called the generalized stable Lie algebra provided that [a, P(L) ] + [b, P(L)] = P(L) with a, b ∈ L implies a + by ∈ P(L) for some y ∈ L.*

*If L is a S-Lie algebra with X a non-empty subset of L with X a K-algebra and SP(L) = {x ∈ X / there exists s, t in L such that s(xt) = 1 or (sx) t = 1}. L is called the Smarandache generalized stable Lie algebra (S-generalized stable Lie algebra) provided [a, SP(L)] + [b , SP(L)] = SP(L) with a, b ∈ X ; a + by ∈ SP(L) for some y ∈ X.*

**DEFINITION 4.4.12:** *Let L be a Lie algebra. L is said to satisfy the n-stable range conditions for whenever $a_1, a_2, \ldots a_{n+1} \in L$ with $[a_1, L] + [a_2, L] + \ldots + [a_n, L] = L$ there exists $b_1, b_2, \ldots, b_n \in L_n$ such that $[a_1 + a_{n+1}b_1, L] + [a_2 + a_{n+1}b_2, L] + \ldots + [a_n + a_{n+1}, b_n, L] = L$.*

*If L is a S-Lie algebra and X ⊂ L is a proper subset which is a K-vectorial algebra and satisfies the Lie identity and $b_1, b_2, \ldots, b_n \in X$ and $a_1, \ldots, a_{n+1}, \in L$ with $[a_1, X] + \ldots + [a_{n+1}, X] = X$ and $[a_1 + a_{n+1} b_1, X] + \ldots + [a_n + a_{n+1} b_n, X] = X$. Then L is said to satisfy Smarandache n-stable range conditions (S-n-stable range conditions).*

**PROBLEMS:**

1. Give an example of a Lie algebra which satisfies the Lie super ore condition.
2. Give an example of a stable Lie algebra and S-stable S-Lie algebra.
3. Construct a Lie algebra which are conjugate subalgebras but has no S-conjugate subalgebras.
4. Give an example of a Lie algebra which is a n-ring ideal and a S-n-ring ideal.
5. Does there exist a generalized stable Lie algebra? a S-generalized stable Lie algebra?

## 4.5 Some New Notions on Smarandache Lie Algebras

In this section we define Smarandache mixed direct product of Lie algebras. Only this concept paves way to define different types of S-Lie algebras and study them. Here notions of Smarandache simple, Smarandache semi-simple, Smarandache strongly simple, Smarandache strongly semi-simple are introduced. The two types of S-mixed direct product of Lie algebras helps us to see simultaneously certain special structures like S-simple and S-semisimple to be present in them.

**DEFINITION 4.5.1**: *Let $SL_A = L_1 \times \ldots \times L_n$ are Lie algebras defined over the same field K where at least one of the Lie algebras $L_i$ is the linear algebra over K. Then $SL_A$ is a direct product of Lie algebras where the Lie conditions are defined component wise for the n-tuples. As they are defined on the same field compactability exist. For example $[(x_1, \ldots, x_n), (y_1, \ldots, y_n)] = ([x_1 y_1], [x_2 y_2], \ldots,*



$[x_n y_n])$ where $x_i$, $y_i \in L_i$ for $i = 1, 2, \ldots, n$. $SL_A$ is called the Smarandache mixed direct product (S-mixed direct product) of Lie algebras of type A.

**THEOREM 4.5.1**: *Let $SL_A$ be the S-mixed direct product Lie algebra of type A. Then $SL_A$ is a S-Lie algebra.*

*Proof*: Follows from the fact that $SL_A$ contains a Lie algebra $L_i$ which is K-algebra.

**DEFINITION 4.5.2**: *Let L be a S-Lie algebra if L has no non-zero S-solvable Lie ideal then we say L is Smarandache semi-simple (S-semi simple). If L has no S-Lie ideals other than 0 and L we call L, Smarandache simple (S-simple).*

**DEFINITION 4.5.3**: *Let L be a Lie algebra we call L a Smarandache Lie simple (S-Lie simple) if L has a S-subalgebra which is S-simple (i.e. if X is a S-subalgebra of L and X has no S-ideals other than 0 and X).*

**DEFINITION 4.5.4**: *Let L be a Lie algebra, we call L a Smarandache Lie semi-simple (S-Lie semi-simple) if L has a S-subalgebra X such that X is S-semi-simple i.e. X is a S-subalgebra of L and X has no non-zero solvable S-Lie ideal.*

**DEFINITION 4.5.5**: *Let L be a Lie algebra. If every proper S-subalgebra X of L is S-Lie simple we say the Lie algebra is Smarandache strongly Lie simple (S-strongly Lie simple).*

**THEOREM 4.5.2**: *If L is S-strongly Lie simple then L is S-Lie simple.*

*Proof*: Direct by the very definitions of these concepts.

**DEFINITION 4.5.6**: *Let L be a Lie algebra if every S-subalgebra X of L is S-Lie semi-simple then we call L Smarandache strongly Lie semi-simple (S-strongly Lie semi-simple).*

Consequent of these definitions we have the following theorem which is left as an exercise for the reader to prove.

**THEOREM 4.5.3**: *Let L be a Lie algebra which is S-strongly Lie semi-simple then L is S-Lie semi-simple.*

Now we proceed on to define S-mixed direct product of type B.

**DEFINITION 4.5.7**: *Let $SL_B = L_1 \times \ldots \times L_n$ be the direct product of Lie algebras $L_1, \ldots, L_n$ defined over the same field K, in which at least one of the $L_i$ is a linear algebra over K and is a S-semi simple Lie Algebra and one of the $L_j$'s ( $i \neq j$) is a linear algebra over K and is S-simple Lie algebra, then we call $SL_B$ the Smarandache mixed direct product of Lie algebras of type B.*

**THEOREM 4.5.4**: *Let $SL_B$ be the S-mixed direct product of Lie algebras of type B then $SL_B$ is*
        *a. S-Lie semi-simple.*
        *b. S-Lie simple.*



*Proof*: Follows directly by the very definitions.

Thus we see that in case of Smarandache structures only we can have simultaneously the Lie algebra to be S-Lie simple and S-Lie semi-simple. Thus this definition of S-mixed direct product of type B helps us to find algebras which is both S-Lie simple and S-Lie semi-simple.

**DEFINITION 4.5.8**: *Let L be a Lie algebra. We call L a Smarandache sympletic Lie algebra (S-sympletic Lie algebra) if L has a S-subalgebra X, which is a sympletic Lie algebra.*

*Similarly one can define Smarandache triangular (S-triangular) and Smarandache orthogonal (S-orthogonal) Lie algebras.*

**DEFINITION 4.5.9**: *Let L be a Lie algebra. X be a S-subalgebra of L. The Smarandache normalizer (S-normalizer) $SN(X) = \{x \in L \;/\; [x, X] \subseteq X\}$ that is $[x, y] \in X$ for every $y \in X$. $SN(X)$ is a subalgebra. We call $SN(X)$ the Smarandache normalizer only when $SN(X)$ is a S-subalgebra.*

**DEFINITION 4.5.10**: *A S-subalgebra X of a Lie algebra L is called a Smarandache Cartan subalgebra (S-Cartan subalgebra) if*

    a. *X is nilpotent.*
    b. *X is its own S-normalizer in L.*

Thus the Smarandache mixed direct products of Lie algebras of type A and B leads us to construct more and more examples of S-Lie algebras, S-Lie simple algebras and S-Lie semi-simple algebras. Further S-strong simple and semi-simple algebras have been also introduced. So given any Lie algebra which is not a S-Lie algebra can be made into a S-Lie algebra by S-mixed direct product of Lie algebras.

**PROBLEMS**:

1. Give an example of a S- Lie simple algebra.
2. Illustrate by an example the notion of S-strong Lie semi-simple algebra.
3. Give an example of a Lie algebra which is both S-Lie simple and S-Lie semi-simple.
4. Show every Cartan subalgebra need not be a S-Cartan subalgebra.
5. Give an example of S-mixed direct product which is a S-strong semi-simple algebra.
6. Can type B S-mixed direct product of Lie algebras be S-strong simple Lie algebras? Justify your claim.
7. Using S-mixed direct product of Lie algebras construct a S-Cartan Lie algebra.



**Chapter 5**

# JORDAN ALGEBRAS AND SMARANDACHE JORDAN ALGEBRAS

This chapter has two sections. The first section is devoted to the recollection and introduction of basic concepts about Jordan algebras. As the class of Jordan algebras is a large one and there are good books which give a complete and an excellent exposition about Jordan algebras [2, 37, 51] we don't intend going deep into the study. As the main aim of this book is only to introduce Smarandache non-associative algebras we don't give a complete analysis of Jordan algebras. Further we are mainly interested in the algebraic aspect of it and not view it as a vector space or its associated transformations. We describe and study chiefly Jordan algebras as non-associative algebras which satisfies the Jordan identity.

Section 2 of this chapter introduces the concept of Smarandache Jordan algebras. We introduce several properties enjoyed by associative rings some of them have been modified and incorporated as Smarandache notions in S-Jordan algebras. To the best of our knowledge such study of introduction of these concepts in Jordan algebras is totally absent. Further all books on non-associative algebras which deals with Jordan algebras studies Jordan algebras in a very different angle. Here our study of Smarandache Jordan algebra in this book is varied and makes an interesting one as the notions/ concepts enjoyed by associative algebras cannot be totally deprived from our study when we assume all the algebraic structure remains the same, except the associativity.

So keeping this in mind we have introduced several new notions like Smarandache Jordan strong ideal property, Smarandache power-joined algebras, Smarandache magnifying and Smarandache shrinking elements in Jordan algebras and other Smarandache concepts.

We have left for the reader to develop each of these concepts and obtain many innovative results about Smarandache strong Jordan algebras, Smarandache commutative Jordan algebra, Smarandache weakly commutative Jordan algebras, Smarandache pseudo commutative Jordan algebras, Smarandache quasi commutative Jordan algebra, Smarandache reduced Jordan algebras, Smarandache periodic Jordan algebras, Smarandache E-Jordan algebras, and Smarandache pre-Boolean algebras. The sole aim of this text is to give innovative and new concepts and notions which can be developed into a very good theory.

## 5.1 Basic properties of Jordan Algebras

In this section we just recall the definition of Jordan algebras and derive some of its basic properties. Several properties about Jordan algebras and its substructures like radical and ideal are given. Examples are also given. For more about Jordan algebras please refer [2, 37, 51].



**DEFINITION [51]**: *Let A be a finite dimensional commutative algebra with identity. A is not assumed to be associative. A is called a Jordan algebra if we have $x^2(xy) = x(x^2y)$ for all $x, y \in A$. An equivalent definition is A is a Jordan algebra if and only if we have for all invertible $x \in A$ and for all $y \in A$, $x^{-1}(xy) = x(x^{-1}y)$.*

**DEFINITION 5.1.1**: *A quadratic Jordan algebra is a triple $J = (V, P, e)$ where V is a finite dimensional vector space over K, P a quadratic map of V into the space end(V) of its endomorphisms and e a non-zero element of V such that the axioms to be stated below hold. Putting $P(x, y) = P(x + y) - P(x) - P(y)$ we obtain a symmetric bilinear map P, of $V \times V$ into End V. The axioms for a quadratic Jordan algebras are as follows:*

1. *$P(e) = e$, $P(x, e)y = P(x, y)e$.*
2. *$P(P(x)y) = P(x)P(y)P(x)$.*
3. *$P(x)P(y, z)x = P(P(x)y, x)z$.*

*J is said to be defined over $k \subset K$ if there exists a k-structure on the vector space V such that P is defined over k and that $e \in V(k)$. $J = (V, P, e)$ is a quadratic Jordan algebra.*

For more properties and results please refer [51] . We do not approach or define Jordan algebras using Zaraski topology. We purely deal it as an algebraic structure or to be more precise as an algebra i.e. like an algebra satisfying the identity $a^2(au) = a(a^2u)$.

**DEFINITION 5.1.2:** *Let A and B be any two Jordan algebras where we assume they are finite dimensional over a field K of characteristic not equal to two. A linear mapping $\sigma : A \to B$ is called a Jordan homomorphism if*

i. *$\sigma(a^2) = (\sigma(a))^2$ for every $a \in A$.*
ii. *$\sigma(a.b.a) = \sigma(a)\sigma(b)\sigma(a)$ for every $a, b \in A$.*

*If B does not contain a zero divisor then $\sigma(a . b) = \sigma(a) . \sigma(b)$ or $\sigma(a . b) = \sigma(b) . \sigma(a)$.*

*We can also derive or define special Jordan algebras using an associative algebra A. Let A be an associative algebra, Define a new multiplication '.' in A by $a . b = (ab + ba)/2$. We then have a Jordan algebra $A^+$. A subalgebra of the Jordan algebra $A^+$ is called the special Jordan algebra. We just recall the notion of free special Jordan algebra. Let $K[x_1, \ldots, x_n]$ be the non-commutative free-ring in the indeterminates $x_1, \ldots, x_n$ (that is $K[x_1, \ldots, x_n]$ is the associative algebra over K that has as its K-basis the free semigroup with identity element 1 over the free generators $x_1, x_2, \ldots, x_n$. The subalgebra $K[x_1, \ldots, x_n]^+$ generated by 1 and the $x_i$'s is called the free special Jordan algebra of n-generators and is denoted by $J_o^{(n)}$. A Jordan algebra A is special if and only if there is an isomorphism from A onto $B_o^+$, where B is some associative algebra. A Jordan algebra that is not special is called exceptional.*

*A Jordan algebra A has a unique solvable ideal N, which contains all nilpotent ideals of A and is called the radical of A. If $N = 0$, A is called semisimple. The*



*quotient A/N is always semisimple. A semisimple Jordan algebra A contains the unit element and can be decomposed into a direct sum, $A = A_1 \oplus \ldots \oplus A_r$ of minimal ideals $A_i$, each $A_i$ is a simple algebra. In particular if K is of characteristic zero there is a semisimple subalgebra S of A such that $A = S \oplus N$. Let e be an idempotent element of A and let $\lambda \in K$. Put $A_e(\lambda) = \{x \mid x \in A, e \cdot x = \lambda x\}$. Then we have $A = A_e(1) \oplus A_e(1/2) \oplus A_e(0)$. This decomposition is called the Pierce decomposition of A relative to e. Suppose that the unity element 1 is expressed as a sum of the mutually orthogonal idempotents of e.*

*A representation S of a Jordan algebra A on a K-linear space M is a K-linear mapping $a \to S_a$ from A into the associative algebra E(M) of all K-endomorphisms of M such that*

1. *$[S_a, S_{b.c}] + [S_b, S_{c.a}] + [S_c, S_{a.b}] = 0$ and*
2. *$S_a S_b S_c + S_c S_b S_a + S_a S_c S_b = S_a S_{b.c} + S_b S_{ac} + S_c S_{a.b}$ for all a, b, c in A.*

*A K-linear space M is called a Jordan module of A if there are given bilinear mappings $M \times A \to M$ (denoted by $(x, a) \to x.a$), $A \times M \to M$ (denoted by $(a, x) \to ax$) such that for every $x \in M$ and for every a, b, c $\in$ A*

i. *$x \cdot a = a \cdot x$.*
ii. *$(x \cdot a) \cdot (b \cdot c) + (x \cdot b)(a \cdot c) + (x \cdot c)(a \cdot b) = (x \cdot (b \cdot c)) \cdot a$*
   *$+ (x \cdot (a \cdot c)) \cdot b + (x \cdot (a \cdot b)) \cdot c$ and*
iii. *$x \cdot a \cdot b \cdot c + x \cdot c \cdot b \cdot a + x \cdot a \cdot c \cdot b = (x \cdot c)(a \cdot b) + (x \cdot a)(b \cdot c)$*
   *$+ (x \cdot a)(a \cdot c)$.*

*As usual there is a natural bijection between the representations of A and the Jordan modules of A. A special representation of a Jordan algebra A is a homomorphism $A \to E^+$, where E is an associative algebra. Among the special representations of A there exists a unique universal one in the following sense. There exists a special representation $S : A \to U^+$ with the following property. For every special representation $\sigma : A \to E^+$ there exists a unique homomorphism $\eta : U^+ \to E^+$ such that $\sigma = \eta S$. The pair (U, S) is uniquely determined. Further more, if A is n-dimensional over K, U is of dimension $\binom{2n+1}{m}$ over K. The pair (U, S) is called the special universal enveloping algebra of A. A is special if and only if $S : A \to U^+$ is injective.*

For more about Jordan algebras please refer [2, 37, 51].

Now we proceed on to study which of the loop rings using the new class of loops are Jordan algebras and which of the groupoid rings are Jordan algebras. We first recall the definition of Jordan loop [66]

**DEFINITION [66]**: *Let (L, ., e) be a loop. We say L is a Jordan loop if*

i) *$ab = ba$.*
ii) *$a^2(ab) = a(a^2b)$ for all a, b $\in$ L.*



***Example 5.1.1***: Let L = {e, $g_1$, $g_2$, …, $g_7$} be a Jordan loop given by the following table:

| . | e | $g_1$ | $g_2$ | $g_3$ | $g_4$ | $g_5$ | $g_6$ | $g_7$ |
|---|---|---|---|---|---|---|---|---|
| e | e | $g_1$ | $g_2$ | $g_3$ | $g_4$ | $g_5$ | $g_6$ | $g_7$ |
| $g_1$ | $g_1$ | e | $g_5$ | $g_2$ | $g_6$ | $g_3$ | $g_7$ | $g_4$ |
| $g_2$ | $g_2$ | $g_5$ | e | $g_6$ | $g_3$ | $g_7$ | $g_4$ | $g_1$ |
| $g_3$ | $g_3$ | $g_2$ | $g_6$ | e | $g_7$ | $g_4$ | $g_1$ | $g_5$ |
| $g_4$ | $g_4$ | $g_6$ | $g_3$ | $g_7$ | e | $g_1$ | $g_5$ | $g_2$ |
| $g_5$ | $g_5$ | $g_3$ | $g_7$ | $g_4$ | $g_1$ | e | $g_2$ | $g_6$ |
| $g_6$ | $g_6$ | $g_7$ | $g_4$ | $g_1$ | $g_5$ | $g_2$ | e | $g_3$ |
| $g_7$ | $g_7$ | $g_4$ | $g_1$ | $g_5$ | $g_2$ | $g_6$ | $g_3$ | e |

Clearly $g_i g_j = g_j g_i$ and $a^2(ab) = a(a^2 b)$ for all a, b ∈ L = {e, $g_1$, …, $g_7$}. We have a class of loops of even order which forms a Jordan loop.

**DEFINITION [66]**: *Let $J_p$ = {e, $g_1$, …, $g_p$} where p is a prime greater than 3. Define '*' on $J_p$ by*

   a. $g_i * g_i = e$.
   b. $e * g_i = g_i * e = g_i$ (e identity element of $J_p$).
   c. $g_i * g_j = g_t$.

*where*

$$t = \left[\frac{(p+1)j}{2} - \frac{(p-1)i}{2}\right] (mod\ p)$$

*for all $g_i$, $g_j$ ∈ $J_p$, $g_i \neq g_j$. Clearly $J_p$ is a loop which is easily verified to be a Jordan loop and order of $J_p$ is even. $C(J_p) = \{J_p / p$ is an odd prime greater than 3\} forms a class of Jordan loops of even order.*

**THEOREM [66]**: *Let $Z_2 = \{0, 1\}$ be the prime field of characteristic two and $J_p$ be a loop. The Jordan loop ring $Z_2 J_p$ is a Jordan ring.*

*Proof*: Given $J_p$ = {e, $g_1$, …, $g_p$}, p > 3, p a prime. $Z_2 = \{0, 1\}$ be the field of characteristic two. $Z_2 J_p$ is a loop-ring. Every element α ∈ $Z_2 J_p$ is of the form $\alpha^2 = 0$ or 1 according as m is even or odd (where $\alpha = \sum_{i=1}^{m} \alpha_i g_i$ i.e. m is the number of terms which are non-zero in α). Suppose $\alpha^2 = 0$, then for any β we have $\alpha^2(\alpha\beta) = 0$ and $\alpha(\alpha^2\beta) = 0$. Thus $\alpha^2(\alpha\beta) = \alpha(\alpha^2\beta)$. Suppose $\alpha^2 = 1$, then for any β ∈ $Z_2 L$ we have $\alpha^2(\alpha\beta) = \alpha\beta$ and $\alpha(\alpha^2\beta) = \alpha\beta$. Thus once again $\alpha^2(\alpha\beta) = \alpha(\alpha^2\beta)$. Hence $Z_2 J_p$ is a Jordan loop ring. Thus we have got a class of Jordan rings of varying order given by $Z_2 J_p$ where $J_p$ ∈ $C(J_p)$.

**THEOREM [66]:** *Let K be a field of characterisitic zero and $J_p$ ∈ $C(J_p)$. The loop ring $KJ_p$ is a Jordan ring.*

*Proof*: Follows from the fact that every $\alpha, \beta \in KJ_p$ we have $\alpha^2(\alpha\beta) = \alpha(\alpha^2\beta)$ as every pair of elements in $J_p$ satisfies the identity $a^2(ab) = a(a^2b)$ for every $a, b \in J_p$.

Now we recall a groupoid in the class of groupoids $Z^{***}(n)$ to be a Jordan groupoid.

**DEFINITION 5.1.3**: *A groupoid $(G, *)$ is said to be Jordan groupoid if*

    a. $ab = ba$.
    b. $a^2(ab) = a(a^2b)$

*for all $a, b \in G$.*

**THEOREM 5.1.1**: *Let $G$ be a commutative groupoid with 1 in which every $g \in G$ is such that $g^2 = 1$. F be any field. Then the groupoid ring FG is a Jordan ring.*

*Proof*: Straightforward by the definitions and simple calculations.

Thus we see groupoid rings and loop rings can also be Jordan rings for suitable loops and groupoids.

**DEFINITION [37]**: *An element $a$ of a Jordan algebra is called quasi invertible if there exists $b$ in the algebra such that $a + b - a \cdot b = 0$ and $a + b - 2a \cdot b - a^2 + a^2 \cdot b = 0$. Then $b$ is called the quasi inverse of $a$. For algebras with 1 this is equivalent to $1 - a$ is invertible with inverse $1 - b$.*

*Further no idempotent is quasi invertible and every nilpotent element is quasi invertible. Further it has been proved that for an algebraic Jordan algebra the following two conditions are equivalent:*

    a. *Every element is nilpotent.*
    b. *Every element is quasi invertible.*

### 5.2 Smarandache Jordan Algebras and its basic properties

Throughout this section we assume the Jordan algebra to be an algebraic Jordan algebra A, i.e, a Jordan algebra is one which satisfies the identity $x(x^2y) = x^2(xy)$, $x, y \in A$. If A is commutative then $ab = ba$, otherwise we say A is a non-commutative Jordan algebra. A is assumed to contain the multiplicative identity 1, i.e. A is assumed to be a finite dimensional algebra with 1. We are only interested in analyzing its algebraic aspects and we do not see its properties via any other mathematical structure. All the more we do not even deal it as a matrix structure. For the main reason is that Smarandache matrix theory and Smarandache vector spaces and Smarandache linear algebras are yet to be completely developed into solid theories.

**DEFINITION 5.2.1**: *Let A be any algebra. We say A is a Smarandache Jordan algebra (S-Jordan algebra) if A has a proper subset X which is a subalgebra with respect to the operation of A and X satisfies the identity $x(x^2y) = x^2(xy)$ for all $x, y \in X$.*

**THEOREM 5.2.1**: *Let A be a Jordan algebra then A is a S-Jordan algebra.*

*Proof*: Direct by the very definition.

**THEOREM 5.2.2**: *Let A be a S-Jordan algebra then A need not be a Jordan algebra.*

*Proof*: By an example.

Let $Z_2L$ be the loop algebra of a loop L over $Z_2$. Here L is a loop having a subloop H such that H is commutative and every element in H is such that $g^2 = 1$ for every $g \in$ H. Clearly the loop algebra $Z_2G$ has $Z_2H$ as a subalgebra and $Z_2H$ satisfies the Jordan identity. $x(x^2y) = x^2(xy)$ for all x, y $\in Z_2H$.

Thus $Z_2G$ is not a Jordan algebra but $Z_2H$ is a Jordan algebra, hence $Z_2G$ is a S-Jordan algebra.

**DEFINITION 5.2.2**: *Let A be a algebra. If every subalgebra of A satisfies the Jordan identity then we call A the Smarandache strong Jordan algebra (S-strong Jordan algebra).*

**THEOREM 5.2.3**: *If A is a S-strong Jordan algebra then A is a S-Jordan algebra.*

*Proof*: Follows by the very definitions directly.

The reader is expected to prove the following theorem:

**THEOREM 5.2.4**: *If A is a S-Jordan algebra then A need not in general be a S-strong Jordan algebra.*

We derive and define several algebraic properties and their Smarandache analogue in case of Jordan algebras. We see in case of Jordan algebras the concept of S-idempotents are defined analogous to S-idempotents in usual associative algebras. We define several nice and new notions and their Smarandache analogue.

**DEFINITION 5.2.3**: *Let A be an algebra. We say a proper subset B of A to be a Smarandache subJordan algebra (S-subJordan algebra) if B $\subset$ X, X $\subset$ A, X is a S-Jordan algebra.*

**DEFINITION 5.2.4**: *Let A be an algebra which is S-Jordan algebra. We say A is a Smarandache commutative Jordan algebra (S-commutative Jordan algebra) if every subalgebra is commutative, We call A to be a Smarandache weakly commutative Jordan algebra (S-weakly commutative Jordan algebra) if at least one of the subalgebra which is a S-Jordan algebra is commutative.*

**THEOREM 5.2.5**: *Every S commutative Jordan algebra is S-weakly commutative Jordan algebra.*

*Proof*: Follows directly by the very definitions.

**THEOREM 5.2.6**: *Let A be an S-commutative Jordan algebra. A need not be a commutative algebra.*

*Proof*: The reader is requested to construct an example to prove the theorem.

We define Smarandache ideals in algebras only algebraically.

**DEFINITION 5.2.5**: *Let A be an algebra. $X \subset A$ be a subset which is a S-Jordan algebra i.e. X satisfies the Jordan identity. A proper subset $J \subset X$ is called the Smarandache Jordan ideal (S-Jordan ideal) if J is a subalgebra of X and for every a $\in$ A we have aj and ja $\in$ J. We can define Smarandache Jordan left ideal or Smarandache Jordan right ideal. So Smarandache Jordan ideal is one which is simultaneously both a S-Jordan left and a S-Jordan right ideal.*

**DEFINITION 5.2.6**: *Let A be a S-Jordan algebra if A has no S- Jordan ideals then we say A is Smarandache simple Jordan algebra (S-simple Jordan algebra).*

**DEFINITION 5.2.7**: *Let A be a Smarandache Jordan algebra. A pair of distinct elements x, y $\in$ X different from identity of R which are such that xy = yx is said to be Smarandache pseudo Jordan commutative pair (S-pseudo Jordan commutative pair) of X, where X is the subalgebra of A which satisfies the Jordan identity if for all a $\in$ A we have x(ay) = (ya)x (or y(ax)) or (xa)y = (ya)x (or y(ax)). If every commuting pair of X happens to be S-pseudo Jordan commutative then we say A is a Smarandache pseudo commutative Jordan algebra (S- pseudo commutative Jordan algebra).*

This concept of pseudo commutative has meaning even if we assume A to be a commutative algebra. Now we proceed on to define quasi commutativity in non-commutative Jordan algebras.

**DEFINITION 5.2.8**: *Let A be non-commutative algebra. A is said to be quasi commutative if $ab = b^\gamma a$ for every pair of elements a, b $\in$ A and $\gamma > 1$. We say the S-Jordan algebra A is Smarandache quasi commutative Jordan algebra (S-quasi commutative Jordan algebra) if every pair of elements in $X \subset A$ where X is a subalgebra, is satisfying the Jordan identity is a quasi commutative algebra.*

*Thus we see if A is a S-quasi commutative Jordan algebra then A in general need not be a quasi commutative ring.*

**DEFINITION 5.2.9**: *Let A be an algebra or a Jordan algebra. An element x $\in$ A is said to be semi-nilpotent if $x^n – x$ is a nilpotent element of A. If $x^n – x = 0$ we say x is trivially semi-nilpotent.*

**THEOREM 5.2.7**: *Let A be a Jordan algebra. Then the following conditions are equivalent*

        a. *Every element x is semi-nilpotent.*
        b. *Every element $x^n – x$ is quasi invertible.*

*Proof*: We know if x is semi-nilpotent then $x^n - x$ is nilpotent so $x^n - x$ is quasi invertible [37].

If A is only a Smarandache Jordan algebra then we have the modified version of the above theorem.

**THEOREM 5.2.8**: *Let A be a S-Jordan algebra i.e. $X \subset A$ satisfies the Jordan identity then the following are equivalent:*

   a. *Every element x in X is semi-nilpotent.*
   b. *Every element $x^n - x$ is quasi invertible.*

*Proof*: Left for the reader to prove.

The concept of semi-nilpotent cannot be defined for Lie algebras as $x^2 = 0$ in case of Lie algebras.

**DEFINITION 5.2.10**: *We call the algebra A to be reduced if R has no non-zero nilpotents. We call A to be Smarandache reduced Jordan algebra (S-reduced Jordan algebra) if the subalgebra $X \subset A$; X satisfies the Jordan identity but X has no non-zero nilpotents.*

*Thus we see when A is a S-reduced Jordan algebras A can have non-zero nilpotents only $X \subset A$ mentioned in the definition does not contain non-zero nilpotents.*

**DEFINITION 5.2.11**: *Let A be a algebra or a S-Jordan algebra, we say A is Smarandache periodic (S-periodic) if $x^n = x$ for some integer n > 1, for all $x \in X \subset A$. X is a subalgebra of A which is a Jordan algebra.*

*Note if n is even in the definition and if A is of characteristic two then we call the ring A, a Smarandache E-Jordan ring (S-E-Jordan ring).*

A quadratic Jordan ring is said to be periodic of every element a in it satisfies the condition $a^n = a$ for some integer n > 1 and has obtained several interesting results in this direction.

**DEFINITION 5.2.12**: *Let A be a Jordan algebra or a S-Jordan algebra ($X \subset A$) we say A is said to be a Smarandache pre J- Jordan algebra (S-pre-J Jordan algebras) if $a^n b = a b^n$ for any pair $a, b \in X \subset A$ and n is a positive integer.*

**DEFINITION 5.2.13**: *Let A be an algebra which is a S-Jordan algebra we say A is a Smarandache pre-Boolean Jordan algebra (S-pre-Boolean Jordan algebra) if for $X \subset A$ we have $xy(x + y) = 0$ for $x, y \in X$.*

**DEFINITION 5.2.14**: *Let A be a algebra. We say A is a Smarandache-J-Jordan algebra (S-J-Jordan algebra) if $Y \subset X \subset A$ (where X is a subalgebra satisfying the Jordan identity) if for all $y \in Y$ we have $y^n = y$, y > 1.*

**DEFINITION 5.2.15**: *Let A be an algebra $\phi \neq X \subset A$ be a S-Jordan algebra; we say A satisfies Smarandache Jordan strong ideal property (S-Jordan strong ideal property) if every pair of ideals of X generates X.*

*We call A a Smarandache Jordan weak ideal ring (S-Jordan weak ideal ring) if there exists at least a pair of distinct ideals of X which generate X.*

The reader is requested to get results related to these notions.

**DEFINITIONS 5.2.16**: *Let A and $A_1$ be two S-Jordan algebras. We say a map $\phi : A \to A_1$ is said to be a Smarandache Jordan algebra homomorphism (S-Jordan algebra homomorphism) if $\phi$ from X to $X_1$ i.e. $\phi : X \to X_1$ (where X and $X_1$ are subalgebras of A and $A_1$ respectively satisfying the Jordan identity) is a linear transformation as subspaces. We are not interested whether $\phi$ is defined entirely on A or $\phi$ is a linear transformation from the algebras A to $A_1$.*

**DEFINITION 5.2.17**: *Let A be an algebra. We say A is Smarandache power-joined Jordan algebra (S-power-joined Jordan algebra) if $a^n = b^m$ for some positive integers m and n and a, b $\in$ X where X $\subset$ A is a Jordan algebra.*

*We say A is Smarandache uniformly power-joined Jordan ring (S-uniformly power-joined Jordan ring) if m = n i.e. $a^n = b^n$ for a, b $\in$ X.*

The reader is expected to obtain some interesting results in this direction.

**DEFINITION 5.2.18**: *Let A be an S-Jordan algebra, i.e. X $\subset$ A be a subalgebra which satisfies the Jordan identity. We say an element x $\in$ A to be Smarandache shrinkable (S-shrinkable) element if xA = X (Ax = X). The study of the nature of x $\in$ A as interesting for all elements in A cannot serve this purpose; only few special elements will serve as shrinkable elements. Further every S-Jordan ring need not in general be shrinkable.*

**DEFINITION 5.2.19**: *Let A be an S-Jordan algebra. Let X $\subset$ A be a proper subset of A which is a subalgebra satisfying the Jordan identity. We say an element m $\in$ A, a Smarandache magnifying element (S-magnifying element) if Xm = A (mX = A). It is to be noted that all S-Jordan algebras need not have S-magnifying elements.*

**DEFINITION 5.2.20**: *Let A be a S-Jordan algebra we say A is Smarandache dispotent Jordan algebra (S-dispotent Jordan algebra) if X ($\subset$ A) has exactly two S-idempotents where X is a subalgebra of A satisfying the Jordan identity. If every subalgebra X of A which satisfies Jordan identity has exactly two S-idempotents then we call A a Smarandache strong dispotent Jordan algebra (S-strong dispotent Jordan algebra).*

**DEFINITION 5.2.21**: *Let A be a S- Jordan algebra. We say A is a Smarandache normal Jordan algebra (S-normal Jordan algebra) if aX = Xa for every a $\in$ A where X is a subalgebra of A which satisfies the Jordan identity.*

**DEFINITION 5.2.22**: *Let A be a S-Jordan algebra. We say A is a Smarandache weakly G-Jordan ring (S-weakly G-Jordan ring) if for every additive S-semigroup P of X (X a subalgebra satisfying the Jordan identity) we have xP = Px for every x ∈ X. If xP = Px = P we say A is a S-G-Jordan ring.*

**DEFINITION 5.2.23**: *Let A be a S-Jordan algebra. An element x ∈ X ⊂ A is said to be a Smarandache clean element (S-clean element) if x can be written as a sum of a S-idempotent and a S-unit in A.*

**DEFINITION 5.2.24**: *Let A be a S-Jordan algebra. A subset M of X, X ⊂ A is said to be Smarandache system of local units (S-system of local units) if and only if*

   a. *M consists of commuting S-idempotents.*
   b. *For any x ∈ X and e ∈ M, ex = xe = x.*

**PROBLEMS:**

1. Give examples of S-G-Jordan ring.
2. Find an example of a S-weakly G-Jordan ring which is not a S-G Jordan ring.
3. What is the order of a smallest S-Jordan algebra?
4. Give an example of a S-strong Jordan algebra?
5. Find an example of a S-Jordan algebra which has S-shrinkable elements.
6. Does every S-Jordan algebra have a S-magnifying element? Justify your answer.
7. Illustrate by an example a S-normal Jordan algebra.
8. Give an example of a S-periodic Jordan ring.
9. Is the loop ring $Z_2L_5(2)$ a S-periodic Jordan ring?
10. Illustrate by an example a S-pre-Boolean Jordan ring.
11. Is $Z_2L_7(3)$ a S-E-Jordan ring?
12. Does the ring $Z_8L_5(3)$ have semi-nilpotent elements?

**Chapter 6**

# SUGGESTIONS FOR FUTURE STUDY

This chapter is unique for it gives suggestions for future study. This chapter is divided into five sections. First section is devoted to the exposition of groupoid rings and the probable study of it. As groupoid rings which is the most generalized structure of non-associative rings till date has not found its place either by way of research publication or in text form, so the author feels that the study of groupoid rings will certainly lead to a lot of research. The second section of this chapter gives a glimpse about the non-associative rings called loop rings.

Research on loop rings is very systematic for there are papers about it but yet there is not even a single textbook fully devoted to loop rings whereas there are several books on group rings. In the third section we give the recent concepts about Lie algebras and Jordan algebras are taken care of in section four. In the final section of this chapter several other recently discovered non-associative rings found in literature are given for the reader to make proper Smarandache developments about them.

## 6.1 Smarandache groupoid rings

The study of groupoids is very meagre. Except for the book [70] where the notions of Smarandache groupoids and groupoids are dealt with extensively. So the groupoid rings which are an enormously a large class of non-associative rings, for they can be constructed using the infinite number of naturally built groupoids using $Z_n$, Z, Q and R over any appropriate or needed ring. Except for a few papers on groupoid rings [31, 64], the study in this direction is totally absent.

Now as groupoid rings are the generalization of group rings, loop rings and semigroup rings and Smarandache groupoid rings contain always semigroup rings as subrings. The study of this concept will be beneficial. As Lie identities and Jordan identities have given rise to two class non-associative algebras, viz., Lie algebras and Jordan algebras. The first thing the reader can do is to develop using groupoids the following types of algebras as they are also governed by identities like Moufang, Bol, etc.

1. Moufang algebras using groupoid rings and excavating the probabilities of giving them vector space representations. The immediate consequence would be the study of Smarandache Moufang algebras using groupoid rings.

2. Construction and study of Bol algebras and Smarandache Bol algebras by using groupoid rings.

3. Construction and study of P-algebras and Smarandache P-algebras using P-identities in groupoid rings.

4. Analysis of alternative algebras and Smarandache alternative algebras using groupoid rings.

5. As groupoids find its application in Smarandache automaton we can certainly see that the groupoid rings will find its application in coding theory.

The final thing of importance is that in this book we have shown that groupoid rings, which are not found to be Jordan algebras, are found to be Smarandache Jordan algebras. Thus the author at this juncture wants to make known to the researchers that by studying the Smarandache notions certainly certain subproperties not found in the total structure can be dealt without complications in a very systematic way by defining the Smarandache analogue.

Finally the author leaves it to the reader to develop these identities using, the program first designed and implemented at Clemson University about 1990 by [31]. They have studied groupoids satisfying the identity $(xy)z = y(zx)$ and it has been analysed.

Thus the reader can use this program in two ways in our groupoids and groupoid rings:

I. Characterize those groupoid rings satisfying the identities.

II. Using these groupoid to build and develop the theory of formal languages for which Smarandache groupoids will serve the proper perspective.

## 6.2 Smarandache loop rings

The study of loop rings has been carried out by several researchers [8 to 16, 21, 24-29, 49, 50, 57-65]. But the only disadvantage is that we do not have any text solely on loop rings as we have for group rings or semigroup rings. As loop rings forms the generalization of group rings and it has more properties in them, as they happen to be non-associative we felt loop rings should be given a good exposition.

For study of loop rings or Smarandache loop rings the most needed concept is the loops and Smarandache loops. They have been introduced in [5, 71, 72]. It is still a pity that we do not have book solely devoted to the study of loops. Hence in this book we have devoted a complete chapter on the study of loop rings and Smarandache loop rings.

As Lie algebras or Jordan algebras are algebras which satisfy certain special identities and are vividly studied, we have several books on them whereas we see that the loop algebras which are more generalized class of non-associative algebras does not find its just place in the mathematical world. Alternative algebras has been intensely studied by [24 to 29] but it could be developed into a book on alternative and Smarandache alternative algebras.

So we have the following suggestions for the future study.

1. Which of the alternative loop rings are Lie or Jordan rings?

2. Construct new class of algebras using loop algebras like "Moufang algebras" in par with Lie or Jordan algebras. Likewise 'Bol algebras', 'Bruck algebras' and 'P-algebras' for all these are also algebras satisfying only a special identity.

It is unknown why such a study has not been carried out in these cases. The author sees that several of the loop algebras which failed to satisfy Jordan identity happen to be Smarandache Jordan algebra. So we predict several of the loop rings which does not satisfy special identities can by all means be Smarandache.

3. Finally the loop rings can be used in the study of algebraic coding theory.

Several of the properties studied for group rings can be studied for Smarandache loop rings and loop rings. Using the Albert program to study non-associative algebras satisfying special identities we can find whether the loop rings /Smarandache loop rings are commutative or associative or neither for special loop rings satisfying special identities like Moufang, Bol, Bruck or P-identity.

### 6.3 Lie algebras and Smarandache Lie algebras

There are several excellent books on Lie algebras; but the concept of Smarandache Lie algebra is introduced only in this book. At the outset we have a lot to develop in case of Smarandache Lie algebras. As nothing has been done only a very few notions have been introduced in this book. Second, if one views a Lie algebra as just a non-associative ring satisfying some special identities then why several of the properties studied about associative algebras are totally absent in the study of Lie algebras? So in this book the author has tried to introduce some of the properties found in associative rings to these algebras.

Further these algebras became little difficult for analysis as $x^2 = 0$ for all elements x in a Lie algebra. Yet we have introduced several concepts in Lie algebras found in general associative algebras and has left for the reader for further development in this subject. In all places the Smarandache analogous wherever possible have been recalled. It will be an ingenious work for any one who takes up solely the study of Smarandache Lie algebras as there is lot of research to be carried out.

We have introduced Smarandache Lie concepts only for basic Lie algebras so it remains open for the reader to explore further studies and define new concepts of Smarandache notions in case of Lie superalgebras graded Lie superalgebras. Further for the study of complex semi simple Lie algebras refer Springer Monographs in mathematics 2001 in [47]. The reader has a lot of problem to study about Smarandache semi simple complex Lie algebras. Also the study of the very recent textbook for advanced graduate students with all recent Lie algorithms in [19] will certainly help a researcher to built Smarandache Lie algebras in a very systematic and in an innovative way.

As there are several text books on Lie algebras, monographs and collected papers which can give an excellent exposition of Lie algebras we felt it would be sufficient to give a simple exposition about Smarandache Lie algebras and some of its basic properties and thus by leaving a greater scope for the reader to develop nice books on Smarandache Lie algebra.

Also it has become inevitable to mention here that for a solid development of Smarandache Lie algebra one should at least develop four Smarandache concepts viz. Smarandache Homological algebra, Smarandache manifolds, Smarandache topology and Smarandache linear algebras. Finally the dictionary on Lie algebras severs the goal of desktop reference [23]. For any beginner in Lie algebra as well as for any researcher this will serve as a boon Thus the reader has lot for his future study on the concept of Smarandache Lie algebras.

## 6.4 Smarandache Jordan Algebras

The study of Jordan algebras started as early as 1940s. A complete and an excellent exposition of Jordan algebras was given in 1968 by the famous algebraist N. Jacobson in his book entitled "*Structure and representations of Jordan algebras*". Till now a lot of research has been carried out and several others have written lecture notes and books on Jordan algebras. So we have not taken up the Smarandache Jordan algebras in an exhaustive way. But we are happy to state that we have introduced several notions which are common to associative algebras and analysed it in the case of Jordan algebras. We have introduced several Smarandache Jordan concepts but we have not developed them into big theories as we felt that any interested reader can develop them as it can be done as a matter of routine.

Further any researcher has always the flexibility to define Smarandache Jordan algebra in a different ways. We at the outset call any non-associative algebra A to be a Smarandache Jordan algebra if it has a proper subset X which is a subalgebra and satisfies the Jordan identity. We can always build non-trivial Smarandache Jordan algebras by using Smarandache mixed direct product of algebras in which we insist that one of the algebras must be a Jordan algebra so that we have a class of Smarandache Jordan algebras.

It is pertinent to mention here that by the Smarandache mixed direct product of algebras we can simultaneously have a algebra $SA_M$ to be both a Smarandache Jordan algebra and Smarandache Lie algebra and so on i.e. by defining Smarandache mixed direct product of type M as follows:

Let $SA_M = L_1 \times \ldots \times L_n$ be the mixed direct product of non-associative algebras $L_1, \ldots, L_n$ defined over the same field K where at least one of the $L_i$'s is a Jordan algebra and one of the $L_j$'s is a Lie algebra then the Smarandache mixed direct product of algebras of type M is a Smarandache Jordan algebra and Smarandache Lie algebra simultaneously. Except for the Smarandache mixed direct product we will not be in a position to see both structure embedded in a single algebra. If in the S-mixed direct product of type M we take one of the $L_s$'s be a Smarandache loop algebra over the field K then $SA_M$ will also be a Smarandache loop algebra. Thus study of such Smarandache mixed direct products will certainly give several nice and interesting theories.

Study of Smarandache Jordan superalgebras have not been carried out by the author. The author suggests the reader to refer the Lecture notes 211 [18] for the study of Jordan superalgebras so that they can develop Smarandache Jordan superalgebras. One of the certain ways of developing is by including in the Smarandache mixed direct product of type M a $L_i$ to be a Jordan superalgebra i.e. we will call a non-associative algebra A to be a Smarandache Jordan superalgebra if A has a proper subset X such that X is a Jordan superalgebra.

Thus we have the algebra $SA_M$ to be simultaneously a Smarandache Jordan algebra and Smarandache Jordan superalgebra if in the product of $L_i$'s we have a $L_t$ and a $L_m$ distinct to be Jordan algebra and a Jordan superalgebra respectively.

On similar lines by taking in the Smarandache mixed direct product of algebras of type M i.e. in $SA_M$ one of the $L_j$'s to be a Lie superalgebra we can say $SA_M$ is a Smarandache Lie superalgebra. Thus we want to emphasize that these Smarandache mixed direct product of non-associative algebras of type M can give an algebra say $SA_M$ which can be at a time a Smarandache Jordan algebra, Smarandache Lie algebra, Smarandache Lie superalgebra, Smarandache Jordan superalgebra and so on. But otherwise getting a same algebra to have all such properties is impossible.

Thus the reader using these concepts and notions is advised to develop a book of research on Smarandache Jordan algebras.

## 6.5 Other non-associative rings and algebras and their Smarandache analogues

This section mainly introduces the concepts of non-associative algebras like Bernstein, Sagle algebra, finite dimensional complex filiform Leibniz algebra, composition algebra, Akisvis algebra and their Smarandache analogue. Thus the researcher has a lot of scope to develop these ideas. One of the recent study in non-associative algebras is the notion of Bernstein algebras

**DEFINITION [18]**: *A basic algebra A over a field F, char F $\neq$ 2 is a non-associative commutative algebra with a non-zero homomorphism w : A $\rightarrow$ F (weight homomorphism). Among basic algebras, those that satisfy the identity $(x^2)^2 = w(x)^2 x^2$ for all x $\in$ A are called Bernstein algebras.*

*It is well known that Bernstein algebra A always has idempotent elements and if $e^2 = e$ is one of them then w(e) = 1 and ker w = U $\oplus$ V where U = {u $\in$ ker w / 2eu = u} and V = {v $\in$ ker w / ev = 0}.*

*The set of idempotent elements of the Bernstein algebra A given by I(A) = {e + u + $u^2$/ u $\in$ U}.*

We define Smarandache Bernstein algebra as follows.

**DEFINITION 6.5.1**: *A basic algebra A over a field F, char F $\neq$ 2 is a non-associative algebra with a non-zero homomorphism w : A $\rightarrow$ F such that the restriction of w from B $\subset$ A where B $\neq \phi$ (B is a commutative algebra) say $w_R$ : B $\rightarrow$ F such that $(x^2)^2 = w_R(x)^2 x^2$ for all x $\in$ B. Then we call A the Smarandache Bernstein algebra (S-Bernstein algebra).*

Several of the properties enjoyed by the Bernstein algebra can be defined and studied in case of Smarandache Bernstein algebras.

**DEFINITION 6.5.2**: *An algebra satisfying the identities $x^2 = 0$, J(x, y, z)w = (J(w, z, xy) + J(w, y, zx) + J(w, x, yz) is called a Sagle algebra where J(x, y, z) = (xy)z + (zx)y + (yz)x.*

We define Smarandache Sagle algebras as follows:

**DEFINITION 6.5.3**: *Let A be any algebra. We call A a Smarandache Sagle algebra (S-Sagle algebras) if A has a proper subset B, (B $\subset$ A, B $\neq \phi$) such that B is a Sagle algebra.*

**THEOREM 6.5.1**: *If A is a Sagle algebra having proper subalgebras then A is a S-Sagle algebra.*

*Proof*: Direct by the very definition.

It is left as an exercise for the reader to prove the following theorem:

**THEOREM 6.5.2**: *Let L be a S-Sagle algebra then L need not in general be a Sagle algebra.*

For Sagle algebras the reader can develop the concepts analogous to any ring theoretical properties.

Here we describe the finite dimensional complex fillform Leibniz algebras as right multiplication algebras that satisfy abelianness condition. These algebras are generalizations of Lie algebras and are defined by the identity [x[y, z] = [[x, y], z] – [[x, z], y] where [x, y] denotes the algebra product.

If L is a Leibniz algebra the series $C^1(L) = [L, L]$ and $C^{k+1}(L) = [C^k(L), L]$ is defined inductively. The algebra L is called nilpotent if $C^k(L) = 0$ for some k, if a stronger condition $\dim(C^k(L)) = \dim(L) – k – 1$ holds for k = 1, 2, …, dim L – 1, then L is called filiform.

The purpose of this paper is to classify those filiform algebra that satisfy an additional condition namely [L, L] is contained in the centre of L.

**DEFINITION 6.5.4**: *We define Smarandache finite dimensional filiform Leibniz algebras (S-finite dimensional filiform Leibniz algebras) as algebras which has a subalgebra which is finite dimensional filiform Leibniz algebras.*

The easiest way of constructing them is by using S-mixed direct product. If we take in the S-mixed direct product type M i.e. in $SA_M$ an algebra $L_i$ to be a filiform Leibniz algebra then the resulting $SA_M$ will be a Smarandache filiform Leibniz algebra.

A non-associative algebra C over a field F is called a composition algebra if there exists a quadratic form $q : C \to F$ such that $q(xy) = q(x)q(y)$ for all $x, y \in C$ whose induced bilinear form $b_q(x, y) = q(x + y) - q(x) - q(y)$ is non-degenerate. Composition algebras containing a unity element are well-known. Here symmetric composition algebras are investigated. These satisfy additional equation $b_q(xy, z) = b_q(x, yz) \; \forall \; x, y, z \in C$ which forces the dimension of C to be either 1, 4, 2, 8.

**DEFINITION 6.5.5**: *Smarandache composition algebras (S-composition algebras) are algebras A which has a subalgebra $X \subset A$ such that X is a composition algebra.*

Smarandache mixed direct product $SA_M$ easily yields such S-composition algebras.

**DEFINITION [48]**: *An Akivis algebra is a vector space V endowed with a skew symmetric bilinear product $(x, y) \mapsto [x, y] : A \times A \to A$ and a trilinear product $(x, y, z) \to \langle x, y, z \rangle : A \times A \times A \to A$ such that these two operations are related by the so-called Akivis identity which is a generalized form of the Jacobi identity.*

*From any non-associative algebra (B, ., +) one obtains a Akivis algebra A(B) by putting $[x, y] = x . y - y . x$ and $\langle x, y, z \rangle = (x . y) . z - x . (y . z)$.*

We define Smarandache Akivis algebras as follows:

**DEFINITION 6.5.6**: *Let V be a Smarandache vector space. If a subalgebra X of V satisfies the Akivis identity then we call V a Smarandache Akivis algebra (S-Akivis algebra).*

**THEOREM 6.5.3**: *If V is a Akivis algebra which has a S-vector space then V is a S-Akivis algebra.*

*Proof*: Follows directly from the very definitions.

The reader can develop all Smarandache concepts in Akivis algebras and study them as further research on these structures. Thus we see that by introduction of several new Smarandache concepts in these non-associative structures any researcher can always develop not only the existing theories but also new Smarandache theories.

For the study of these structures the available references has been provided but we don't claim that we have exhausted all the references. Only what was available to us has been provided.

Finally these will find its main application in algebraic coding theory, which is a very fast developing area of research. So when Smarandache concepts are introduced it will be much more helpful. Thus this book totally provides immense scope for the reader to develop the subject.

A few simple problems are suggested for the reader.

**PROBLEMS:**

1. Give an example of a Akisvis algebra.
2. Does there exist a Akisvis algebra which is not a S-Akisvis algebra.
3. Give an example of a S-Sagle algebra.
4. Does there exist a finite dimensional Bernstein algebra?

**Chapter 7**

# SUGGESTED PROBLEMS

This chapter suggests 150 problems on non-associative rings and related concepts. Certainly these problems will not only make the researcher get more insight into the subject but also help him to get fixed ideas by tackling these problems. Some of these problems can even be very challenging research problems. Whereas problems of introducing and characterizing Moufang algebras, Bol algebras or Bruck algebras will give us bulk results, as such this type of trial to the best of my knowledge has not been carried out by anyone. Further some of the solutions to these problems can generate substantial material for research. As understanding of any subject lies in the ability to solve related problems, we have suggested the following problems. We wish to state that we have not included all the problems, which we have mentioned in our suggestions for future study.

1. Let R and $R_1$ be two SNA-rings. $\phi : R$ to $R_1$ be a S-homomorphism find conditions on $\phi$ so that Ker $\phi$ is an ideal of R. (Hint: We don't assume Ker $\phi$ = (0), we assume Ker $\phi \neq$ (0)).

2. Find the SNA-ring of least order.

3. Is the smallest non-associative ring a SNA-ring?

4. Can on $Z_n$ be defined two binary operations so that $Z_n$ is a non-associative ring?

5. Find the smallest SNA-ring, which is a SNA-Bol ring.

6. Does there exist SNA-rings other than the ones got from

    a. loop rings?
    b. Groupoid rings?

7. Find a SNA-ring in which every ideal is a SNA-ideal of R.

8. Find conditions on a SNA-ring R so that every subring of R is a SNA-subring of R.

    a. Characterize SNA-rings which has ideals none of which are SNA-ideals
    b. Characterize SNA-rings which has subrings none of which are SNA-rings.

9. Characterize those SNA-rings which are SNA-simple.

10. Let R be SNA-ring having a maximal SNA-ideal say I relative to a SNA-subring B of R. Find conditions on R and B so that R/I is simple.

11. For the loop ring $Z_2 L_n(m)$ n > 3, n odd and m ≤ n, (m, n) = 1 = (m, n-1) (for a fixed n)

    a. How many S-right ideal exists in $Z_2L_n(m)$?
    b. What is the lattice structure of S-right ideals?
    c. How many right ideals exists in $Z_2 L_n(m)$?
    d. What is the lattice structure of right ideals?
    e. Compare the lattice structures in (b) and (d).

12. Study the same problem for

    a. left ideals.
    b. SNA-left ideals.
    c. ideals and
    d. SNA-ideals.

13. Let $Z_t L_n(m)$ be the loop ring. Study problems (11) and (12) when

    a. (t, n + 1) = 1
    b. (t, n + 1) = t; t a prime
    c. (t, n) = 1
    d. (t, n) = t
    e. (t, n) = d
    f. (t, n + 1) = d

14. Characterize those units in $ZL_n(m)$; $L_n(m) \in L_n$ which are S-NA units.

15. Find all units in $Z_p L_n(m)$; $L_n(m) \in L_n$

    a. Which are not SNA-units
    b. Which are SNA-units.

16. Find conditions on n and p so that the loop ring $Z_p L_n(m)$ has S-idempotents where (n, p) = 1 or (n, p) = d.

17. Find necessary and sufficient conditions for the loop ring $Z_n L_n(m)$ where $L_n(m) \in L_n$ to have nontrivial S-quasi regular elements.

18. Does $Z_2L_n(m)$ have S-quasi regular elements x where |supp x| is an even number?

19. Find $S(J(Z_2L_n(m)))$.

20. Is $SJ(Z_pL_n(m)) \neq \phi$, if (n, p) = 1? Justify your answer.

21. Find $SJ(ZL_n(m))$.

22. Characterize those NA rings in which every quasi regular element is a S-quasi regular element.

23. Can the loops in the class $L_n$ for any suitable ring or field contribute quasi regular elements which are also S-quasi regular elements?

24. Characterize those loop rings or NA rings in which $J(R) = SJ(R)$.

25. Find a necessary and sufficient condition for any NA-ring to be both strictly right ring and S-strictly right ring.

26. Find condition for the loop ring $Z_n L_n(m)$ to be

    a. strictly right ring.
    b. S-strictly right ring.

27. Classify those NA-rings which satisfy only strictly right ring property and does not satisfy the S-strictly right ring property.

    Note: Study the same problem just in case of ideals.

28. Give an example of a S-ring which is not a S-strictly right ring.

29. Determine all S-semi-idempotents of $ZL_n(m)$; $L_n(m) \in L_n$.

30. Find conditions on n and t so that the collection of S-semi-idempotents in $Z_t L_n(m)$ is non-empty.

31. When will $SI(R) = SSI(R)$? Find conditions on the NA ring R.

32. Can $Z_p L_n(m)$ when $L_n(m)$ is a non-commutative loop be a S-Marot ring? Justify your answer. Characterize those loops in $L_n$ which are not commutative still they happen to be

    a. S-weak Marot ring.
    b. S-Marot ring.

33. Let L be any loop, K any field and KL the loop ring of the loop L over K. Find all idempotents in $L^*$, can KL have idempotents which are Smarandache idempotents? Find a necessary and sufficient condition on loop ring KL, so that

    a. $L^*$ is a loop
    b. $L^*$ is only a groupoid

    Can $L^*$ every be a semigroup? Justify your claim.

34. Characterize those loops L so that $L^* = S(L^*)$.

35. Let $L_n(m) \in L_n$ and K be any field. Find conditions on $KL_n(m)$ so that

a. $S((L_n(m))^*) \neq 0$
   b. Find the number of elements in $S((L_n(m))^*)$ or equivalently find the number of S-subloops in $KL_n(m)$.
   c. When will $S(L_n(m))^*)$ be a loop?
   d. When will $S((L_n(m))^*)$ be only a groupoid?

36. In a loop L, is the concept of t.u.p loops and u.p loop equivalent or distinct?

37. Let L be a loop which is not a u.p loop. Can L be a S.u.p loop? Find conditions on L so that L is a S-strong u.p loop but is not a u.p loop.

38. Find those loops in $L_n(m) \in L_n$ which are

   a. S-u.p loops.
   b. S-S.u.p loops.

39. Give an example of a Smarandache S.u.p loop which is not a u.p loop.

40. Can we say if L is a u.p loop or a t.u.p loop, L can be S.u.p loop, or S.S. u.p loop or S.t.u.p loop or S.S.t.u.p loop? (<u>Hint</u>: Illustrate by examples. This question is vital as we may have S-subloops A of L. It may happen L is a u.p loop still $A \subset L$ will not be a u.p loop).

41. If L is a S.S.u.p loop or a S.S.t.u.p loop will L be a u.p loop or a t.u.p loop?

42. Which of the loops in $L_n(m)$ are S-Hamiltonian loops?

43. Characterize those loops of odd order which are S-Hamiltonian.

44. Let $Z_n$ be a ring and L be a loop of order m;

   a. If $(n, m) = 1$ when does the loop ring $Z_nL$ contain nontrivial orthogonal ideals?
   b. If $(m, n) = d$, can we say $Z_nL$ contains always a pair of nontrivial orthogonal ideals?

45. Obtain a necessary and sufficient condition on the loop ring FL so that every pair of ideals in FL is orthogonal.

46. Characterize those loop rings in which every pair of orthogonal ideals is also a S-orthogonal ideal.

47. Characterize those loop rings which are

   a. Strictly right commutative.
   b. Strongly right commutative.
   c. Right commutative.

48. Let $Z_pL_n(m)$ be the loop ring $L_n(m) \in L_n$. Can the loop ring $Z_pL_n$ be

a. S-strongly right commutative
    b. S-right commutative.

49. Characterize those loop rings which are inner commutative.

50. Which of the loop rings $ZL_n(m)$ are inner commutative? Characterize those loop rings which are inner commutative but not S-inner commutative.

51. Characterize those loop rings which are inner commutative and S-inner commutative.

52. Obtain a necessary and sufficient condition for the loop ring RL to have only one S(P(RL)). (<u>Hint</u>: It is assumed that RL has more than one S-subring).

53. Characterize those loop rings RL which are such that each S-subring A of L give a distinct S(P(RL)).

54. Characterize those loop rings RL in which P(RL) = S(P(RL)).

55. Characterize those loop rings in which all ideals and S-ideals coincide.

56. Characterize those loop rings in which the S-ideals form a

    a. Super modular lattice.
    b. Strongly modular lattice.

57. Let $L_n(m) \in L_n$ and Z the ring of integers Find the lattice structure of

    a. ideals of $ZL_n(m)$.
    b. S-ideals of $ZL_n(m)$.
    c. S-subrings of $ZL_n(m)$ and
    d. Subrings of $ZL_n(m)$.

    Study problem (57) in which Z is replaced by $Z_p$.

58. Find conditions on the loop ring $ZL_n(m)$ or $Z_p L_n(m)$ so that the

    a. $SN (ZL_n(m)) = N (ZL_n(m))$.
    b. $SN (Z_p L_n(m)) = N (Z_p L_n(m))$.
    c. $Z (ZL_n(m)) = S (Z (ZL_n(m)))$.
    d. $Z (Z_p L_n(m)) = SZ (Z_p L_n(m))$.

59. Does the center and nucleus of $RL_n(m)$ coincide?

60. Does the S-center and S-nucleus coincide on $RL_n(m)$?

61. Find all units and S-units in the loop ring $ZL_n(m)$. Do we have loop ring so that both the sets are identical?

62. Find conditions on the loop ring RL so that S-units in RL are central units.

63. Characterize those vector spaces V which has S-basis which are S-strong basis of V.

64. Characterize those vector spaces V for which every basis is a S-basis. Also characterize those vector spaces which has no S-basis.

65. Characterize those vector spaces which has S-basis.

66. Can we have any vector space V in which every basis is a S-basis?

67. Study Zassenhaus conjecture for the class of loop rings $ZL_n(m)$.

68. Verify S-Zassenhaus conjecture for the new class of loop rings $ZL_n(m)$.

69. Find the S-normalized units of $ZL_n(m)$.

70. Find S U $(ZL_n(m))$.

71. Find $V(ZL_n(m))$.

72. Find all groupoid rings RG (necessary and sufficient condition) so that every idempotent in RG is a S-idempotent of RG. Another extreme case is no idempotent in RG is a S-idempotent of RG.

73. Find whether the groupoid rings ZG where G is a groupoid built using modulo integers be a

    a. n-ideal ring.
    b. S-n-ideal ring.

74. Study the above question when Z is replaced by $Z_n$, where $Z_n$ is:

    a. Any prime field $Z_p$.
    b. Any ring $Z_n$.

75. Characterize or prove or disprove the existence of groupoid rings in which the notion of S-n-ideal ring and n-ideal ring coincides.

76. Let $Z_nG$ be a groupoid ring. Find conditions on G and $Z_n$ so that the groupoid ring $Z_nG$ has every generalized semi-ideal to be a S-generalized semi-ideal of $Z_nG$.

77. Find necessary and sufficient condition for the groupoid ring RG to have no generalized semi-ideal to be a S-generalized semi-ideal of RG.

78. Find the number of S-mod p envelope $SG^*$ for the groupoids in $Z(n)$, $Z^*(n)$, $Z^{**}(n)$ and $Z^{***}(n)$, where n is a fixed number.

79. Prove or disprove the number of elements in $SG^*$ depends also on the field $Z_p$ we choose over which the groupoid rings are defined.

80. Does there exist a groupoid G which has only one S-subgroupoid such that $SG^* = G^*$? Justify your claim.

81. Classify those groupoid rings using the groupoids from $Z(n)$, $Z^*(n)$, $Z^{**}(n)$ and $Z^{***}(n)$ in which for suitable rings

    a. The minimal and maximal ideal coincide.
    b. The S- maximal ideals
    c. The S- minimal ideals.
    d. Can S-maximal ideals and S-minimal ideals be coincident in ZG where $G \in Z^{***}(n)$?

82. Does the groupoid ring $Z_nG$ where $G \in Z^{***}(n)$ have the set of all ideals to form a modular lattice?

83. Study the problem 82 in case of

    a. S-ideals
    b. S-pseudo ideals

84. Find necessary and sufficient conditions so that the groupoid ring KG satisfies

    a. A.C.C
    b. D.C.C.

85. Find a necessary and sufficient conditions for the groupoid rings RG to satisfy

    a. S-A.C.C
    b. S-D.C.C

    on S-ideals of RG.

86. Characterize those groupoid rings which are S-E-rings.

87. Obtain a necessary and sufficient condition for the groupoid ring to be a S-p-ring.

88. Find a class of groupoid rings which are S-pre-J-rings.

89. Show that the class of groupoid rings which are S-reduced is nonempty.

90. Determine a class of S-quasi commutative groupoid rings.

91. Characterize those groupoid rings which are Lin rings.

92. Characterize those groupoid rings RG in which every ideal is a S-ideal. Hence or otherwise, prove that RG is both a ideally obedient ring and S-ideally obedient ring.

93. Characterize the class of groupoid rings which are S-I∗- rings.

94. Study problem (93) for the groupoid ring RG where $G \in Z^{***}(n)$.

95. Find necessary and sufficient condition for the groupoid rings $Z_mG$ where $G \in Z^{***}(n)$ to be S-quotient-ring (S-Q ring) with $Z_sG$ where $G \in Z^{***}(r)$.

96. Characterize those groupoid rings RG (where $G \in Z^{***}(n)$) so that

    a. RG is never a F-ring.
    b. RG is always a F-ring.

97. Find conditions on G and on the ring R so that the groupoid ring RG is a $\gamma_n$-ring.

98. For the groupoid $L = \{ [-\infty, \infty], *, (m, n)\}$ find associator and commutator.

99. Using the groupoid L given in the above problem. Study the groupoid ring $Z_pL$, ZL and any RL.

100. Find a necessary and sufficient condition for a groupoid ring RG to be a S-normal groupoid ring.

101. Does the new class of groupoids $L = \{[-\infty, \infty], *, (m, n)\}$ and their groupoid rings enjoy any new property if

    a. $(m, n) = 1$.
    b. if $n = m - 1$ or $m = (n - 1)$.
    c. if m and n are both primes.

102. Obtain a necessary and sufficient condition on the groupoid ring RG to satisfy the identity $(xy)x = x(yx)$ for all $x, y \in RG$.

103. Is it possible to find S-exponentiation ring using groupoids without unit? Or if the groupoids G are replaced by S-groupoid will we have RG to be a S-exponentiation ring? Justify your answer.

104. Classify those groupoid rings which has n-capacitor subgroups.

105. Characterize those commutative groupoid rings which do not have n-capacitor subgroups.

106. Characterize those groupoid rings RG which are

    a. essential.
    b. S- essential.

    for $R = Z_n$ and $G \in Z^{**}(n)$.

107. Find the class of groupoid rings $Z_nG$ where $G \in Z^{***}(n)$ which are

   a. not essential
   b. not S-essential

108. Characterize those groupoid rings $Z_pG$ which has

   a. non trivial radical.
   b. non trivial S-radical.

   where $G \in Z^{***}(n)$ for

   a. $(n, p) = 1$.
   b. $(n, p) = p$.
   c. $n = p$.

109. Find a necessary and sufficient condition on the groupoid rings $Z_mG$ where m is a prime or a composite number and $G \in Z^{***}(n)$ to be S-Lie ring.

110. Find necessary and sufficient condition on the loop ring $Z_kL_n(m)$ to a S-Lie ring where K may be a prime or a composite number and $L_n(m) \in L_n$.

111. The problems (109) and (110) can be reformulated as follows: Let $Z_mG$ or $Z_kL_n^{(m)}$ be a groupoid ring of a loop ring. Find a necessary and sufficient condition for the subring U =
$\{\alpha = \Sigma \alpha_i g_i \mid g_i \in G \text{ or } g_i \in L_n(m) \text{ such that } \Sigma \alpha_i = 0\}$ to be

   a. Smarandache NA ring
   b. Satisfy the identities $x^2 = 0$ and $(xy)z + (yz)x + (zx)y = 0$ i.e. U satisfies Lie identities. Is $Z_mG$ or $Z_kL_n(m)$ is a S-Lie algebra?

112. Characterize those groupoid rings $Z_nG$ where $G \in Z^{***}(n)$ such that $U = \{\alpha = \Sigma \alpha_i g_i \in Z_mG / \Sigma \alpha_i = 0\}$ is a SNA ring.

113. Study the same problem 112 in the context of loop rings.

114. Characterize those S-Lie algebras which have S-Lie maximal ideals which are not maximal ideals in general.

115. Does there exist examples of S-Lie minimal ideals which are not minimal ideals. Justify your answer.

116. Describe those S-Lie algebras

   a. which has no S-pseudo Lie ideals
   b. which has S-pseudo Lie ideals.

117. Is the S-pseudo quotient ring L/X a S-Lie ring?

118. What is the structure of L/X described in problem 117? Is it at least a Lie ring? Justify your answers.

119. Characterize those S-Lie rings which are S-filial Lie rings.

120. Let L be a Lie ring. Find necessary and sufficient conditions for L to be a S-filial Lie ring.

121. Characterize those Lie rings which are

    a. S-Lie subsemi-ideal ring.
    b. Subsemi ideal ring.

122. Characterize those Lie rings which are

    a. S-semiprime.
    b. Semiprime.

123. Characterize those Lie rings which are

    a. S-Chinese Lie rings.
    b. Chinese Lie rings.

124. Characterize those Lie rings which are never Chinese Lie rings or S-Chinese Lie rings.

125. Characterize those Lie algebras which are

    a. Lie multiplication ideals
    b. S-Lie multiplication ideals

126. Characterize those Lie algebras which has no Lie multiplication ideals.

127. Study those Lie algebras which has S-Lie multiplication ideals but has no Lie multiplication ideals.

128. Characterize those Lie algebras which satisfy

    a. Super ore condition.
    b. S-super ore condition.

129. Describe those Lie algebras which are

    a. Obedient Lie rings.
    b. S-obedient Lie rings.

130. Characterize those Lie algebras in which every subalgebra of same dimension are

    a. Conjugate.
    b. S-conjugate.

131. Characterize those Lie algebras which are generalized stable algebras.

132. Classify those Lie algebras which are S-generalized stable Lie algebras.

133. Obtain a necessary and sufficient condition for Lie algebras to be

    a. n-ring ideal.
    b. S-n-ring ideal.

134. Derive a necessary and sufficient condition for a Lie algebra to be Smarandache semi-simple.

135. Characterize those Lie algebras which are S-simple.

136. Obtain some interesting results on

    a. S-Lie simple algebras.
    b. S-Lie semi-simple algebras.

137. Characterize those Lie algebras which are S-strongly Lie simple.

138. Characterize those Lie algebras which are S-strongly Lie semi-simple.

139. Does there exist any relations or marked differences between S-strongly Lie simple algebras and S-strongly Lie semi-simple algebra? Justify your answer with explicit examples.

140. Can we prove any S-ideal is a S-semi-simple Lie algebra of characteristic 0 is S-semi-simple.

141. Classify those Lie algebras which has S-Cartan subalgebras.

142. Characterize those Lie algebras which are both S-Lie simple and S-Lie semi-simple.

143. Obtain a necessary and sufficient condition so that every S-Jordan algebra has a S-magnifying element.

144. Characterize those S-Jordan algebras which are S-shrinkable.

145. Find a necessary and sufficient condition for a S-Jordan algebra to be S-dispotent Jordan algebra.

146. Obtain a necessary and sufficient condition for a S-Jordan ring to be

    a. S-weakly G-Jordan ring.
    b. S-G-Jordan rings.

147. Characterize those S-Jordan rings which has at least one S-clean elements.

148. Can there be a S-Jordan ring in which every element can be a S-clean element?

149. Characterize those S-Jordan rings which has no S-clean elements.

150. Characterize those Akisvis algebras which are S-Akisvis algebras.

# INDEX

# About the Author

Dr. W. B. Vasantha is an Associate Professor in the Department of Mathematics, Indian Institute of Technology Madras, Chennai, where she lives with her husband Dr. K. Kandasamy and daughters Meena and Kama. Her current interests include Smarandache algebraic structures, fuzzy theory, coding/ communication theory. In the past decade she has guided seven Ph.D. scholars in the different fields of non-associative algebras, algebraic coding theory, transportation theory, fuzzy groups, and applications of fuzzy theory to the problems faced in chemical industries and cement industries. Currently, six Ph.D. scholars are working under her guidance. She has to her credit 241 research papers of which 200 are individually authored. Apart from this she and her students have presented around 262 papers in national and international conferences. She teaches both undergraduate and post-graduate students at IIT and has guided 41 M.Sc. and M.Tech projects. She has worked in collaboration projects with the Indian Space Research Organization and with the Tamil Nadu State AIDS Control Society.

She can be contacted at vasantha@iitm.ac.in
You can visit her on the web at http://mat.iitm.ac.in/~wbv